\input ./style/arxiv-general.cfg
\documentclass[aap,MSNbibl,dvips]{arximspdf}
\makeatletter
   \@ifpackageloaded{graphicx}{}{\usepackage{graphicx}}
\makeatother

\doi{10.1214/15-AAP1111}
\volume{26}
\issue{2}
\pubyear{2016}
\firstpage{1029}
\lastpage{1081}
\docsubty{FLA}

\makeatletter
\newcommand{\rrvert}{\vert}
\newcommand{\rrVert}{\Vert}
\newcommand{\llvert}{\vert}
\newcommand{\llVert}{\Vert}
\newtheorem{lemma}{Lemma}[section]
\newproclaim{defn}[lemma]{Definition}

\newtheorem{theorem}[lemma]{Theorem}
\newtheorem{proposition}[lemma]{Proposition}
\newtheorem{corollary}[lemma]{Corollary}
\newproclaim{rem}[lemma]{Remark}
\newproclaim{notn}[lemma]{Notation}
\newproclaim{example}[lemma]{Example}
\makeatother

\begin{document}
\begin{frontmatter}

\title{A consistency estimate for Kac's model of elastic collisions in
a dilute gas}
\runtitle{A consistency estimate for Kac's model}

\begin{aug}
\author[A]{\fnms{James}~\snm{Norris}\corref{}\thanksref{T1}\ead[label=e1]{J.R.Norris@statslab.cam.ac.uk}}
\runauthor{J. Norris}
\affiliation{University of Cambridge}
\address[A]{Statistical Laboratory\\
Centre for Mathematical Sciences\\
Wilberforce Road\\
Cambridge CB3 0WB\\
United Kingdom\\
\printead{e1}}
\end{aug}
\thankstext{T1}{Supported by EPSRC Grant EP/103372X/1.}

%
\received{\smonth{5} \syear{2014}}
%
\revised{\smonth{3} \syear{2015}}

%
\begin{abstract}
An explicit estimate is derived for Kac's mean-field model of colliding
hard spheres,
which compares, in a Wasserstein distance, the empirical velocity
distributions for two versions of the
model based on different numbers of particles.
For suitable initial data, with high probability, the two processes
agree to within a tolerance of order $N^{-1/d}$,
where $N$ is the smaller particle number and $d$ is the dimension,
provided that $d\ge3$.
From this estimate we can deduce that the spatially homogeneous
Boltzmann equation is well posed in
a class of measure-valued processes and provides a good approximation
to the Kac process when the number of particles is large.
We also prove in an appendix a basic lemma on the total variation of
time-integrals of time-dependent signed measures.
\end{abstract}

%
\end{frontmatter}

\section{Kac process}\label{KAC}
\setcounter{footnote}{1}

Kac \cite{MR0084985} proposed in 1954 a random process to model the
dynamics of a dilute gas.
The process models the velocities of $N$ particles in $\mathbb{R}^d$
as they
evolve under elastic collisions.
The case $d=3$ is of main interest, but we will allow any $d\ge2$.
Since no account is taken of particle positions, any physical
justification for the model relies on assumptions of spatial
homogeneity and rapid mixing.
It is thus impossible to give a physical meaning to the number of
particles $N$.
Yet, on the mathematical side, we have to make a choice.
Hence it is of interest to show consistency for sufficiently large
values of $N$.

Kac's process depends on a choice of collision kernel $B$.
This is a finite measurable kernel $B(v,d\sigma)$ on $\mathbb
{R}^d\times S^{d-1}$
which is chosen to model physical characteristics of the gas.
The collision kernel specifies the rate for collisions of pairs of
particles with incoming relative velocity $v$ and outgoing direction of
separation $\sigma$.
Since collisions are assumed to conserve momentum and energy, for a
pair of particles with pre-collision velocities $v$ and $v_*$, and
hence relative velocity $v-v_*$,
the post-collision velocities $v'=v'(v,v_*,\sigma)$ and
$v'_*=v'_*(v,v_*,\sigma)$ are determined by the direction of
separation through
\[
v'+v'_*=v+v_*, \qquad  v'-v'_*=
\llvert v-v_*\rrvert \sigma.
\]
We will often write $u$ for the direction of approach, given by
$u=(v-v_*)/|v-v_*|$.
We assume throughout that, for all $u\in S^{d-1}$, $B(u,\cdot)$ is a
probability measure, supported on $S^{d-1}\setminus\{-u,u\}$,
and that the following standard scaling and symmetry properties hold.
For ${\lambda}\in[0,\infty)$ and $u\in S^{d-1}$, and for any
isometry $R$ of $S^{d-1}$, we have
\begin{equation}
\label{LBP}
B({\lambda}u,\cdot)={\lambda}B(u,\cdot),\qquad  B(Ru,\cdot)=B(u,\cdot)\circ
R^{-1}.
\end{equation}
Our main results require further that the map $u\mapsto B(u,\cdot)$ is
Lipschitz on $S^{d-1}$ for the total variation norm on measures on $S^{d-1}$.
Then there is a constant $\kappa\in[1,\infty)$ such that, for all
$v,v'\in\mathbb{R}^d$,
\begin{equation}
\label{LCP}
\bigl\llVert B(v,\cdot)-B\bigl(v',\cdot\bigr)\bigr\rrVert \le
\kappa\bigl\llvert v-v'\bigr\rrvert.
\end{equation}
Here and throughout, we denote the total variation norm by $\|\cdot\|$.
The Boltzmann sphere $\mathcal{S}$ is the set of probability measures
$\mu$
on $\mathbb{R}^d$ such that\footnote{Here, on the left-hand side, and where convenient below, we use $v$ to
denote the identity function on $\mathbb{R}^d$.}
\[
\langle v,\mu\rangle=\int_{\mathbb{R}^d}v\mu(dv)=0,\qquad \bigl\langle
\llvert v\rrvert ^2,\mu\bigr\rangle=\int_{\mathbb{R}
^d}
\llvert v\rrvert ^2\mu(dv)=1.
\]
For $N\in\mathbb{N}$, write $\mathcal{S}_N$ for the subset of
$\mathcal{S}$ of normalized
empirical measures of the form $N^{-1}\sum_{i=1}^N\delta_{v_i}$.
The Kac process with collision kernel $B$ and particle number $N$ is
the Markov chain in $\mathcal{S}_N$ with generator $\mathcal{G}$
given on bounded
measurable functions $F$ by
\[
\mathcal{G}F(\mu)=N\int_{\mathbb{R}^d\times\mathbb{R}^d\times
S^{d-1}}\bigl\{F\bigl(
\mu^{v,v_*,\sigma
}\bigr)-F(\mu)\bigr\}\mu(dv)\mu(dv_*)B(v-v_*,d \sigma),
\]
where
\[
\mu^{v,v_*,\sigma}=\mu+N^{-1}\{\delta_{v'}+
\delta_{v'_*}-\delta _v-\delta_{v_*}\}.
\]
The choice of state-space $\mathcal{S}_N$ is possible because in each
collision the number of particles, the momentum $v+v_*$ and the energy
$|v|^2+|v_*|^2$ are conserved.
There is no Kac process on $\mathcal{S}_1$ because this set is empty.
For $N\ge2$, the transition rates of the Kac process are bounded by
$2N$ on $\mathcal{S}_N$.
Hence, by the elementary theory of Markov chains, given any initial
state $\mu_0^N\in\mathcal{S}_N$, there exists a Kac process $(\mu
_t^N)_{t\ge0}$ in $\mathcal{S}_N$ starting from $\mu^N_0$,
the law of this process is unique, and almost surely it takes only
finitely many values in any compact time interval.

It is of special interest to model particles colliding as hard spheres.
Under plausible physical assumptions, this leads, by a well-known
calculation, to the choice of kernel $B(v,d\sigma)\propto|v|\sin
^{3-d}(\theta/2)\,d \sigma$,
where $\theta\in[0,\pi]$ is given by $\cos\theta=u\cdot\sigma$
and $d\sigma$
is the uniform distribution on $S^{d-1}$.
It is straightforward to check that (\ref{LBP}) and (\ref{LCP}) hold
in this case for all $d\ge2$.
Indeed, for $d=3$, we can take $\kappa=1$, and the dynamics have a
particularly simple description:
for every pair of particles $(v,v_*)$, at rate $|v-v_*|/N$, consider
the sphere with poles at $v$ and $v_*$; choose randomly a new axis for
the sphere,
label the poles $v'$ and $v'_*$ and replace $v$ and $v_*$ by $v'$ and~$v_*'$.

Consider the set $\mathcal{F}$ of functions $f$ on $\mathbb{R}^d$
such that
\[
\bigl\llvert \hat f(v)\bigr\rrvert \le1,  \qquad \bigl\llvert \hat f(v)-\hat f
\bigl(v'\bigr)\bigr\rrvert \le\bigl\llvert v-v'\bigr\rrvert
\]
for all $v,v'$, where\footnote{The notation is chosen as a reminder of
the shape of the weight function $1/(1+|v|^2)$.}
\[
\hat{f}(v)=f(v)/\bigl(1+|v|^2\bigr).
\]
Define a distance function $W$ on $\mathcal{S}$ by
\[
W(\mu,\nu)=\sup_{f\in\mathcal{F}}\langle f,\mu-\nu\rangle.
\]
Then $W$ makes $\mathcal{S}$ into a complete separable metric space.
This is shown in Section~\ref{W}, along with the convergence of a
natural approximation scheme by random samples in $(\mathcal{S},W)$.
Our first main result is the following consistency estimate for Kac
processes with different numbers of particles.
We make no assumption on the joint law of the processes.
They could, for example, be independent.

\begin{theorem}\label{MR}
Assume that the collision kernel $B$ satisfies conditions (\ref{LBP}) and (\ref{LCP}).
Let $\varepsilon\in(0,1]$, ${\lambda}\in[1,\infty)$, $p\in
(2,\infty)$
and $T\in[0,\infty)$.
Then there exist constants $\alpha(d,p)>0$ and $C(B,d,\varepsilon,{\lambda
},p,T)<\infty$ with the following property.
Let $N,N'\in\mathbb{N}$ with $N\le N'$, and let $(\mu^N_t)_{t\ge0}$ and
$(\mu_t^{N'})_{t\ge0}$ be Kac processes in $\mathcal{S}_N$ and
$\mathcal{S}_{N'}$
such that
\begin{equation}
\label{UMC}
\bigl\langle|v|^p,\mu_0^N\bigr
\rangle\le{\lambda},\qquad \bigl\langle|v|^p,\mu _0^{N'}
\bigr\rangle\le {\lambda}.
\end{equation}
Then, with probability exceeding $1-\varepsilon$, for all $t\in
[0,T]$, we have
\[
W\bigl(\mu^N_t,\mu_t^{N'}\bigr)\le
C\bigl(W\bigl(\mu_0^N,\mu_0^{N'}
\bigr)+N^{-\alpha}\bigr).
\]
\end{theorem}

We have not found a way to prove a similar estimate for $p=2$.
This is consistent with the current theory for the spatially
homogeneous Boltzmann equation where also, for $p=2$, there is no
quantitative stability estimate.
We can improve the rate of convergence at the cost of a stronger moment
condition.

\begin{theorem}\label{MR10}
Assume further that $(\mu^N_t)_{t\ge0}$ and $(\mu_t^{N'})_{t\ge0}$
are adapted as Markov processes to a common filtration.
For $p>8$ and $d\ge3$, we can take $\alpha=1/d$ in Theorem~\ref{MR}.
Also, for $p>8$ and $d=2$, we can replace $N^{-\alpha}$ in Theorem~\ref
{MR} by $N^{-1/2}\log N$.
\end{theorem}

The theorems could be considered as providing a measure of accuracy for
a Monte Carlo scheme, using say $N$ computational particles, for the
evolution of a Kac process having a much larger number of particles $N'$.

The rates of convergence in Theorem~\ref{MR10} are known to be optimal
for the convergence of sample empirical distributions in Wasserstein distance.
Indeed, there is no discrete approximation\vspace*{1pt} scheme for a smooth measure
which achieves a rate better than $N^{-1/d}$.
So it seems unlikely that the rates in can be improved in this context.
Our need for the condition $p>8$ can be traced to the stochastic
convolution estimates in Section~\ref{SCB}.
We show in Section~\ref{W} that, for laws in $\mathcal{S}$ having a finite
$p$th moment, their sample empirical distributions converge in the
metric $W$ with optimal rates if $p>3d/(d-1)$, but this can fail if
$p<3d/(d-1)$.
This makes it plausible that some moment condition beyond $p>2$ is
necessary for the conclusions of Theorem~\ref{MR10}, but we do not
know whether this is so.

By combining Theorem~\ref{MR10} with Proposition~\ref{ME} below, we
obtain the following estimate.
\begin{theorem}
For $d\ge3$, for all $\varepsilon\in(0,1]$ and all $\tau,T\in
(0,\infty)$
with $\tau\le T$, there is a constant $C(B,d,\varepsilon,\tau
,T)<\infty$ such that,\vspace*{1pt}
for all $N,N'\in\mathbb{N}$ with $N\le N'$ and any Kac processes
$(\mu
_t^N)_{t\ge0}$ in $\mathcal{S}_N$ and $(\mu_t^{N'})_{t\ge0}$ in
$\mathcal{S}_{N'}$,
with probability exceeding $1-\varepsilon$, for all $t\in[\tau,T]$,
we have
\[
W\bigl(\mu^N_t,\mu_t^{N'}\bigr)\le
C\bigl(W\bigl(\mu_\tau^N,\mu_\tau^{N'}
\bigr)+N^{-1/d}\bigr).
\]
\end{theorem}
Note that $\tau$ can be arbitrarily small, and we obtain here the
optimal rate $N^{-1/d}$ without the supplementary moment condition
(\ref{UMC}).
Thus it is only for the initial evolution of the processes that
consistency may rely on a such a moment condition.

We have avoided so far any mention of the Boltzmann equation, which
classically is the starting point for kinetic theory.
We shall show in our other main results, Theorem~\ref{CBE} and
Corollaries \ref{COR} and \ref{COR2},
that the consistency estimate leads quickly to existence and uniqueness
of measure solutions for the spatially homogeneous Boltzmann equation,
and convergence to such solutions of the Kac process in the large $N$ limit.
Indeed, we obtain a more precise estimate of this convergence than was
previously known.
This was the original motivation for our work.

In the next two sections, we identify martingales of the Kac process,
and we derive some moment estimates.
The difference of two Kac processes, with the same collision kernel but
different numbers of particles, satisfies a noisy version of a
linearized Boltzmann equation.
In Section~\ref{BPR} we develop a representation formula for solutions
of this equation in terms
of an auxiliary branching process, which we call the linearized Kac process.
We use coupling arguments for this process to develop some estimates.
The proof of Theorem~\ref{MR} is given in Section~\ref{PMR}.

We develop in Section~\ref{MKP} some further continuity estimates for
the linearized Kac process,
and in Section~\ref{SCB} some maximal inequalities for stochastic
convolutions appearing in the representation formula.
These are then used in Section~\ref{MRD} to prove Theorem~\ref{MR10}.
The relation of our estimates to prior work on the Kac process and the
spatially homogeneous Boltzmann equation is discussed in Section~\ref{SHBE}.
The final section is a self-contained appendix, proving a basic result
on the evolution of signed measures, which is used in Sections~\ref{BPR} and \ref{SHBE}.

\section{Martingales of the Kac process}\label{MARTK}
We compute the martingale decomposition for linear functions of the Kac
process $(\mu^N_t)_{t\ge0}$.
Set $E=\mathbb{R}^d\times\mathbb{R}^d\times S^{d-1}\times(0,\infty)$.
Denote by $m$ the un-normalized empirical measure on $E$ of the set
of all random vectors $(V,V_*,\Sigma,T)$ such that there is a
collision at time $T$ in the particle system $(\mu^N_t)_{t\ge0}$
of a velocity pair $(V,V_*)$ with direction of separation $\Sigma$.
Denote by $\bar m$ the random measure on $E$ given by
\[
\bar m(dv,dv_*,d\sigma,dt)=N\mu^N_{t-}(dv)\mu
^N_{t-}(dv_*)B(v-v_*,d\sigma)\,dt.
\]
Define a random signed measure $M^N$ on $(0,\infty)\times\mathbb
{R}^d$ by specifying,
for bounded measurable functions $f$ on $(0,\infty)\times\mathbb
{R}^d$, the
integral\footnote{We will sometimes write $f_s$ for $f(s,\cdot)$.}
\begin{eqnarray}
M^{N,f}_t &=&\int_0^t
\bigl\langle f_s,dM^N_s\bigr\rangle =\int
_{(0,t]\times\mathbb{R}^d}f(s,v)M^N(ds,dv)
\nonumber
\\[-2pt]
\label{KLL}
 &=&\frac{1}N\int_E\bigl\{f_s\bigl(v'\bigr)+f_s
\bigl(v'_*\bigr)-f_s(v)-f_s(v_*)\bigr\}
\\[-2pt]
\nonumber
&&\hspace*{3pt}\qquad{}\times 1_{(0,t]}(s) (m-\bar m) (dv,dv_*,d\sigma,ds).
\end{eqnarray}
Then, by standard results for Markov chains, the process
$(M^{N,f}_t)_{t\ge0}$ is a martingale.
We use the same notation also in the case where $f$ has no dependence
on the time parameter.
Define for finite measures $\mu,\nu$ on $\mathbb{R}^d$ a signed measure
$Q(\mu,\nu)$ on $\mathbb{R}^d$ by specifying, for
bounded measurable functions $f$ of compact support in $\mathbb{R}^d$,
the\vspace*{-2pt} integral
\begin{eqnarray}
 \bigl\langle f,Q(\mu,\nu)\bigr\rangle
&=&\int_{\mathbb{R}^d\times\mathbb
{R}^d\times S^{d-1}}
\bigl\{ f\bigl(v'\bigr)+f\bigl(v'_*\bigr)-f(v)-f(v_*)
\bigr\}
\nonumber
\\[-9pt]
\label{TBO}
\\[-9pt]
\nonumber
&&\hspace*{57pt}{}\times\mu(dv)\nu(dv_*)B(v-v_*,d\sigma).
\end{eqnarray}
Then the martingale decomposition for $(\langle f,\mu^N_t\rangle
)_{t\ge0}$ is
given by
\begin{equation}
\label{MD} \bigl\langle f,\mu^N_t\bigr\rangle=\bigl
\langle f,\mu_0^N\bigr\rangle+M^{N,f}_t+
\int_0^t\bigl\langle f,Q\bigl(
\mu^N_s,\mu _s^N\bigr)\bigr
\rangle \,ds.
\end{equation}

We note for later use the following estimates. First, by Doob's
$L^2$-inequality,\footnote{For the same calculation in a general
setting, see, for example, \cite{MR2395153}, Proposition~8.7.}
\begin{eqnarray}
&&\mathbb{E} \Bigl(\sup_{s\le t}\bigl\llvert M_s^{N,f}
\bigr\rrvert ^2 \Bigr)\nonumber\\
&&\label{DME}\qquad \le  \frac{4}{N^2}\mathbb{E}\int
_E\bigl\{f_s\bigl(v'
\bigr)+f_s\bigl(v'_*\bigr)-f_s(v)-f_s(v_*)
\bigr\} ^21_{(0,t]}(s)\bar m(dv,dv_*,d\sigma,ds)
\\
\nonumber
&&\qquad\le  128\llVert f\rrVert _\infty t/N.
\end{eqnarray}
Next, for the total variation measure $|M^N|$ of $M^N$, we have
\begin{eqnarray}
&&\mathbb{E}\int_{(0,t]\times\mathbb{R}^d}\bigl(1+|v|^2
\bigr)\bigl\llvert M^N(ds,dv)\bigr\rrvert\nonumber \\
&&\qquad\le  \mathbb{E}\int
_E\bigl(4+2\llvert v\rrvert ^2+2\llvert v_*
\rrvert ^2\bigr)1_{(0,t]}(s) (m+\bar m) (dv,dv_*,d\sigma ,ds)
\nonumber
\\[-8pt]
\label{MDV}
\\[-8pt]
\nonumber
&&\qquad= \mathbb{E}\int_0^t\!\!\int_{\mathbb{R}^d\times\mathbb
{R}^d}
\bigl(8+4|v|^2+4\llvert v_*\rrvert ^2\bigr)\llvert v-v_*
\rrvert \mu _s^N(dv)\mu_s^N(dv_*)\,ds
\\
&&\qquad\le  24\mathbb{E}\int_0^t\bigl
\langle1+|v|^3,\mu_s^N\bigr\rangle \,ds.
\nonumber
\end{eqnarray}
We\vspace*{1pt} used $|v-v_*|\le|v|+|v_*|$ and the fact that $\mu_s^N\in\mathcal
{S}$ for
the second inequality.
Finally, for any interval $(s,s']$ during which $(\mu_t^N)_{t\ge0}$
does not jump, there is no contribution to the left-hand side of (\ref
{MDV}) from $m$,
so the same calculation yields the following pathwise estimate:
\begin{equation}
\label{TVMST}
\int_{(s,s']\times\mathbb{R}^d}\bigl(1+|v|^2\bigr)\bigl
\llvert M^N(dr,dv)\bigr\rrvert \le12\int_s^{s'}
\bigl\langle1+|v|^3,\mu_r^N\bigr\rangle \,dr.
\end{equation}

\section{Moment estimates for the Kac process}\label{MOM}
We derive some moment inequalities for the Kac process, which we shall
use later.
The basic arguments are standard for the Boltzmann equation and are
applied to the Kac process in \cite{MR3069113}, Lemma~5.4.
We have quantified the moment-improving property and added some maximal
inequalities.
We begin with the Povzner inequality.
For all $p\in(2,\infty)$, there is a constant $\beta(B,p)>0$ such that,
for all $v,v_*\in\mathbb{R}^d$ and for $u=(v-v_*)/|v-v_*|$,
\begin{eqnarray}
&& \int_{S^{d-1}}\bigl\{\bigl\llvert v'
\bigr\rrvert ^p+\bigl\llvert v_*'\bigr\rrvert
^p-|v|^p-\llvert v_*\rrvert ^p\bigr\}B(u,d
\sigma)
\nonumber
\\[-8pt]
\label{POV}
\\[-8pt]
\nonumber
&&\qquad\le-\beta \bigl(\llvert v\rrvert ^p+\llvert v_*\rrvert
^p\bigr)+\beta^{-1}\bigl(\llvert v\rrvert \llvert v_*
\rrvert ^{p-1}+\llvert v\rrvert ^{p-1}\llvert v_*\rrvert \bigr).
\end{eqnarray}
Here is a proof for the class of collision kernels we consider.
Note first that
\begin{eqnarray}
\bigl\llvert v'\bigr\rrvert ^p+\bigl
\llvert v_*'\bigr\rrvert ^p &\le & \bigl(\llvert v\rrvert
^2+\llvert v_*\rrvert ^2\bigr)^{p/2}
\nonumber
\\[-8pt]
\label{CPP}
\\[-8pt]
\nonumber
&\le & \llvert v
\rrvert ^p+\llvert v_*\rrvert ^p+C(p) \bigl(\llvert v
\rrvert \llvert v_*\rrvert ^{p-1}+|v|^{p-1}\llvert v_*\rrvert
\bigr).
\end{eqnarray}
It suffices by symmetry to consider the case $|v_*|\le|v|$.
Set $y=|v-v_*|(u+\sigma)/2$, then $v'=v_*+y$ and $|y|^2=|v-v_*|^2t$, where
$t=(1+u\cdot\sigma)/2$.
Note that $t\in(0,1)$ for $B(u,\cdot)$-almost all $\sigma$.
We use the inequalities $|v'|\le|y|+|v_*|$ and $|v-v_*|\le|v|+|v_*|$
to see that, for all $\delta\in(0,1]$,
\[
\bigl\llvert v'\bigr\rrvert ^2\le(1+
\delta)|y|^2+\bigl(1+\delta^{-1}\bigr)\llvert v_*\rrvert
^2\le(1+\delta )^2t\llvert v\rrvert ^2+2
\bigl(1+\delta ^{-1}\bigr)\llvert v_*\rrvert ^2.
\]
From this inequality and a similar one for $|v'_*|^2$, we deduce that,
for some $C(p)<\infty$,
\begin{eqnarray*}
\bigl\llvert v'\bigr\rrvert ^p &\le & (1+
\delta)^{p+1}t^{p/2}|v|^p+C(p)\delta^{-1}
\llvert v_*\rrvert ^p,\\
 \bigl\llvert v'_*\bigr\rrvert
^p &\le & (1+\delta)^{p+1}(1-t)^{p/2}\llvert v\rrvert
^p+C(p)\delta^{-1}\llvert v_*\rrvert ^p.
\end{eqnarray*}
Then
\begin{eqnarray}
&& \bigl\llvert v'\bigr\rrvert ^p+\bigl
\llvert v_*'\bigr\rrvert ^p-\llvert v\rrvert
^p-\llvert v_*\rrvert ^p
\nonumber
\\[-8pt]
\label{PPOV}
\\[-8pt]
\nonumber
&&\qquad\le-\beta(\delta ,t) \bigl(
\llvert v\rrvert ^p+\llvert v_*\rrvert ^p\bigr)+C(p)\delta
^{-1}\bigl(\llvert v\rrvert \llvert v_*\rrvert ^{p-1}+\llvert
v\rrvert ^{p-1}\llvert v_*\rrvert \bigr),
\end{eqnarray}
where $\beta(\delta,t)=(1-(1+\delta)^{p+1}(t^{p/2}+(1-t)^{p/2}))^+/2$.
Set $\beta(\delta)=(\delta/C(p))\wedge\int_{S^{d-1}}\beta(\delta
,t)B(u,d\sigma)$.
Then\vspace*{1.5pt} we obtain (\ref{POV}) for $u$ with $\beta=\beta(\delta)$ by
integrating~(\ref{PPOV}).
But $\beta(\delta)$ does not depend on $u$ by the symmetry condition
(\ref{LBP}) and $\beta(\delta)>0$ for all sufficiently small $\delta
$, so we are done.

\begin{proposition}\label{ME}
Let $(\mu^N_t)_{t\ge0}$ be a Kac process with collision kernel $B$
satisfying (\ref{LBP}).
Let $p\in[2,\infty)$ and $q\in(2,\infty)$ with $p\le q$.
There exists a constant $C(B,p,q)<\infty$ such that, for all $t\ge0$,
we have
\begin{equation}
\label{NMG} \mathbb{E} \bigl(\bigl\langle\llvert v\rrvert ^q,
\mu_t^N\bigr\rangle \bigr)\le C\bigl(1+t^{p-q}
\bigr)\bigl\langle\llvert v\rrvert ^p,\mu_0^N
\bigr\rangle.
\end{equation}
Moreover, there is a constant $C(B,q)<\infty$ such that, for all $t\ge0$,
\begin{equation}
\label{NME}
\mathbb{E} \Bigl(\sup_{s\le t}\bigl\langle\llvert
v\rrvert ^q,\mu_s^N\bigr\rangle \Bigr)
\le(1+Ct)\bigl\langle\llvert v\rrvert ^q,\mu _0^N
\bigr\rangle,
\end{equation}
and there is a constant $C(B,p,q)<\infty$ such that, for all $t\ge0$,
\begin{equation}
\label{NMF} \mathbb{E} \Bigl(\sup_{s\le t}\bigl\langle1+
\llvert v\rrvert ^p,\bigl\llvert \mu_s^N-\mu
_0^N\bigr\rrvert \bigr\rangle \Bigr)\le C
\bigl(t+t^{q-p}\bigr)\bigl\langle\llvert v\rrvert ^q,
\mu_0^N\bigr\rangle.
\end{equation}
%
\end{proposition}
\begin{pf}
By the Povzner inequality, there are constants $\beta(B,q)>0$ and
$C(B,q)<\infty$ such that, for all $v,v_*\in\mathbb{R}^d$,
\begin{eqnarray*}
&&\int_{S^{d-1}}\bigl\{\bigl\llvert v'\bigr\rrvert
^q+\bigl\llvert v'_*\bigr\rrvert ^q-\llvert
v\rrvert ^q-\llvert v_*\rrvert ^q\bigr\}B(v-v_*,d\sigma)
\\
&&\qquad\le-\beta\llvert v-v_*\rrvert \bigl(\llvert v\rrvert ^q+\llvert
v_*\rrvert ^q\bigr)+\beta ^{-1}\llvert v-v_*\rrvert \bigl(
\llvert v\rrvert \llvert v_*\rrvert ^{q-1}+\llvert v\rrvert
^{q-1}\llvert v_*\rrvert \bigr)
\\
&&\qquad \le-\beta \bigl(\llvert v\rrvert ^{q+1}+\llvert v_*\rrvert
^{q+1}\bigr)+C\bigl(\llvert v\rrvert ^q\bigl(1+\llvert v_*
\rrvert \bigr)+\bigl(1+\llvert v\rrvert \bigr)\llvert v_*\rrvert ^q
\bigr).
\end{eqnarray*}
Set $f_q(t)=\mathbb{E}(\langle|v|^q,\mu_t^N\rangle)$ and
$f_{q,p}(t)=f_q(t)/f_p^*$,
where $f_p^*=\sup_{t\ge0}f_p(t)$.
Since $\langle|v|^q,\mu\rangle\le N^{q/2}$ for all $\mu\in\mathcal
{S}$,\vspace*{1pt} we have
$f_q(t)\le N^{q/2}<\infty$ for all $t$.
The process $(\langle|v|^q,\mu_t^N\rangle)_{t\ge0}$ makes\vspace*{1pt} jumps of
size $\{
|v'|^q+|v'_*|^q-|v|^q-|v_*|^q\}/N$ at rate $N\mu_{t-}^N(dv)\mu
_{t-}^N(dv_*)B(v-v_*,d\sigma)\,dt$.
Hence
\begin{eqnarray*}
\nonumber
f_q(t) &=& f_q(0)\\
&&\!{}+\mathbb{E}\int
_0^t\!\!\int\bigl\{\bigl\llvert v'
\bigr\rrvert ^q+\bigl\llvert v'_*\bigr\rrvert
^q-\llvert v\rrvert ^q-\llvert v_*\rrvert ^q
\bigr\}\mu _s^N(dv)\mu_s^N(dv_*)B(v-v_*,d
\sigma)\,ds
\\
&\le&  f_q(0)-2\beta\int_0^tf_{q+1}(s)\,ds+2C
\int_0^tf_q(s)\,ds.\label{FPPP}
\end{eqnarray*}
By H\"older's inequality, we have $f_q(t)^{q-p+1}\le
f_{q+1}(t)^{q-p}f_p(t)$, so we deduce that
\[
f_{q,p}(t)\le f_{q,p}(0)-2\beta\int_0^t
\bigl(f_{q,p}(s)\bigr)^{1+1/(q-p)}\,ds+2C\int_0^tf_{q,p}(s)\,ds
\]
which implies by standard arguments that, for some $C(B,p,q)<\infty$
and all $t\ge0$, we have
\begin{equation}
\label{FPPL}
f_{q,p}(t)\le C\bigl(1+f_{q,p}(0)\wedge
t^{p-q}\bigr).
\end{equation}
Now $f^*_2=1$, so by taking $p=2$, we obtain (\ref{NMG}), for the
cases $p=2$ and $p=q$.
In particular, this shows that $f_p^*\le C\langle|v|^p,\mu_0^N\rangle
$ for all
$p$, so (\ref{FPPL}) implies (\ref{NMG}) also for $p\in(2,q)$.

Consider\vspace*{1.5pt} the process $(A_t)_{t\ge0}$ starting from $0$ which jumps by
$\{|v||v_*|^{q-1}+|v|^{q-1}|v_*|\}/N$ when $(\langle|v|^q,\mu
_t^N\rangle)_{t\ge
0}$ jumps by $\{|v'|^q+|v'_*|^q-|v|^q-|v_*|^q\}/N$.
Then
\begin{eqnarray*}
\label{FPSS}
\mathbb{E}(A_t) &=& \mathbb{E}\int_0^t\!\!
\int_{\mathbb{R}^d\times
\mathbb{R}^d} \bigl\{ \llvert v\rrvert \llvert v_*\rrvert
^{q-1}+\llvert v\rrvert ^{q-1}\llvert v_*\rrvert \bigr\}\llvert
v-v_*\rrvert \mu_s^N(dv)\mu ^N_s(dv_*)\,ds
\\
&\le& 4\int_0^tf_q(s)\,ds.
\end{eqnarray*}
Now $\sup_{s\le t}\langle|v|^q,\mu_s^N\rangle\le\langle|v|^q,\mu
_0^N\rangle+C(q)A_t$
for all $t$, where $C(q)$ is the constant from (\ref{CPP}).
Hence we obtain (\ref{NME}) by taking expectations and using the case
$p=q$ of (\ref{NMG}) to estimate $f_q$.

The\vspace*{1pt} process $(\langle1+|v|^p,|\mu_t^N-\mu_0^N|\rangle)_{t\ge0}$
jumps by at
most $\{4+|v'|^p+|v'_*|^p+|v|^p+|v_*|^p\}/N$ at each jump of $(\mu
_t^N)_{t\ge0}$.
Consider the process $(B_t)_{t\ge0}$ starting from $0$ which jumps by
$\{1+|v|^p+|v_*|^p\}/N$ at the same times.
Then $\langle1+|v|^p,|\mu_s^N-\mu_0^N|\rangle\le2^pB_t$ whenever
$s\le t$ and
\begin{eqnarray*}
\label{FPSR} \mathbb{E}(B_t) &=& \mathbb{E}\int_0^t\!\!
\int_{\mathbb{R}^d\times\mathbb{R}^d} \bigl\{ 1+\llvert v\rrvert ^p+\llvert
v_*\rrvert ^p \bigr\} \llvert v-v_*\rrvert \mu_s^N(dv)
\mu^N_s(dv_*)\,ds \\
&\le & 6\int_0^tf_{p+1}(s)\,ds.
\end{eqnarray*}
So (\ref{NMF}) follows from (\ref{NMG}).
%
\end{pf}

\section{Linearized Kac process and representation formula}\label{BPR}
In this section we introduce a branching process of signed particles in
$\mathbb{R}^d$ which may be considered as a linearization of the Kac process.
A particular case of this process allows us to write a representation
formula for the difference of two Kac processes $(\mu^N_t)_{t\ge0}$
and $(\mu^{N'}_t)_{t\ge0}$.
We use coupling arguments to obtain continuity estimates for the
branching process, which are later used to control $\mu_t^N-\mu_t^{N'}$.
The representation formula rests only on the fact that $(\mu_t^N-\mu
_t^{N'})_{t\ge0}$ solves the linear equation (\ref{LEB}) below.
It seems possible that the same conclusions can be reached by a direct
analysis of this equation, but we have not done this.

The branching process will have ``positive'' and ``negative''
particles, making the following general notation convenient.
Given a set $V$, we denote by $V^*$ the signed space $V\times\{-1,1\}
=V^-\cup V^+$, by $\pi$ the projection $V^*\to V$ and by $\pi_\pm$
the bijections $V^\pm\to V$.
Note that $*$ does not signify the dual space.
From now on, we set $V=\mathbb{R}^d$.

The data for our branching process are an initial time $s\in[0,\infty
)$ and an initial type $v\in V^*$,
together with a process $(\rho_t)_{t\ge0}$ of measures on $\mathbb{R}^d$
such that, for all~$t$,
\begin{equation}
\label{SKP} \langle1,\rho_t\rangle\le1, \qquad  \bigl\langle|v|^2,
\rho_t\bigr\rangle\le1.
\end{equation}
The case $\rho_t=(\mu^N_t+\mu_t^{N'})/2$ will be of main interest later.
Consider the continuous-time branching particle system\footnote{%
The dynamics of the branching process can be motivated as follows.
Fix a large integer $N$, and suppose that $(N\rho_t)_{t\ge0}$ evolves
as an unnormalized Kac process on $N$ particles.
Consider the perturbed process obtained by introducing one additional
particle of velocity $v$ at time $s$,
where the pairwise collision rules are unchanged and where transitions
are coupled as far as possible with the original.
The discrepancy between the original and the perturbed systems will
grow over time approximately as the branching process $(\Lambda
^*_t)_{t\ge s}$,
a ``negative'' particle in $V^-$ corresponding to one present in the
original system but removed by collision in the perturbed system.
Formally, the approximation becomes exact as $N\to\infty$. We do not
rely on this.
The construction of $(\Lambda^*_t)_{t\ge s}$ does not require $(\rho
_t)_{t\ge0}$ to be a Kac process.}
with types in $V^*$ where each particle of type $v$ in $V^\pm$, at
rate $2\rho_t(dv_*)B(v-v_*,d\sigma)\,dt$ for $v_*\in\mathbb{R}^d$ and
$\sigma\in S^{d-1}$,
dies and is replaced by three particles $v'(v,v_*,\sigma)$ and
$v'_*(v,v_*,\sigma)$ in $V^\pm$ and $v_*$ in $V^\mp$.
More properly, the rate is $2\rho_t(dv_*)B(\pi(v)-v_*,d\sigma)\,dt$ and
the offspring are
$(v'(\pi(v),v_*,\sigma),1)$, $(v'_*(\pi(v),v_*,\sigma),1)$ and $(v_*,-1)$
when $v\in V^+$, and
$(v'(\pi(v),v_*,\sigma),-1)$, $(v'_*(\pi(v),v_*,\sigma),-1)$ and $(v_*,1)$
when $v\in V^-$.
We assume throughout that, for all $t\ge0$,
\begin{equation}
\label{BTM}
\int_0^t\bigl\langle\llvert v
\rrvert ^3,\rho_s\bigr\rangle \,ds<\infty.
\end{equation}
We will show that (\ref{BTM}) ensures there is no explosion; that is,
the time $T_n$ of the $n$th branching event tends to $\infty$ almost surely.
So the process is well defined for all time by the specification of its
branching rates, and consists at all times $t\ge s$ of a finite number
of particles.
Write $(\Lambda^*_t)_{t\ge s}$ for the associated process of un-normalized
empirical measures on $V^*$.
We call this process the \textit{linearized Kac process in environment
$(\rho_t)_{t\ge0}$ starting from $v$ at time $s$}.

Set $\Lambda_t=\Lambda_t^*\circ\pi^{-1}$.
Then $(\Lambda_t)_{t\ge s}$ is itself the empirical process of a branching
process in $V$, in which we forget the book-keeping exercise of giving
a sign to each particle.
Write $E_{(s,v)}$ for the expectation over $(\Lambda^*_t)_{t\ge s}$ to
recall that $\Lambda^*_s=\delta_v$
and that this is not the full expectation in the case that $(\rho
_t)_{t\ge0}$ is itself random.
Given an initial type $v\in\mathbb{R}^d$, without a sign, we will by default
start the process $(\Lambda^*_t)_{t\ge s}$ with the positive type $(v,1)$.

\begin{proposition}\label{JPE}
There is almost surely no explosion in the branching construction
described above.
Moreover, for all $p\in[2,\infty)$, there is a constant $c(p)<\infty$
such that, for all $v_0\in\mathbb{R}^d$ and all $t\ge s$, we have
\[
E_{(s,v_0)}\bigl\langle1+\llvert v\rrvert ^p,
\Lambda_t\bigr\rangle\le\bigl(1+\llvert v_0\rrvert
^p\bigr)\exp \biggl\{c(p)\int_s^t
\bigl\langle 1+\llvert v\rrvert ^{p+1},\rho_r\bigr\rangle
\,dr \biggr\}.
\]
In particular we can take $c(2)=8$.
\end{proposition}
We will reserve the notation $c(p)$ for this constant throughout.
We will also use throughout the notation
\[
\tilde\Lambda_t=\Lambda_t^+-\Lambda_t^-,\qquad
\Lambda_t^\pm =\Lambda^*_t\circ
\pi_\pm^{-1}.
\]
Thus $\Lambda_t^+$ and $\Lambda_t^-$ are random measures on $\mathbb
{R}^d$, which are
the empirical distributions of positive and negative particles,
and $\tilde\Lambda_t$ is a random signed measure on $\mathbb{R}^d$.
Note that $\Lambda_t=\Lambda_t^++\Lambda_t^-$.
By Proposition~\ref{JPE}, we can define, for any $s,t\ge0$ with $s\le
t$, a linear map $E_{st}$ on the set of measurable functions of
quadratic growth on $\mathbb{R}^d$ by
\[
E_{st}f(v)=E_{(s,v)}\langle f,\tilde\Lambda_t
\rangle.
\]
Note that, by the Markov property, we have $E_{st}E_{tu}=E_{su}$.
We will write $f_{st}$ for $E_{st}f$ and sometimes just $f_s$ when the
value of $t$ is understood.
We will use the same notation for functions $f$ of polynomial growth,
whenever $(\rho_t)_{t\ge0}$ has sufficient moments for this to make
sense using Proposition~\ref{JPE}.
We base our main argument on the following representation formula,
which is proved at the end of this section.

\begin{proposition}\label{FMU}
In the case where $\rho_t=(\mu_t^N+\mu_t^{N'})/2$ for all $t$, we have
\[
\bigl\langle f,\mu_t^N-\mu_t^{N'}
\bigr\rangle=\bigl\langle f_{0t},\mu_0^N-\mu
_0^{N'}\bigr\rangle+\int_0^t
\bigl\langle f_{st},dM_s^N\bigr\rangle-\int
_0^t\bigl\langle f_{st},dM^{N'}_s
\bigr\rangle.
\]
\end{proposition}

We will use the following two estimates expressing continuity of the
linearized Kac process in its initial data.
Write $\|f\|$ for the smallest constant such that $|\hat f(v)|\le\|f\|
$ and $|\hat f(v)-\hat f(v')|\le\|f\||v-v'|$ for all $v,v'\in\mathbb{R}^d$.
Thus $f\in\mathcal{F}$ if and only if $\|f\|\le1$.

\begin{proposition}\label{FSE}
Assume condition (\ref{LCP}).
Then
\[
\llVert E_{st}f\rrVert \le3\bigl(1+6\kappa(t-s)\bigr)\exp \biggl\{
\int_s^t8\bigl\langle 1+\llvert v\rrvert
^3,\rho_r\bigr\rangle \,dr \biggr\}\llVert f\rrVert.
\]
\end{proposition}

\begin{proposition}\label{FSF}
For all $v\in\mathbb{R}^d$ and all $s,s'\in[0,t]$ with $s\le s'$, we have
\begin{eqnarray*}
&& \bigl\llvert E_{st}f(v)-E_{s't}f(v)\bigr\rrvert \\
&&\qquad \le5
\bigl(1+\llvert v\rrvert ^3\bigr)\exp \biggl\{\int
_s^t8\bigl\langle 1+\llvert v\rrvert
^3,\rho_r\bigr\rangle \,dr \biggr\}\llVert f\rrVert \int
_s^{s'}\bigl\langle1+\llvert v\rrvert
^3,\rho _r\bigr\rangle \,dr.
\end{eqnarray*}
\end{proposition}

\begin{pf*}{Proof of Proposition~\ref{JPE}}
Consider first the case $p=2$ and $s=0$.
Fix $v_0\in\mathbb{R}^d$, and consider the branching particle system
$(\Lambda
_t)_{t<{\zeta}}$ starting from $\delta_{v_0}$ at time $0$ and run up to
explosion ${\zeta}=\sup_nT_n$.
Note that, at a branching event with colliding particle velocity $v_*$,
the total number of particles in the system increases by $2$, and the
total kinetic energy increases by $|v'|^2+|v'_*|^2+|v_*|^2-|v|^2=2|v_*|^2$.
Hence $\langle1+|v|^2,\Lambda_t\rangle$ makes jumps of size
$2(1+|v_*|^2)$ at rate
$2|v-v_*|\Lambda_{t-}(dv)\rho_t(dv_*)\,dt$.
Set $S_n=\inf\{t<{\zeta}\dvtx \langle1+|v|^2,\Lambda_t\rangle\ge n\}$,
and set
\[
g(t)=E_{(0,v_0)}\bigl\langle1+\llvert v\rrvert ^2,
\Lambda_{t\wedge S_n}\bigr\rangle.
\]
Note that $S_n\le T_n$.
We use the estimate
\[
\bigl(1+\llvert v_*\rrvert ^2\bigr)\llvert v-v_*\rrvert \le2\bigl(1+
\llvert v\rrvert ^2\bigr) \bigl(1+\llvert v_*\rrvert ^3
\bigr)
\]
to see that
\[
\int_{\mathbb{R}^d\times\mathbb{R}^d}\bigl(1+\llvert v_*\rrvert ^2\bigr)
\llvert v-v_*\rrvert \Lambda _t(dv)\rho_t(dv_*)\le
2m_3(t)\bigl\langle1+\llvert v\rrvert ^2,
\Lambda_t\bigr\rangle,
\]
where $m_3(t)=\langle1+|v|^3,\rho_t\rangle$.
Hence, by optional stopping, the process
\[
\bigl\langle1+\llvert v\rrvert ^2,\Lambda_{t\wedge S_n}\bigr
\rangle-\int_0^{t\wedge
S_n}8m_3(s)\bigl
\langle1+\llvert v\rrvert ^2,\Lambda _s\bigr\rangle \,ds
\]
is a supermartingale.
On taking expectations, we obtain
\begin{eqnarray*}
g(t) &\le & 1+\llvert v_0\rrvert ^2+E_{(0,v_0)}\int
_0^{t\wedge S_n}8m_3(s)\bigl\langle 1+\llvert v
\rrvert ^2,\Lambda _s\bigr\rangle \,ds\\
&\le & 1+\llvert
v_0\rrvert ^2+\int_0^t8m_3(s)g(s)\,ds
\end{eqnarray*}
so $g(t)<\infty$ and then
\[
g(t)\le\bigl(1+\llvert v_0\rrvert ^2\bigr)\exp \biggl\{
\int_0^t8m_3(s)\,ds \biggr\}.
\]
The right-hand side does not depend on $n$, so we must have $S_n\to
\infty$ almost surely as $n\to\infty$.
Hence $T_n\to\infty$ almost surely, and the claimed estimate follows
by monotone convergence.

For $p\in(2,\infty)$, there is a constant $C(p)<\infty$ such that
\begin{eqnarray}
\bigl\llvert v'\bigr\rrvert ^p+\bigl
\llvert v'_*\bigr\rrvert ^p+\llvert v_*\rrvert
^p-\llvert v\rrvert ^p &\le & \bigl(\llvert v\rrvert
^2+\llvert v_*\rrvert ^2\bigr)^{p/2}+\llvert v_*
\rrvert ^p-\llvert v\rrvert ^p
\nonumber
\\[-8pt]
\label{POF}
\\[-8pt]
\nonumber
 &\le &  C(p) \bigl(\llvert v
\rrvert ^{p-2}\llvert v_*\rrvert ^2+\llvert v_*\rrvert
^p\bigr)
\end{eqnarray}
and then, for another constant $c(p)<\infty$,
\begin{equation}
\label{POE}
\hspace*{8pt}2\llvert v-v_*\rrvert \bigl(2+\bigl\llvert v'\bigr
\rrvert ^p+\bigl\llvert v'_*\bigr\rrvert ^p+
\llvert v_*\rrvert ^p-\llvert v\rrvert ^p\bigr)\le c(p)
\bigl(1+\llvert v\rrvert ^p\bigr) \bigl(1+\llvert v_*\rrvert
^{p+1}\bigr).\hspace*{-20pt}
\end{equation}
The argument used for $p=2$ then gives the desired estimate in the case $s=0$.
The argument is the same for $s\ge0$.
\end{pf*}

We now describe a coupling of linearized Kac processes starting from
different initial velocities, constructed to
branch at the same times and with the same sampled velocities $v_*$ and
angles $\sigma$, as far as possible.
To simplify, we begin without the signs.
Define sets
\begin{equation}
\label{TS} V_0=\mathbb{R}^d\times
\mathbb{R}^d,\qquad V_1=\mathbb{R}^d\times\{ 1\},
\qquad V_2=\mathbb{R} ^d\times\{2\}
\end{equation}
which we treat as disjoint.
Consider the continuous-time branching process in $V_0\cup V_1\cup V_2$
with the following branching mechanism.
For each particle (of type) $(v_1,v_2)\in V_0$, there are three
possible transitions.
First, at rate $2B(v_1-v_*,d\sigma)\wedge B(v_2-v_*,d\sigma)\rho_t(dv_*)\,dt$
for $v_*\in\mathbb{R}^d$ and $\sigma\in S^{d-1}$,
the particle $(v_1,v_2)$ dies and is replaced by three particles
$(v_*,v_*)$, $(v'_1,v'_2)$ and $(v'_{1*},v'_{2*})$ in $V_0$.
Here we are writing $v'_k$ for $v'(v_k,v_*,\sigma)$ and $v'_{k*}$ for
$v'_*(v_k,v_*,\sigma)$ for short.
Call this a coupled transition.
Second, at rate $2(B(v_1-v_*,d\sigma)-B(v_2-v_*,d\sigma))^+\rho_t(dv_*)\,dt$,
the particle $(v_1,v_2)$ dies and is replaced by four particles $v_*$,
$v'_1$ and $v'_{1*}$ in $V_1$ and $v_2$ in $V_2$.
Third, at rate $2(B(v_1-v_*,d\sigma)-B(v_2-v_*,d\sigma))^-\rho_t(dv_*)\,dt$,
the particle $(v_1,v_2)$ dies and is replaced by $v_*$, $v'_2$ and
$v'_{2*}$ in $V_2$ and $v_1$ in $V_1$.\vspace*{1pt}
The second and third will be called decoupling transitions.
Finally, for $k=1,2$, each particle $v_k$ in $V_k$, at rate
$2B(v_k-v_*,d\sigma)\rho_t(dv_*)\,dt$,
dies and is replaced by three particles, $v_*$, $v'_k$ and $v'_{k*}$ in $V_k$.
It is easy to check, by the triangle inequality, that in each coupled
transition, we have $|v'_1-v'_2|\le|v_1-v_2|$ and
$|v'_{1*}-v'_{2*}|\le|v_1-v_2|$.

Fix $v_1,v_2\in\mathbb{R}^d$, and suppose we start with one particle
$(v_1,v_2)\in V_0$ at time $0$.
Write $(\Gamma^0_t,\Gamma^1_t,\Gamma^2_t)_{t\ge0}$ for the
empirical\vspace*{1.5pt} process of
particle types on $ V_0\cup V_1\cup V_2$.
Then, inductively, $\Gamma^0_t$ is supported on pairs $(u_1,u_2)$
with $|u_1-u_2|\le|v_1-v_2|$.
For $k=1,2$, write $p_k$ for the projection\vspace*{1pt} to the $k$th component
$V_0\to\mathbb{R}^d$, and write $\pi_k$ for the bijection $V_k\to
\mathbb{R}^d$.
Define a measure $\Lambda_t^k$ on $\mathbb{R}^d$ by
\[
\Lambda_t^k=\Gamma_t^0\circ
p_k^{-1}+\Gamma_t^k\circ
\pi_k^{-1}.
\]
It is straightforward to check that $(\Lambda_t^1)_{t\ge0}$ and
$(\Lambda
_t^2)_{t\ge0}$ are copies of the Markov
process $(\Lambda_t)_{t\ge0}$ starting from $\delta_{v_1}$ and
$\delta_{v_2}$,
respectively.

For $k=0,1,2$, consider the signed space $V^*_k=V_k^-\cup
V_k^+=V_k\times\{-1,1\}$.
The process $(\Gamma^0_t,\Gamma_t^1,\Gamma_t^2)_{t\ge0}$ lifts in
an obvious
way to a branching process $(\Gamma^{0,*}_t,\Gamma_t^{1,*},\Gamma
_t^{2,*})_{t\ge0}$
in $V^*_0\cup V_1^*\cup V_2^*$ starting from $((v_1,v_2),1)$ in
$V^+_0$, where the ``$v_*$'' offspring switch signs, just as in
$(\Lambda
^*_t)_{t\ge0}$.
By lift we mean that $\Gamma_t=\Gamma_t^*\circ\pi^{-1}$ for the projection
$\pi \dvtx  V_k^*\to V_k$.
We write $E_{(0,v_1,v_2)}$ for the expectation over this process.
For $k=1,2$, set
\[
\Lambda_t^{k,*}=\Gamma_t^{0,*}\circ
p_k^{-1}+\Gamma_t^{k,*}\circ
\pi_k^{-1}.
\]
Then $(\Lambda_t^{k,*})_{t\ge0}$ is a linearized Kac process with
environment $(\rho_t)_{t\ge0}$ starting from $(v_k,1)$.

\begin{lemma}\label{GHU}
Assume condition (\ref{LCP}).
Then
\begin{eqnarray*}
&& E_{(0,v_1,v_2)}\bigl\langle1+\llvert v\rrvert ^2,
\Gamma^1_t+\Gamma_t^2\bigr\rangle
\\
&&\qquad \le 6\kappa t\bigl(2+\llvert v_1\rrvert ^2+\llvert
v_2\rrvert ^2\bigr)\llvert v_1-v_2
\rrvert \exp \biggl\{\int_0^t8\bigl\langle 1+
\llvert v\rrvert ^3,\rho_s\bigr\rangle \,ds \biggr\}.
\end{eqnarray*}
Moreover, for all $p\in(2,\infty)$, there is a constant $C(p)<\infty
$ such that
\begin{eqnarray*}
&& E_{(0,v_1,v_2)}\bigl\langle1+\llvert v\rrvert ^p,
\Gamma^1_t+\Gamma_t^2\bigr\rangle
\\
&&\qquad \le C(p)\kappa t\bigl(1+\llvert v_1\rrvert ^p+\llvert
v_2\rrvert ^p\bigr)\llvert v_1-v_2
\rrvert \exp \biggl\{\int_0^tc(p)\bigl\langle 1+
\llvert v\rrvert ^{p+1},\rho_s\bigr\rangle \,ds \biggr\}.
\end{eqnarray*}
\end{lemma}
\begin{pf}
The decoupling transition $(u_1,u_2)\to(u'_1,u'_{1*},v_*;u_2)$ occurs
at rate
\[
2\Gamma^0_{t-}\bigl(d(u_1,u_2)
\bigr) \bigl(B(u_1-v_*,d\sigma)-B(u_2-v_*,d\sigma )
\bigr)^+\rho_t(dv_*)\,dt
\]
and increases $\langle1+|v|^2,\Gamma_t^1+\Gamma_t^2\rangle$ by
$4+2|v_*|^2+|u_1|^2+|u_2|^2$.
On adding the rate for the other decoupling transition
$(u_1,u_2)\to(u_1;u'_2,u'_{2*},v_*)$, we see that a decoupling
transition which increases
$\langle1+|v|^2,\Gamma_t^1+\Gamma_t^2\rangle$ by
$4+2|v_*|^2+|u_1|^2+|u_2|^2$ occurs at
total rate
\[
2\Gamma^0_{t-}\bigl(d(u_1,u_2)
\bigr)\bigl\llVert B(u_1-v_*,\cdot)-B(u_2-v_*,\cdot)\bigr\rrVert
\rho_t(dv_*)\,dt.
\]
By condition (\ref{LCP}),
\[
\bigl\llVert B(u_1-v_*,\cdot)-B(u_2-v_*,\cdot)\bigr\rrVert \le
\kappa\llvert u_1-u_2\rrvert \le\kappa\llvert
v_1-v_2\rrvert
\]
for all pairs $(u_1,u_2)$ in the support of $\Gamma^0_t$ for all $t$.
Hence the drift of $\langle1+|v|^2,\Gamma_t^1+\Gamma_t^2\rangle$
due to decoupling
transitions is no greater than
\[
6\kappa\llvert v_1-v_2\rrvert \bigl\langle2+\llvert
u_1\rrvert ^2+\llvert u_2\rrvert
^2,\Gamma^0_{t-}\bigr\rangle.
\]
On the other hand, by the same estimates used in Proposition~\ref
{JPE}, the drift of $\langle1+|v|^2,\Gamma_t^1+\Gamma_t^2\rangle$
due to branching of uncoupled particles is no greater than
\[
8m_3(t)\bigl\langle1+\llvert v\rrvert ^2,
\Gamma_{t-}^1+\Gamma_{t-}^2\bigr
\rangle.
\]
Hence the following process is a supermartingale:
\begin{eqnarray*}
&& \bigl\langle1+\llvert v\rrvert ^2,\Gamma_t^1+
\Gamma_t^2\bigr\rangle -6\kappa\llvert
v_1-v_2\rrvert \int_0^t
\bigl\langle2+\llvert u_1\rrvert ^2+\llvert
u_2\rrvert ^2,\Gamma^0_s\bigr
\rangle \,ds \\
&&\qquad{}-\int_0^t8m_3(s)\bigl
\langle1+\llvert v\rrvert ^2,\Gamma_s^1+
\Gamma_s^2\bigr\rangle \,ds.
\end{eqnarray*}
Set $g(t)=E_{(0,v_1,v_2)}(\langle1+|v|^2,\Gamma_t^1+\Gamma
_t^2\rangle)$.
Since $\Gamma^0_t\circ p_k^{-1}\le\Lambda_t^k$, by Proposition~\ref{JPE}
\[
E_{(0,v_1,v_2)}\bigl\langle2+\llvert u_1\rrvert ^2+
\llvert u_2\rrvert ^2,\Gamma^0_t
\bigr\rangle\le \bigl(2+\llvert v_1\rrvert ^2+\llvert
v_2\rrvert ^2\bigr)\exp \biggl\{\int_0^t8m_3(s)\,ds
\biggr\}.
\]
Then
\begin{eqnarray*}
g(t) &\le & 6\kappa\bigl(2+\llvert v_1\rrvert ^2+\llvert
v_2\rrvert ^2\bigr)\llvert v_1-v_2
\rrvert \int_0^t\exp \biggl\{\int
_0^s8m_3(r)\,dr \biggr\}\,ds\\
&&{}+\int
_0^t8m_3(s)g(s)\,ds,
\end{eqnarray*}
and the first of the claimed estimates follows by Gronwall's lemma.
For $p>2$, a straightforward modification of this argument, using
$|v'|^p+|v'_*|^p\le C(p)(|v|^p+|v_*|^p)$ and (\ref{POE}), leads to the
second estimate.
\end{pf}

\begin{pf*}{Proof of Proposition~\ref{FSE}}
For all $f\in\mathcal{F}$ and all $v,v'\in\mathbb{R}^d$, we have
\[
\bigl\llvert f(v)\bigr\rrvert \le1+\llvert v\rrvert ^2,\qquad \bigl\llvert
f(v)-f\bigl(v'\bigr)\bigr\rrvert \le\bigl(2+\llvert v\rrvert
^2+\bigl\llvert v'\bigr\rrvert ^2\bigr)
\bigl\llvert v-v'\bigr\rrvert.
\]
To see the second inequality, note that
\begin{eqnarray*}
\bigl\llvert f(v)-f\bigl(v'\bigr)\bigr\rrvert &=& \bigl\llvert
\bigl(1+\llvert v\rrvert ^2\bigr)\hat f(v)-\bigl(1+\bigl\llvert
v'\bigr\rrvert ^2\bigr)\hat f\bigl(v'
\bigr)\bigr\rrvert
\\
&\le & \bigl(1+\llvert v\rrvert ^2\bigr)\bigl\llvert \hat f(v)-\hat f
\bigl(v'\bigr)\bigr\rrvert +\bigl\llvert \llvert v\rrvert
^2-\bigl\llvert v'\bigr\rrvert ^2\bigr
\rrvert \bigl\llvert \hat f\bigl(v'\bigr)\bigr\rrvert\\
& \le & \bigl(1+
\llvert v\rrvert ^2+\llvert v\rrvert +\bigl\llvert v'
\bigr\rrvert \bigr)\bigl\llvert v-v'\bigr\rrvert
\end{eqnarray*}
and then symmetrize.
We write the proof for the case $\|f\|=1$ and $s=0$. Set $f_0=E_{0t}f$.
By Proposition~\ref{JPE}, for all $v\in\mathbb{R}^d$,
\begin{equation}
\label{FHG} \bigl\llvert f_0(v)\bigr\rrvert \le\bigl(1+\llvert v
\rrvert ^2\bigr)\exp \biggl\{\int_0^t8m_3(s)\,ds
\biggr\}.
\end{equation}
We have
\[
f_0(v_1)-f_0(v_2)=E_{(0,v_1,v_2)}
\bigl(\bigl\langle f\circ p_1-f\circ p_2,\tilde \Gamma
^0_t\bigr\rangle+\bigl\langle f\circ\pi_1,
\tilde\Gamma_t^1\bigr\rangle-\bigl\langle f\circ
\pi_2,\tilde\Gamma_t^2\bigr\rangle\bigr).
\]
So, since $|u_1-u_2|\le|v_1-v_2|$ for all $(u_1,u_2)\in\operatorname
{supp}\Gamma^0_t$,
\[
\bigl\llvert f_0(v_1)-f_0(v_2)
\bigr\rrvert \le E_{(0,v_1,v_2)}\bigl(\bigl\langle2+\llvert u_1
\rrvert ^2+\llvert u_2\rrvert ^2,
\Gamma^0_t\bigr\rangle \llvert v_1-v_2
\rrvert +\bigl\langle1+\llvert v\rrvert ^2,\Gamma _t^1+
\Gamma_t^2\bigr\rangle\bigr).
\]
By Proposition~\ref{JPE},
\begin{eqnarray*}
&& E_{(0,v_1,v_2)}\bigl\langle2+\llvert u_1\rrvert ^2+
\llvert u_2\rrvert ^2,\Gamma^0_t
\bigr\rangle
\\
&&\qquad \le E_{(0,v_1)}\bigl\langle1+\llvert v\rrvert ^2,
\Lambda_t\bigr\rangle +E_{(0,v_2)}\bigl\langle1+\llvert v\rrvert
^2,\Lambda _t\bigr\rangle \\
&&\qquad\le  \bigl(2+\llvert
v_1\rrvert ^2+\llvert v_2\rrvert
^2\bigr)\exp \biggl\{\int_0^t8m_3(s)\,ds
\biggr\}.
\end{eqnarray*}
We combine this with Lemma~\ref{GHU} to obtain
\[
\bigl\llvert f_0(v_1)-f_0(v_2)
\bigr\rrvert \le(1+6\kappa t) \bigl(2+\llvert v_1\rrvert
^2+\llvert v_2\rrvert ^2\bigr)\llvert
v_1-v_2\rrvert \exp \biggl\{\int_0^t8m_3(s)\,ds
\biggr\},
\]
which implies that
\[
\bigl\llvert \hat f_0(v_1)-\hat f_0(v_2)
\bigr\rrvert \le3(1+6\kappa t)\llvert v_1-v_2\rrvert \exp
\biggl\{\int_0^t8m_3(s)\,ds \biggr\}
\]
and in conjunction with (\ref{FHG}) gives the claimed estimate.
\end{pf*}

\begin{pf*}{Proof of Proposition~\ref{FSF}}
It will suffice to consider the case $\|f\|=1$. Write $f_s=E_{st}f$.
Let $(\Lambda_t)_{t\ge s}$ and $(\Lambda'_t)_{t\ge s'}$ be independent
linearized Kac processes starting
from $\delta_{v_0}$ at times $s$ and $s'$, respectively.
Write $T$ for the first branch time of $(\Lambda_t)_{t\ge s}$ and
$V_*,V',V'_*$ for the velocities
of the new particles formed in $(\Lambda_t)_{t\ge s}$ at time $T$.
By the Markov property of the branching process and using Proposition~\ref{JPE}, on the event $\{T\le s'\}$,
\begin{eqnarray*}
&&\bigl\llvert E\bigl(\bigl\langle f,\tilde\Lambda_t-\tilde
\Lambda'_t\bigr\rangle|T,V_*,V',V'_*
\bigr)\bigr\rrvert \\
&&\qquad= \bigl\llvert f_T\bigl(V'
\bigr)+f_T\bigl(V_*'\bigr)-f_T(V_*)-f_{s'}(v_0)
\bigr\rrvert
\\
&&\qquad\le \bigl(4+\llvert v_0\rrvert ^2+\llvert V_*\rrvert
^2+\bigl\llvert V'\bigr\rrvert ^2+\bigl
\llvert V_*'\bigr\rrvert ^2\bigr) \exp \biggl\{\int
_s^t8m_3(r)\,dr \biggr\}
\end{eqnarray*}
while
\[
E\bigl(\bigl\langle f,\tilde\Lambda_t-\tilde\Lambda'_t
\bigr\rangle|T>s'\bigr)=0.
\]
Now $|V'|^2+|V_*'|^2=|v_0|^2+|V_*|^2$, so
\begin{eqnarray*}
\bigl\llvert f_s(v_0)-f_{s'}(v_0)
\bigr\rrvert &=& \bigl\llvert E\bigl\langle f,\tilde\Lambda_t-\tilde
\Lambda _t'\bigr\rangle\bigr\rrvert \\
&\le&  E \bigl(
\bigl(4+2\llvert v_0\rrvert ^2+2\llvert V_*\rrvert
^2\bigr)1_{\{T\le s'\}} \bigr)\exp \biggl\{\int
_s^t8m_3(r)\,dr \biggr\},
\end{eqnarray*}
and, using the inequality $|v-v_*|(4+2|v|^2+2|v_*|^2)\le
5(1+|v|^3)(|1+|v_*|^3)$, we have
\begin{eqnarray*}
&& E \bigl(\bigl(4+2\llvert v_0\rrvert ^2+2\llvert V_*
\rrvert ^2\bigr)1_{\{T\le s'\}} \bigr)\\
 &&\qquad\le  \int_s^{s'}\!\!
\int_{\mathbb
{R}^d}\llvert v_0-v_*\rrvert \bigl(4+2\llvert
v_0\rrvert ^2+2\llvert v_*\rrvert ^2\bigr)
\rho _r(dv_*)\,dr
\\
&&\qquad\le  5\bigl(1+\llvert v_0\rrvert ^3\bigr)\int
_s^{s'}m_3(r)\,dr,
\end{eqnarray*}
so
\[
\bigl\llvert f_s(v_0)-f_{s'}(v_0)
\bigr\rrvert \le5\bigl(1+\llvert v_0\rrvert ^3\bigr)\exp
\biggl\{\int_s^t8m_3(r)\,dr \biggr\}
\int_s^{s'}m_3(r)\,dr.
\]
\upqed\end{pf*}

\begin{pf*}{Proof of Proposition~\ref{FMU}}
Recall that now $\rho_t=(\mu_t^N+\mu_t^{N'})/2$ and $m_3(t)=\langle
1+|v|^3,\mu_t^N+\mu_t^{N'}\rangle/2$.
In particular $m_3(t)\le1+(N')^{3/2}<\infty$ for all~$t$.
Set $M=M^N-M^{N'}$, and write $M^\pm$ for the positive and negative
parts of the signed measure $M$ on $[0,\infty)\times\mathbb{R}^d$.
Consider a branching particle system in $V^*$, with the same branching
rules as $(\Lambda_t^*)_{t\ge s}$ above, but where,
instead of starting with just one particle at time $s$,
we initiate particles randomly in the system according to a Poisson
random measure on $[0,\infty)\times V^*$ of intensity
\[
\theta(dt,dv)=
\cases{ \delta_0(dt)\mu_0^N(dv)+M^+(dt,dv),&
\quad$\mbox{on }[0,\infty)\times V^+$,\vspace*{3pt}
\cr
\delta_0(dt)
\mu_0^{N'}(dv)+M^-(dt,dv),&\quad $\mbox{on }[0,\infty)\times
V^-$.}
\]
We use the same notation as above for the empirical measures associated
to the branching process, and signify the new rule for initiating particles
by writing now $E$ for the expectation.
Define, for $t\ge0$, a signed measure ${\lambda}_t$ on $\mathbb
{R}^d$ by
\begin{equation}
\label{defl} {\lambda}_t=E(\tilde\Lambda_t)=\int
_{[0,t]\times V}E_{(s,v)}(\tilde \Lambda _t)
\theta(ds,dv).
\end{equation}
Then, by Proposition~\ref{JPE},
\[
\bigl\langle1+\llvert v\rrvert ^2,\llvert {\lambda}_t
\rrvert \bigr\rangle \le\exp \biggl\{\int_0^t8m_3(s)\,ds
\biggr\}\int_{[0,t]\times\mathbb{R}
^d}\bigl(1+\llvert v\rrvert ^2
\bigr)\bigl\llvert \theta(ds,dv)\bigr\rrvert
\]
and, by estimate (\ref{MDV}),
\begin{eqnarray*}
&&\int_{[0,t]\times\mathbb{R}^d}\bigl(1+\llvert v\rrvert ^2\bigr)\bigl
\llvert \theta(ds,dv)\bigr\rrvert \\
&&\qquad\le\bigl\langle1+\llvert v\rrvert
^2,\mu_0^N+\mu_0^{N'}
\bigr\rangle+\int_{[0,t]\times
\mathbb{R}
^d}\bigl(1+\llvert v\rrvert ^2
\bigr)\bigl\llvert M(ds,dv)\bigr\rrvert <\infty.
\end{eqnarray*}
We see in particular that $\langle1+\llvert v\rrvert ^2,\llvert {\lambda}_t\rrvert \rangle$ is
bounded on
compacts in $t$.

Under $E_{(s,v)}$, the pair of empirical processes of positive and
negative particles $(\Lambda_t^+,\Lambda_t^-)_{t\ge s}$ evolves as a
Markov chain,
which makes jumps $(\delta_{v'}+\delta_{v'_*}-\delta_v,\delta
_{v_*})$ at rate $2\Lambda
_{t-}^+(dv)\rho_t(dv_*)B(v-v_*,d\sigma)\,dt$
and makes jumps $(\delta_{v_*},\delta_{v'}+\delta_{v'_*}-\delta_v)$
at rate $2\Lambda
_{t-}^-(dv)\rho_t(dv_*)B(v-v_*,d\sigma)\,dt$.
So, using Proposition~\ref{JPE} for integrability,
under $E_{(s,v)}$, for any bounded measurable function $f$, the
following process is a martingale:
\[
\langle f,\tilde\Lambda_t\rangle-\int_s^t
\bigl\langle f,2Q(\rho _r,\tilde\Lambda_r)\bigr\rangle
\,dr,\qquad  t\ge s.
\]
Taking expectations and setting $f_{st}(v)=E_{st}f(v)=E_{(s,v)}\langle
f,\tilde\Lambda_t\rangle$, we obtain
\begin{equation}
\label{FST}
f_{st}(v)=f(v)+\int_s^t
\bigl\langle f,2Q\bigl(\rho_r,E_{(s,v)}(\tilde\Lambda
_r)\bigr)\bigr\rangle \,dr.
\end{equation}
Then
\begin{eqnarray*}
\langle f,{\lambda}_t\rangle &=& \int_{[0,t]\times V}E_{(s,v)}
\langle f,\tilde\Lambda_t\rangle \theta(ds,dv)
\\
&=& \bigl\langle f_{0t},\mu_0^N-
\mu_0^{N'}\bigr\rangle+\int_0^t
\langle f_{st},dM_s\rangle
\\
&=&\bigl\langle f,\mu_0^N-\mu_0^{N'}
\bigr\rangle+\langle f,M_t\rangle+\int_0^t
\bigl\langle f,2Q(\rho_r,{\lambda }_r)\bigr\rangle \,dr.
\end{eqnarray*}
Here, we used (\ref{defl}) for the first equality, and for the third
we substituted for $f_{0t}$ and $f_{st}$ using (\ref{FST}) and then
rearranged the integrals using Fubini to make ${\lambda}_r$, as given
by (\ref{defl}), appear on the inside.
Since $f$ is an arbitrary bounded measurable function, we have shown that
\begin{equation}
\label{LEB}
{\lambda}_t=\bigl(\mu_0^N-
\mu_0^{N'}\bigr)+M_t+\int_0^t2Q(
\rho_s,{\lambda}_s)\,ds.
\end{equation}
Note the estimate of total variation,
\[
\bigl\llVert Q(\rho_t,{\lambda}_t)\bigr\rrVert \le4\int
_{\mathbb{R}^d\times\mathbb
{R}^d}\llvert v-v_*\rrvert \rho _t(dv)\llvert {
\lambda}_t\rrvert (dv_*)\le6\bigl\langle1+\llvert v\rrvert
^2,\llvert {\lambda}_t\rrvert \bigr\rangle.
\]
For the second inequality, we used $\langle1,\rho_t\rangle=\langle
|v|^2,\rho_t\rangle=1$
and $2|v-v_*|\le2+|v|^2+|v_*|^2$.
For any interval $(s,t]$ on which neither $(\mu^N_t)_{t\ge0}$ nor
$(\mu^{N'}_t)_{t\ge0}$ jump,
using estimate (\ref{TVMST}), we deduce that
\[
\bigl\langle1+\llvert v\rrvert ^2,\llvert M_t-M_s
\rrvert \bigr\rangle\le24\int_s^tm_3(r)\,dr.
\]
On the other hand
\[
\biggl\langle1+\llvert v\rrvert ^2,\int_s^t
\bigl\llvert Q(\rho_r,{\lambda}_r)\bigr\rrvert \,dr
\biggr\rangle\to0
\]
as $t\downarrow s$, for all $s\ge0$. Hence from equation (\ref{LEB}) we
deduce that $(1+|v|^2){\lambda}_t$ is right continuous
in total variation.

Set ${\lambda}_t'=\mu_t^N-\mu_t^{N'}$, and note from (\ref{MD})
that $({\lambda}'_t)_{t\ge0}$ also satisfies (\ref{LEB}).
We subtract to see that $\delta_t={\lambda}_t-{\lambda}'_t$ satisfies
\[
\delta_t=\int_0^t2Q(
\rho_s,\delta_s)\,ds.
\]
Set $\nu_t=2Q(\rho_t,\delta_t)$.
Then\vspace*{1pt} $\|\nu_t\|\le12\langle1+|v|^2,|\delta_t|\rangle$ and,
on any interval $(s,t]$ when neither $(\mu^N_t)_{t\ge0}$ nor $(\mu
^{N'}_t)_{t\ge0}$ jump,
\[
\llVert \nu_t-\nu_s\rrVert =\bigl\llVert 2Q(
\rho_s,{\lambda}_t-{\lambda}_s)\bigr\rrVert
\le 12\bigl\langle 1+\llvert v\rrvert ^2,\llvert {
\lambda}_t-{\lambda}_s\rrvert \bigr\rangle.
\]
The process of signed measures $(\nu_t)_{t\ge0}$ is thus locally
bounded and right continuous in total variation.
Hence the measure $\int_0^T|\nu_t|\,dt+|\nu_T|$ is finite, and $\nu
_t$ is absolutely continuous with respect to this measure for all $t\in[0,T]$.

We apply Lemma~\ref{BML}, with $\mu_0=0$, to obtain a measurable map
$\sigma \dvtx  [0,\infty)\times V\to\{-1,0,1\}$
such that $\delta_t=\sigma_t|\delta_t|$ and $|\delta_t|=\int_0^t\sigma_s\nu_s\,ds$.
Set $\check\sigma_s(v)=(1+|v|^2)\sigma_s(v)$. Then
\begin{eqnarray*}
&&\bigl\langle1+\llvert v\rrvert ^2,\llvert \delta_t
\rrvert \bigr\rangle \\
&&\qquad=\int_E 2\bigl\{\check
\sigma_s\bigl(v'\bigr)+\check\sigma_s
\bigl(v_*'\bigr)-\check\sigma _s(v)-\check\sigma
_s(v_*)\bigr\}\\
&&\qquad\qquad\hspace*{1.5pt}{}\times 1_{(0,t]}(s)B(v-v_*,d\sigma)
\rho_s(dv_*)\delta_s(dv)\,ds
\\
&&\qquad\le \int_0^t\!\!\int_{\mathbb{R}^d\times\mathbb
{R}^d}4
\bigl(1+\llvert v_*\rrvert ^2\bigr)\llvert v-v_*\rrvert \rho
_s(dv_*)\llvert \delta_s\rrvert (dv)\,ds\\
&&\qquad\le\int
_0^t4\bigl\langle1+\llvert v\rrvert
^2,\llvert \delta _s\rrvert \bigr\rangle
m_3(s)\,ds.
\end{eqnarray*}
But $\int_0^tm_3(s)\,ds<\infty$, so $\delta_t=0$, for all $t$.
\end{pf*}

\section{Proof of Theorem \texorpdfstring{\protect\ref{MR}}{1.1}}\label{PMR}
We will write the proof for $d\ge3$, leaving the minor modifications
necessary for $d=2$ to the reader.
Fix $p\in(2,\infty)$ and ${\lambda}\in[1,\infty)$.
Suppose that $(\mu^N_t)_{t\ge0}$ and $(\mu^{N'}_t)_{t\ge0}$ are Kac
processes in $\mathcal{S}_N$ and $\mathcal{S}_{N'}$, respectively,
with $\langle|v|^p,\mu_0^N\rangle\le{\lambda}$ and $\langle
|v|^p,\mu_0^{N'}\rangle\le
{\lambda}$.
Set $\rho_t=(\mu^N_t+\mu^{N'}_t)/2$.
Fix $t\in[0,T]$ and a function $f_t\in\mathcal{F}$.
Define a random function $f$ on $[0,T]\times\mathbb{R}^d$ by
\[
f(s,v)=f_s(v)=E_{(s\wedge t,v)}\langle f_t,\tilde
\Lambda_t\rangle,
\]
where $(\Lambda_t^*)_{t\ge s}$ is a linearized Kac process in environment
$(\rho_t)_{t\ge0}$.
Note that we have extended $f$ as a constant in time from $t$ to $T$.
We have, by Proposition~\ref{FMU},
\begin{equation}
\label{FMUE}
\bigl\langle f_t,\mu_t^N-
\mu_t^{N'}\bigr\rangle=\bigl\langle f_0,
\mu_0^N-\mu _0^{N'}\bigr\rangle+\int
_0^t\bigl\langle f_s,dM_s^N
\bigr\rangle-\int_0^t\bigl\langle
f_s,dM^{N'}_s\bigr\rangle.
\end{equation}
Write $m_q(t)=\langle1+|v|^q,\rho_t\rangle$, as above.
By Proposition~\ref{ME}, for $q<p+1$, there is a constant
$C(B,p,q)<\infty$ such that
\begin{equation}
\label{HFMA} \mathbb{E}\int_0^Tm_q(s)\,ds
\le C\bigl(T^{p+1-q}+T\bigr){\lambda}.
\end{equation}
Set
\begin{equation}
\label{FMA}
A=3(1+6\kappa T)\exp \biggl\{\int_0^T8m_3(s)\,ds
\biggr\}.
\end{equation}
By Proposition~\ref{FSE}, for all $s\le t$, we have
\begin{equation}
\label{LE} \llVert f_s\rrVert \le A
\end{equation}
so
\begin{equation}
\label{IV} \bigl\langle f_0,\mu_0^N-
\mu_0^{N'}\bigr\rangle\le AW\bigl(\mu_0^N,
\mu_0^{N'}\bigr).
\end{equation}
The main step of the proof is to bound the second and third terms on
the right in~(\ref{FMUE}), uniformly in $t\in[0,T]$ and $f_t\in
\mathcal{F}$.
We will derive estimates for the second term, which then apply also to
the third, because $N\le N'$.
The notation conceals the fact that the integrand $f_s$ depends on the
terminal time $t$.
Worse, $f_s$ depends on $(\mu^N_r+\mu_r^{N'})_{s\le r\le t}$, so is
anticipating, and martingale estimates cannot be applied directly even
at the individual time $t$.

For $p\ge3$, set $\beta=1$ and $Z=\sup_{t\in[0,T]}m_3(t)$.
By Propositions \ref{ME} and \ref{FSF}, we have $\mathbb{E}(Z)\le
1+(1+CT){\lambda}$ and, for all $v\in\mathbb{R}^d$ and $s,s'\in
[0,T]$ with
$s\le s'$,
\begin{equation}
\label{LTE} \bigl\llvert f_s(v)-f_{s'}(v)\bigr\rrvert
\le A'\bigl(1+\llvert v\rrvert ^3\bigr)
\bigl(s'-s\bigr)^\beta,
\end{equation}
where
\begin{equation}
\label{FMAP} A'=5Z\exp \biggl\{\int_0^T8m_3(s)\,ds
\biggr\}.
\end{equation}

For $p\in(2,3)$, set $\beta=(p-2)/2$, and set
\[
Z=2\sup_{t\in[1,T]}m_3(t)+\sum
_{\ell\in\mathbb{N}}2^{(\beta
-1)\ell+1}\beta ^{-1}\sup
_{t\in[2^{-\ell},2^{-\ell+1}]}m_3(t).
\]
By Proposition~\ref{ME}, there is a constant $C(B,p)<\infty$ such
that, for $t\le T$,
\[
\mathbb{E} \Bigl(\sup_{s\in[t,T]}m_3(s) \Bigr)\le C
\bigl(t^{p-3}\vee T\bigr){\lambda}
\]
so
\begin{equation}
\label{FMAQ} \mathbb{E}(Z)\le C{\lambda} \biggl(T+\sum
_{\ell\in\mathbb
{N}}2^{(\beta-1)\ell
}2^{-\ell(p-3)} \biggr)=C{\lambda}
\bigl(T+1/\bigl(2^\beta-1\bigr)\bigr).
\end{equation}
Note that $m_3(t)\le(\beta t^{\beta-1}+1)Z/2$ for all $t\le T$, so for
$s\le s'\le T$ with $s'-s\le1$,
\[
\int_s^{s'}m_3(t)\,dt\le\bigl(
\bigl(s'-s\bigr)^\beta+\bigl(s'-s\bigr)
\bigr)Z/2\le\bigl(s'-s\bigr)^\beta Z.
\]
Hence, by Proposition~\ref{FSF}, (\ref{LTE}) remains valid for $p\in
(2,3)$, provided $s'-s\le1$ and $\beta$ and $Z$ have their new meanings.

Fix $r\in(0,1]$ and $R\in[1,\infty)$ such that $T/r$ and $R/r$ are
integers, and set $n=(T/r)\times(R/r)^d$. Set $B(R)=(-R,R]^d$.
There exist $s_1,\dots,s_n\in[0,T)$ and $v_1,\dots,v_n\in B(R)$ such
that $B_1\cup\cdots\cup B_n=(0,T]\times B(R)$, where
$B_k=(s_k,v_k)+(0,r]\times(-r,r]^d$.
Write
\[
f=\sum_{k=1}^na_kf^{(k)}+g,
\]
where $a_k$ is the average value of $\hat f$ on $B_k$ and
$f^{(k)}(s,v)=\check1_{B_k}(s,v)=(1+|v|^2)1_{B_k}(s,v)$.
Then
\begin{equation}
\label{IF} \int_0^t\bigl\langle
f_s,dM^N_s\bigr\rangle=\sum
_{k=1}^na_kM^{(k)}_t+
\int_0^t\bigl\langle g_s,dM_s^N
\bigr\rangle,
\end{equation}
where
\[
M^{(k)}_t=\int_0^t\bigl
\langle f^{(k)}_s,dM^N_s\bigr
\rangle.
\]
Now, by (\ref{LE}), for all $k$, we have $|a_k|\le A$ and, for
$v,v'\in B_k$,
\[
\bigl\llvert \hat f_s(v)-\hat f_s\bigl(v'
\bigr)\bigr\rrvert \le A\bigl\llvert v-v'\bigr\rrvert \le2
\sqrt{d}Ar.
\]
By (\ref{LTE}), we have, for $(s,v)\in B_k$,
\[
\bigl(1+\llvert v\rrvert ^2\bigr)\bigl\llvert \hat
f_s(v)-a_k(v)\bigr\rrvert \le A'\bigl(1+
\llvert v\rrvert ^3\bigr)r^\beta,
\]
where $a_k(v)$ is the average value of $\hat f$ on $(s_k,s_k+r]\times\{
v\}$.
Hence, for $(s,v)\in B_k$,
\begin{eqnarray*}
\bigl\llvert g_s(v)\bigr\rrvert &=& \bigl(1+\llvert v\rrvert
^2\bigr)\bigl\llvert \hat f_s(v)-a_k(v)+a_k(v)-a_k
\bigr\rrvert \\
&\le&  A'\bigl(1+\llvert v\rrvert ^3
\bigr)r^\beta +2\sqrt dA\bigl(1+\llvert v\rrvert ^2\bigr)r.
\end{eqnarray*}
On the other hand, $|g_s(v)|\le A(1+|v|^2)$ for all $v\in\mathbb
{R}^d\setminus B(R)$.

Set $Q_t=Q_t^N=\sum_{k=1}^n|M^{(k)}_t|^2$. Then
\begin{equation}
\label{IG} \Biggl\llvert \sum_{k=1}^na_kM^{(k)}_t
\Biggr\rrvert \le A\sqrt{nQ_t}.
\end{equation}
Fix $q\in(3,p+1)$.
Note that, for $s\in(0,T]$,
\begin{eqnarray*}
&&\sum_{k=1}^n\bigl\{
f^{(k)}_s\bigl(v'\bigr)+f^{(k)}_s
\bigl(v'_*\bigr)-f^{(k)}_s(v)-f^{(k)}_s(v_*)
\bigr\}^2
\\
&&\qquad  \le4\sum_{k=1}^n\bigl\{\check{
\check1}_{B_k}\bigl(s,v'\bigr)+\check {\check1}
_{B_k}\bigl(s,v_*'\bigr)+\check{\check1}_{B_k}(s,v)+
\check{\check 1}_{B_k}(s,v_*)\bigr\}
\\
&&\qquad=4\bigl\{\check{\check1}_{B(R)}\bigl(v'\bigr)+\check{
\check 1}_{B(R)}\bigl(v_*'\bigr)+\check{
\check1}_{B(R)}(v)+\check{\check1} _{B(R)}(v_*)\bigr\}\\
&&\qquad \le
CR^{(5-q)^+}\bigl(1+\llvert v\rrvert ^{q-1}+\llvert v_*\rrvert
^{q-1}\bigr)
\end{eqnarray*}
for some constant $C<\infty$, depending only on $d$ and $q$.
So, by Doob's $L^2$-inequality,
\begin{eqnarray}
\nonumber
\mathbb{E} \Bigl(\sup_{t\le T}Q_t \Bigr) &
\le&\frac{4}{N^2}\sum_{k=1}^n
\mathbb{E}\int_E\bigl\{ f^{(k)}_s
\bigl(v'\bigr)+f_s^{(k)}\bigl(v'_*
\bigr)-f_s^{(k)}(v)-f_s^{(k)}(v_*)
\bigr\} ^2 \\
&&\label{QE}\hspace*{34pt}\qquad{}\times 1_{(0,T]}(s)\bar m(dv,dv_*,d\sigma,ds)
\\
\nonumber
&\le&\frac{CR^{(5-q)^+}}N\mathbb{E}\int_0^T\bigl
\langle\llvert v\rrvert ^q,\mu ^N_s\bigr
\rangle \,ds.
\end{eqnarray}
On the other hand
\begin{equation}
\label{IH} \quad\biggl\llvert \int_0^t\bigl
\langle g_s,dM^N_s\bigr\rangle\biggr\rrvert
\le C\bigl(A'r^\beta R^{(4-q)^+}+Ar+AR^{3-q}
\bigr)\int_0^T\bigl\langle1+\llvert v\rrvert
^{q-1},\bigl\llvert dM^N_s\bigr\rrvert \bigr
\rangle\hspace*{-18pt}
\end{equation}
and
\begin{eqnarray}
&&\mathbb{E}\int_0^T\bigl\langle1+\llvert v
\rrvert ^{q-1},\bigl\llvert dM^N_s\bigr\rrvert
\bigr\rangle
\nonumber
\\
&&\qquad \le\frac{2}N\mathbb{E}\int_E\bigl\{ 4+\bigl
\llvert v'\bigr\rrvert ^{q-1}+\bigl\llvert
v'_*\bigr\rrvert ^{q-1}+\llvert v\rrvert ^{q-1}+
\llvert v_*\rrvert ^{q-1}\bigr\}
\nonumber
\\[-8pt]
\label{AE}
\\[-8pt]
\nonumber
&&\hspace*{22pt}\qquad\qquad{}\times1_{(0,T]}(s)\bar m(dv,dv_*,d
\sigma,ds)
\\
\nonumber
&&\qquad \le C\mathbb{E}\int_0^T\bigl\langle\llvert
v\rrvert ^q,\mu^N_s\bigr\rangle
\,ds.
\end{eqnarray}

We combine (\ref{FMUE}), (\ref{HFMA}), (\ref{FMA}), (\ref{IV}),
(\ref{FMAP}), (\ref{FMAQ}), (\ref{IF}), (\ref{IG}), (\ref{QE}),
(\ref{IH}) and (\ref{AE})
to see that, for all $\varepsilon\in(0,1]$, there is a constant
$C<\infty$,
depending only on $B,d,\varepsilon,{\lambda},p,q$ and $T$, such that,
for all
$N,N'\in\mathbb{N}$ with $N\le N'$, with probability exceeding
$1-\varepsilon$, we
have, for all $t\in[0,T]$,
\[
W\bigl(\mu_t^N,\mu_t^{N'}\bigr)\le
C\bigl(W\bigl(\mu^N_0,\mu^{N'}_0
\bigr)+\sqrt {R^{(5-q)^+}n/N}+rR^{(4-q)^+}+R^{3-q}\bigr).
\]
An optimization over $q$, $r$ and $R$ now shows the existence of an
$\alpha
(d,p)>0$ for which the estimate claimed in Theorem~\ref{MR} holds.

For large $p$, the reader may check the optimization yields a value for
$\alpha(d,p)$ close to $1/(d+3)$.
The proof given can be varied by replacing the one-step discrete
approximation by a chaining argument.
See the proof of Proposition~\ref{DAT} for this idea in a simple context.
This gives $\alpha(d,p)=1/(d+1)$ for $p$ sufficiently large.
We omit the details because Theorem~\ref{MR10} gives a stronger result.
Here $d+1$ is the dimension of space--time, reflecting the fact that we
maximize over a class of functions on $[0,T]\times\mathbb{R}^d$.
This is wasteful because, in fact, we only need to maximize over a
certain process of functions $(f_s \dvtx  s\in[0,t])$ associated to $t\in
[0,T]$ and $f=f_t$
and then over a class of functions $f$ on $\mathbb{R}^d$.
In the next three sections, we exploit the structure of the process
$(f_s \dvtx  s\in[0,t])$ to obtain an improved bound.

\section{Continuity of the linearized Kac process in its
environment}\label{MKP}
We\break showed in Propositions \ref{FSE} and \ref{FSF} that the linearized
Kac process is continuous in its initial data.
For the proof of our main estimate with optimal rate $N^{-1/d}$, we
will need also continuity in the environment.
The following notation will be convenient.
For $p\in[2,\infty)$ and a function $f$ on $\mathbb{R}^d$, we will write
$\hat f^{(p)}$ for the reweighted function $\hat
f^{(p)}(v)=f(v)/(1+|v|^p)$ and write $\|f\|_{(p)}$ for the smallest
constant such that, for all $v,v'\in\mathbb{R}^d$, we have
\[
\bigl\llvert \hat f^{(p)}(v)\bigr\rrvert \le\llVert f\rrVert
_{(p)},\qquad  \bigl\llvert \hat f^{(p)}(v)-\hat f^{(p)}
\bigl(v'\bigr)\bigr\rrvert \le\llVert f\rrVert _{(p)}\bigl
\llvert v-v'\bigr\rrvert.
\]
Denote by $\mathcal{F}(p)$ the set of all functions $f$ on $\mathbb
{R}^d$ with $\|f\|
_{(p)}\le1$.
We earlier wrote $\mathcal{F}$ for $\mathcal{F}(2)$ and $\|f\|$ for
$\|f\|_{(2)}$.
We will use the cases $p=2$ and $p=3$.
Suppose that $(\rho_t^1)_{t\ge0}$ and $(\rho_t^2)_{t\ge0}$ are
processes of measures on $\mathbb{R}^d$, both satisfying (\ref{SKP}).
Given $t\ge0$ and a function $f$ of quadratic growth on $\mathbb{R}^d$,
define for $s\in[0,t]$ and $v\in\mathbb{R}^d$, and for $j=1,2$,
\[
E_{st}^jf(v)=E_{(s,v)}\bigl\langle f,\tilde
\Lambda^j_t\bigr\rangle,
\]
where $(\Lambda^{j,*}_t)_{t\ge s}$ is a linearized Kac process with
environment $(\rho^j_t)_{t\ge0}$ starting from $v$ at time $s$.
We will use the following notation:
\[
d_p(t)=\bigl\langle1+\llvert v\rrvert ^p,\bigl\llvert
\rho_t^1-\rho_t^2\bigr\rrvert
\bigr\rangle, \qquad\bar m_p(t)=\bigl\langle 1+\llvert v\rrvert
^p,\rho_t^1+\rho_t^2
\bigr\rangle.
\]

\begin{proposition}\label{FSG}
For all $p\in[2,\infty)$, there is a constant $C(p)<\infty$ with the
following properties.
Let $f$ be a function on $\mathbb{R}^d$ with $|f(v)|\le1+|v|^p$ for
all~$v$.
Then, for all $s,t\ge0$ with $s\le t$ and all $v\in\mathbb{R}^d$, we have
\begin{eqnarray*}
&&\bigl\llvert E_{st}^1f(v)-E_{st}^2f(v)
\bigr\rrvert \\
&&\qquad\le C(p) \bigl(1+\llvert v\rrvert ^{p+1}\bigr)\exp \biggl
\{C(p)\int_s^t\bar m_{p+2}(r)\,dr
\biggr\}\int_s^td_{p+1}(r)\,dr.
\end{eqnarray*}
Assume that the collision kernel satisfies condition (\ref{LCP}).
Then $C(p)$ may be chosen so that
\[
\bigl\llVert E_{st}^1f-E_{st}^2f
\bigr\rrVert _{(p+1)}\le C(p)\kappa\llVert f\rrVert _{(p)}\exp
\biggl\{ C(p)\int_s^t\bar m_{p+2}(r)\,dr
\biggr\}\int_s^td_{p+1}(r)\,dr.
\]
\end{proposition}

Our first step toward a proof of Proposition~\ref{FSG} is to describe
a coupling of $(\Lambda^{1,*}_t)_{t\ge0}$ and $(\Lambda
^{2,*}_t)_{t\ge0}$
when both process start from $v$ at time $0$.
For this, we take as type space the set $\hat V_0\cup\hat V_1\cup\hat
V_2$, where $\hat V_j=\mathbb{R}^d\times\{j\}$ for $j=0,1,2$.
Particles with types in $\hat V_0$ are called coupled, the others are uncoupled.
Consider the branching process with the following branching transitions.
For a particle $v$ in $\hat V_0$, there are three possible transitions.
First, at rate $2(\rho^1_t\wedge\rho^2_t)(dv_*)B(v-v_*,d\sigma)\,dt$, we
replace $v$ by three particles $v_*$, $v'$ and $v_*'$ in $\hat V_0$.
Second, at rate $2(\rho_t^1-\rho_t^2)^+(dv_*)B(v-v_*,d\sigma)\,dt$, we
replace $v$ by three particles $v_*$, $v'$ and $v_*'$ in $\hat V_1$ and
one particle $v$ in $\hat V_2$.
Third, at rate $2(\rho_t^2-\rho_t^1)^+(dv_*)B(v-v_*,d\sigma)\,dt$, we
replace $v$ by one particle $v$ in $\hat V_1$ and three particles
$v_*$, $v'$ and $v_*'$ in $\hat V_2$.
The second and third are called uncoupling transitions.
The transitions for uncoupled particles are as in the original
branching process; that is,
for $j=1,2$ and $v$ in $\hat V_j$, at rate $2\rho
_t^j(dv_*)B(v-v_*,d\sigma
)\,dt$ we replace $v$ by particles $v_*$, $v'$ and $v_*'$, also in $\hat V_j$.
For $v\in\mathbb{R}^d$ and $s\ge0$, and for $j=0,1,2$, write $\hat
\Gamma_t^j$
for the un-normalized empirical distribution of particles in $\hat V_j$
when we
initiate the branching process with a single particle $v$ in $\hat V_0$
at time $s$.
Define analogously the lifted processes $(\hat\Gamma_t^{j,*})_{t\ge s}$
in $\hat V_j^*$.
For $j=1,2$, set
\[
\hat\Lambda_t^{j,*}=\hat\Gamma_t^{0,*}
\circ\hat\pi^{-1}_0+\hat \Gamma _t^{j,*}
\circ\hat\pi^{-1}_j,
\]
where $\hat\pi_j$ is the bijection $\hat V_j^*\to V^*$.
It is straightforward to check that $(\hat\Lambda_t^{j,*})_{t\ge s}$
is a
linearized Kac process with environment $(\rho^j_t)_{t\ge0}$
starting from $v$ in $V^+$ at time $s$.
We have burdened the notation with hats so that we can later refer
simultaneously to this coupling and to the coupling for two different
starting points.

\begin{lemma}\label{DCL}
For all $p\in[2,\infty)$, there is a constant $C(p)<\infty$ such that
\begin{eqnarray*}
&& E_{(0,v_0)}\bigl\langle1+\llvert v\rrvert ^p,\hat
\Gamma^1_t+\hat\Gamma^2_t\bigr
\rangle\\
&&\qquad\le C(p) \bigl(1+\llvert v_0\rrvert ^{p+1}\bigr)\exp
\biggl\{\int_0^tC(p)\bar m_{p+2}(s)\,ds
\biggr\} \int_0^td_{p+1}(s)\,ds.
\end{eqnarray*}
\end{lemma}
\begin{pf}
The process $\langle1+|v|^p,\hat\Gamma_t^1+\hat\Gamma_t^2\rangle$
starts at $0$ and
makes jumps both at uncoupling transitions and due to the branching\vspace*{1pt} of
uncoupled particles.
Uncoupling transitions occur at rate $2B(v-v_*,d\sigma)\hat\Gamma
^0_{t-}(dv)|\rho^1_t-\rho^2_t|(dv_*)\,dt$ and result in jumps of
$4+|v'|^p+|v_*'|^p+|v|^p+|v_*|^p$.
Uncoupled particles branch at rate $2B(v-v_*,d\sigma)(\hat\Gamma
^1_{t-}(dv)\rho_t^1(dv_*)+\hat\Gamma^2_{t-}(dv)\rho_t^2(dv_*))\,dt$
and result in jumps of $2+|v'|^p+|v_*'|^p+|v_*|^2-|v|^p$.
We use the inequalities
\[
4+\bigl\llvert v'\bigr\rrvert ^p+\bigl\llvert
v_*'\bigr\rrvert ^p+\llvert v\rrvert ^p+
\llvert v_*\rrvert ^p\le C(p) \bigl(1+\llvert v\rrvert ^p+
\llvert v_*\rrvert ^p\bigr)
\]
and
\[
\bigl(1+\llvert v\rrvert ^p+\llvert v_*\rrvert ^p\bigr)
\llvert v-v_*\rrvert \le C(p) \bigl(1+\llvert v\rrvert ^{p+1}\bigr)
\bigl(1+\llvert v_*\rrvert ^{p+1}\bigr)
\]
to see that the drift of $\langle1+|v|^p,\hat\Gamma_t^1+\hat\Gamma
_t^2\rangle$ due to
uncoupling transitions is no greater than $C(p)d_{p+1}(t)\langle
1+|v|^{p+1},\hat\Gamma^0_{t-}\rangle$.
On the other hand, inequalities (\ref{POF}) and (\ref{POE}) show that
the drift of $\langle1+|v|^p,\hat\Gamma_t^1+\hat\Gamma_t^2\rangle
$ due to branching
of uncoupled particles
is no greater than $c(p)\bar m_{p+1}(t)\langle1+|v|^p,\hat\Gamma
_{t-}^1+\hat
\Gamma_{t-}^2\rangle$.
Hence the following process is a supermartingale:
\begin{eqnarray*}
&& \bigl\langle1+\llvert v\rrvert ^p,\hat\Gamma_t^1+
\hat\Gamma_t^2\bigr\rangle-C(p)\int_0^t
\bigl\langle1+\llvert v\rrvert ^{p+1},\hat \Gamma^0_s
\bigr\rangle d_{p+1}(s)\,ds\\
&&\qquad{}-c(p)\int_0^t
\bigl\langle1+|v|^p,\hat \Gamma_s^1+\hat
\Gamma_s^2\bigr\rangle \bar m_{p+1}(s)\,ds.
\end{eqnarray*}
On taking expectations, we obtain
\[
g(t)\le\int_0^t\bigl\{C(p)f_{p+1}(s)d_{p+1}(s)+c(p)
\bar m_{p+1}(s)g(s)\bigr\}\,ds,
\]
where $g(t)=E_{(0,v_0)}\langle1+|v|^p,\hat\Gamma_t^1+\hat\Gamma
_t^2\rangle$ and
$f_p(t)=E_{(0,v_0)}\langle1+|v|^p,\hat\Gamma^0_t\rangle$.
By Proposition~\ref{JPE}, we have
\[
f_{p+1}(t)\le\bigl(1+|v_0|^{p+1}\bigr)\exp\int
_0^tc(p+1)\bar m_{p+2}(s)\,ds,
\]
so, for some constant $C(p)<\infty$,
\[
g(t)\le C(p) \bigl(1+|v_0|^{p+1}\bigr)\exp \biggl\{C(p)\int
_0^t\bar m_{p+2}(r)\,dr \biggr\}\int
_0^td_{p+1}(s)\,ds.
\]
\upqed\end{pf}

The second ingredient needed for Proposition~\ref{FSG} is a coupling
of four linearized Kac processes,
with environments $(\rho^1_t)_{t\ge0}$ and $(\rho^2_t)_{t\ge0}$ and
with starting points $v_1$ and $v_2$.
We will specify this coupling in detail, at the cost of some heaviness
of notation,
the understanding of which may be guided by the thought that the
coupling is the obvious one for typed branching at different rates
and is an elaboration of the couplings described in the case of a
single environment or a single starting point above.
To define the coupled processes, we consider a type space which is the
disjoint union of nine sets,
\[
(V_{00}\cup V_{01}\cup V_{02})
\cup(V_{10}\cup V_{11}\cup V_{12})\cup
(V_{20}\cup V_{21}\cup V_{22}).
\]
Here, for $k=0,1,2$, we take $V_{0k}=V_k$ as at (\ref{TS}) and, for
$j=1,2$, we take $V_{jk}=V_k\times\{j\}$.
The first index refers to the environment, a $0$ indicating a particle
present in the branching process in both environments.
The second index refers to the starting point.

The branching rules for a particle $(v_1,v_2)$ in $V_{00}$ are as follows.
There are nine possible transitions.
First, at rate $2(\rho^1_t\wedge\rho^2_t)(dv_*)(B(v_1-v_*,d\sigma
)\wedge B(v_2-v_*,d\sigma))\,dt$
(for all $v_*\in\mathbb{R}^d$ and all $\sigma\in S^{d-1}$), we replace
$(v_1,v_2)$ by three particles $(v_*,v_*)$, $(v_1',v_2')$ and
$(v_{1*}',v_{2*}')$ in $V_{00}$.
As above, we are writing $v_k'$ for $v'(v_k,v_*,\sigma)$ and $v'_{k*}$ for
$v'_*(v_k,v_*,\sigma)$.
Second, at rate $2(\rho_t^1-\rho_t^2)^+(dv_*)(B(v_1-v_*,d\sigma
)\wedge
B(v_2-v_*,d\sigma))\,dt$, we replace $(v_1,v_2)$ by three particles
$(v_*,v_*)$, $(v_1',v_2')$ and $(v_{1*}',v_{2*}')$ in $V_{10}$ and one
particle $(v_1,v_2)$ in $V_{20}$.
Third, at rate $2(\rho_t^1\wedge\rho_t^2)(dv_*)(B(v_1-v_*,d\sigma
)-B(v_2-v_*,d\sigma))^+dt$, we replace $(v_1,v_2)$ by three particles
$v_*$, $v_1'$ and $v_{1*}'$ in $V_{01}$ and one particle $v_2$ in $V_{02}$.
Fourth, at rate $2(\rho_t^1-\rho_t^2)^+(B(v_1-v_*,d\sigma
)-B(v_2-v_*,d\sigma
))^+\,dt$, we replace $(v_1,v_2)$ by three particles
$v_*$, $v_1'$ and $v_{1*}'$ in $V_{11}$, one particle $v_2$ in $V_{12}$
and one particle $(v_1,v_2)$ in $V_{20}$.
The second and third transitions each have an obvious counterpart by
swapping $1$ and $2$, while there are three variants of the fourth transition
by swapping $1$ and $2$ in the environment or in the collision
intensity or in both.

On leaving $V_{00}$, either the coupling with respect to environment is
broken, or that with respect to the starting points.
This corresponds to transitions on the one hand to $V_{1k}$ or $V_{2k}$
for some $k$, or on the other hand to $V_{j1}$ or $V_{j2}$ for some
$j$, respectively.
Once the environment coupling is broken, a particle branches as in the
starting point coupling, while if the starting point coupling is
broken, a
particle branches as in the environment coupling.
Thus the transitions in $V_{jk}$ for $j=1,2$ are as described above for
$V_k$, while those in $V_{jk}$ for $k=1,2$ are as described above for
$\hat V_j$.

For $j,k=0,1,2$, write $(\Xi_t^{jk})_{t\ge0}$ for the empirical
distribution of particles in $V_{jk}$ when we initiate the branching process
just described with a single particle $(v_1,v_2)$ in $V_{00}$ at time $0$.
Write $q_{jk}$ for the bijection $V_{jk}\to V_k$.
For $k=1,2$, write $p_{jk}$ for the projection $(v_1,v_2,j)\mapsto
(v_k,j) \dvtx  V_{j0}\to\hat V_j$, and write $\hat q_{jk}$ for the bijection
$V_{jk}\to\hat V_j$.
Note that $\hat\pi_j\circ p_{jk}=p_k\circ q_{j0}$ on $V_{j0}$ and
$\hat\pi_j\circ\hat q_{jk}=\pi_k\circ q_{jk}$ on $V_{jk}$ for
$j=0,1,2$ and $k=1,2$.
For $j,k=1,2$, set
\[
\Gamma_t^{j0}=\Xi_t^{00}\circ
q_{00}^{-1}+\Xi_t^{j0}\circ
q_{j0}^{-1},\qquad \Gamma_t^{jk}=
\Xi_t^{0k}\circ q_{0k}^{-1}+
\Xi_t^{jk}\circ q_{jk}^{-1}
\]
and
\[
\hat\Gamma_t^{0k}=\Xi_t^{00}\circ
p_{0k}^{-1}+\Xi_t^{0k}\circ\hat
q_{0k}^{-1}, \qquad\hat\Gamma_t^{jk}=
\Xi_t^{j0}\circ p_{jk}^{-1}+
\Xi_t^{jk}\circ\hat q_{jk}^{-1},
\]
and set
\[
\Lambda_t^{jk}=\Gamma_t^{j0}\circ
\pi_k^{-1}+\Gamma_t^{jk}\circ
p_k^{-1}=\hat\Gamma _t^{0k}\circ\hat
\pi_0^{-1}+\hat\Gamma_t^{jk}\circ\hat
\pi_j^{-1}.
\]
It can\vspace*{1pt} be checked that $(\Lambda_t^{jk})_{t\ge0}$ is a copy of the
branching process
$(\Lambda_t)_{t\ge0}$ starting from $v_k$ at time $0$ in the environment
$(\rho_t^j)_{t\ge0}$.
Moreover $(\Gamma_t^{j0},\Gamma_t^{j1},\break \Gamma_t^{j2})_{t\ge0}$ is a
copy of the
starting point coupling in environment $(\rho_t^j)_{t\ge0}$,
and $(\hat\Gamma_t^{0k},\hat\Gamma_t^{1k},\hat\Gamma
_t^{2k})_{t\ge0}$ is a
copy of the environment coupling\vspace*{1pt} with starting point $v_k$.
As in the earlier constructions, we lift to processes $(\Xi
_t^{jk,*})_{t\ge0}$ in the signed spaces $V_{jk}^*=V_{jk}^-\cup
V_{jk}^+=V_{jk}\times\{-1,1\}$,
initiating with a particle\vspace*{-1,5pt} $(v_1,v_2)$ in $V_{00}^+$ and with the
`$v_*$' particles switching signs.
Then, for $j,k=1,2$ the associated process $(\Lambda_t^{jk,*})_{t\ge
0}$ in
$V^*$ is a linearized Kac process with environment $(\rho_t^j)_{t\ge
0}$ starting from $v_k$.

\begin{lemma}\label{DCL2}
For all $p\in[2,\infty)$, there is a constant $C(p)<\infty$ such that
\begin{eqnarray*}
&& E_{(0,v_1,v_2)}\bigl\langle1+\llvert v\rrvert ^p,
\Xi_t^{11}+\Xi_t^{12}+
\Xi_t^{21}+\Xi _t^{22}\bigr\rangle
\\
&&\qquad\le C(p)\kappa\bigl(1+\llvert v_1\rrvert ^{p+1}+\llvert
v_2\rrvert ^{p+1}\bigr)\llvert v_1-v_2
\rrvert \exp \biggl\{ C(p)\int_0^t\bar
m_{p+2}(s)\,ds \biggr\}\\
&&\qquad\quad{}\times\int_0^td_{p+1}(s)\,ds.
\end{eqnarray*}
\end{lemma}
\begin{pf}
It will suffice by symmetry to consider $\langle1+|v|^p,\Xi
_t^{11}\rangle$.
The process $\langle1+|v|^p,\Xi_t^{11}\rangle$ makes jumps due to uncoupling
transitions from $V_{01}$ and $V_{10}$ and also directly from $V_{00}$,
and it makes further jumps due to the branching of particles in $V_{11}$.
Jumps of $3+|v'_1|^p+|v_{1*}'|^p+|v_*|^p$ occur at rate
\begin{eqnarray*}
&& \!\!2B(v_1-v_*,d\sigma)\Xi_{t-}^{01}(dv_1)
\bigl(\rho^1_t-\rho _t^2
\bigr)^+(dv_*)\,dt
\\
&&\!\!\qquad {}+2\bigl(B(v_1-v_*,d\sigma)-B(v_2-v_*,d\sigma)\bigr)^+
\Xi _{t-}^{10}(dv_1,dv_2)
\rho^1_t(dv_*)\,dt
\\
&&\!\!\qquad{}+2\bigl(B(v_1-v_*,d\sigma)-B(v_2-v_*,d\sigma)\bigr)^+
\Xi _{t-}^{00}(dv_1,dv_2) \bigl(
\rho^1_t-\rho_t^2\bigr)^+(dv_*)\,dt.
\end{eqnarray*}
Jumps of $1+|v_1|^p$ occur at rate
\begin{eqnarray*}
&&\!\!2B(v_1-v_*,d\sigma)\Xi_{t-}^{01}(dv_1)
\bigl(\rho^1_t-\rho _t^2
\bigr)^-(dv_*)\,dt
\\
&&\!\!\qquad {}+2\bigl(B(v_1-v_*,d\sigma)-B(v_2-v_*,d\sigma)\bigr)^-
\Xi _{t-}^{10}(dv_1,dv_2)
\rho^1_t(dv_*)\,dt
\\
&&\!\!\qquad{}+2\bigl(B(v_1-v_*,d\sigma)-B(v_2-v_*,d\sigma)\bigr)^-
\Xi _{t-}^{00}(dv_1,dv_2) \bigl(
\rho^1_t-\rho_t^2\bigr)^+(dv_*)\,dt.
\end{eqnarray*}
Jumps of $2+|v'|^p+|v_*'|^p+|v_*|^p-|v|^p$ occur at rate
\[
2B(v-v_*,d\sigma)\Xi_{t-}^{11}(dv)\rho_t^1(dv_*)\,dt.
\]
Fix a starting point $(v_1,v_2)$ in $V_{00}$.
Recall that $\Xi_t^{00}$ and $\Xi_t^{10}$ are supported on pairs
$(u_1,u_2)$ with $|u_1-u_2|\le|v_1-v_2|$.
We use the inequalities
\begin{equation}
\label{TVPO} 3+\bigl\llvert v'\bigr\rrvert ^p+\bigl
\llvert v_*'\bigr\rrvert ^p+\llvert v_*\rrvert
^p\le C(p) \bigl(1+\llvert v\rrvert ^p+\llvert v_*\rrvert
^p\bigr)
\end{equation}
and
\[
\bigl(1+\llvert v\rrvert ^p\bigr)\llvert v-v_*\rrvert \le\bigl(1+
\llvert v\rrvert ^p+\llvert v_*\rrvert ^p\bigr)\llvert
v-v_*\rrvert \le C(p) \bigl(1+\llvert v\rrvert ^{p+1}\bigr) \bigl(1+
\llvert v_*\rrvert ^{p+1}\bigr)
\]
to see that the drift of $\langle1+|v|^p,\Xi_t^{11}\rangle$ due to uncoupling
transitions from $V_{01}$ is no greater than
\[
C(p)\bigl\langle1+\llvert v\rrvert ^{p+1},\Xi^{01}_{t-}
\bigr\rangle d_{p+1}(t).
\]
We use (\ref{TVPO}) and the inequalities
\[
1+\llvert v\rrvert ^p\le1+\llvert v\rrvert ^p+\llvert
v_*\rrvert ^p\le\bigl(1+\llvert v\rrvert ^p\bigr)
\bigl(1+\llvert v_*\rrvert ^p\bigr)
\]
to see that the drift of $\langle1+|v|^p,\Xi_t^{11}\rangle$ due to uncoupling
transitions from $V_{10}$ is no greater than
\[
C(p)\kappa\llvert v_1-v_2\rrvert \bigl\langle1+\llvert v
\rrvert ^p,\Xi^{10}_{t-}\circ
p_{11}^{-1}\bigr\rangle\bar m_p(t)
\]
while the drift of $\langle1+|v|^p,\Xi_t^{11}\rangle$ due to uncoupling
transitions from $V_{00}$ is no greater than
\[
C(p)\kappa\llvert v_1-v_2\rrvert \bigl\langle1+\llvert v
\rrvert ^p,\Xi^{00}_{t-}\circ
p_{01}^{-1}\bigr\rangle d_p(t).
\]
Finally, by (\ref{POF}) and (\ref{POE}), the drift of $\langle
1+|v|^p,\Xi
_t^{11}\rangle$ due to branching in $V_{11}$ is no greater than
\[
C(p)\bigl\langle1+\llvert v\rrvert ^p,\Xi^{11}_{t-}
\bigr\rangle\bar m_{p+1}(t).
\]
Set
\[
E_{p+2}(t)=C(p)\exp \biggl\{\int_0^tC(p)
\bar m_{p+2}(s)\,ds \biggr\},
\]
where $C(p)<\infty$ remains to be chosen.
By Lemma~\ref{GHU}, we can choose $C(p)$ so that
\begin{eqnarray*}
E_{(0,v_1,v_2)}\bigl\langle1+\llvert v\rrvert ^{p+1},
\Xi_t^{01}\bigr\rangle &\le &  E_{(0,v_1,v_2)}\bigl\langle 1+
\llvert v\rrvert ^{p+1},\Gamma_t^{11}\bigr\rangle
\\
&\le & \kappa t\bigl(1+\llvert v_1\rrvert ^{p+1}+\llvert
v_2\rrvert ^{p+1}\bigr)E_{p+2}(t)\llvert
v_1-v_2\rrvert.
\end{eqnarray*}
By Lemma~\ref{DCL}, we can choose $C(p)$ so that, moreover,
\begin{eqnarray*}
&& E_{(0,v_1,v_2)}\bigl\langle1+\llvert v\rrvert ^p,
\Xi^{10}_t\circ p_{11}^{-1}\bigr\rangle\\
&&\qquad\le   E_{(0,v_1)}\bigl\langle 1+\llvert v\rrvert ^p,\hat
\Gamma_t^{11}\bigr\rangle\\
&&\qquad\le  \bigl(1+\llvert v_1
\rrvert ^{p+1}\bigr)E_{p+2}(t)\int_0^td_{p+1}(s)\,ds.
\end{eqnarray*}
By Proposition~\ref{JPE}, we can choose $C(p)$ so that, moreover,
\[
E_{(0,v_1,v_2)}\bigl\langle1+\llvert v\rrvert ^p,
\Xi^{00}_t\circ p_{01}^{-1}\bigr\rangle
\le E_{(0,v_1)}\bigl\langle 1+\llvert v\rrvert ^p,
\Lambda_t^{11}\bigr\rangle\le\bigl(1+\llvert v_1
\rrvert ^p\bigr)E_{p+2}(t).
\]
Set $g(t)=E_{(0,v_1,v_2)}\langle1+|v|^p,\Xi_t^{11}\rangle$.
The three estimates just obtained give us control of the expected drift
of $\langle1+|v|^p,\Xi_t^{11}\rangle$, so we obtain a constant
$C(p)<\infty$
such that
\begin{eqnarray*}
g(t) &\le &  C(p)\kappa\bigl(1+\llvert v_1\rrvert ^{p+1}+\llvert
v_2\rrvert ^{p+1}\bigr)\llvert v_1-v_2
\rrvert E_{p+2}(t)\int_0^td_{p+1}(s)\,ds\\
&&{}+C(p)
\int_0^t\bar m_{p+1}(s)g(s)\,ds,
\end{eqnarray*}
which gives the claimed inequality by Gronwall's lemma.
\end{pf}

\begin{pf*}{Proof of Proposition~\ref{FSG}}
It will suffice to consider the case where $s=0$.
Set $f_0^j=E_{0t}^jf$.
We use the coupling of linearized Kac processes for environments $(\rho
_t^1)_{t\ge0}$ and $(\rho_t^2)_{t\ge0}$ described above.
By Lemma~\ref{DCL},
\begin{eqnarray}
\nonumber
&& \bigl\llvert f_0^1(v)-f_0^2(v)
\bigr\rrvert\\
\label{DEEE}
&&\qquad =\bigl\llvert E_{(0,v)}\bigl\langle f,\tilde
\Lambda_t^1-\tilde\Lambda_t^2\bigr
\rangle\bigr\rrvert
 \le  E_{(0,v)}\bigl\langle\llvert f\rrvert ,\hat
\Gamma_t^1+\hat\Gamma_t^2\bigr
\rangle \le E_{(0,v)}\bigl\langle1+\llvert v\rrvert ^p,\hat
\Gamma_t^1+\hat\Gamma_t^2\bigr
\rangle
\\
\nonumber
 &&\qquad \le   C(p) \bigl(1+\llvert v\rrvert ^{p+1}\bigr)\exp
\biggl\{\int_0^tC(p)\bar m_{p+2}(s)\,ds
\biggr\}\int_0^td_{p+1}(s)\,ds.
\end{eqnarray}
Assume now that $f\in\mathcal{F}(p)$.
Then
\[
\bigl\llvert f(v_1)-f(v_2)\bigr\rrvert \le p\bigl(1+
\llvert v_1\rrvert ^p+\llvert v_2\rrvert
^p\bigr)\llvert v_1-v_2\rrvert
\]
for all $v_1,v_2\in\mathbb{R}^d$.
On the other hand, if $g$ satisfies this inequality, together with
$|g(v)|\le1+|v|^p$, then $\|g\|_{(p)}\le3p$.
We now use the coupling of four linearized Kac processes for the two
environments and two starting points $v_1$ and $v_2$.
For $j=1,2$, the measure $\Xi_t^{j0,*}$ is supported on pairs
$(u_1,u_2)$ with $|u_1-u_2|\le|v_1-v_2|$.
So
\begin{eqnarray*}
\bigl\langle f\circ p_1-f\circ p_2,\Xi_t^{j0,*}
\bigr\rangle &\le &  p\llvert v_1-v_2\rrvert \bigl\langle 1+
\llvert u_1\rrvert ^p+\llvert u_2\rrvert
^p,\Xi_t^{j0}\bigr\rangle \\
&\le &  p\llvert
v_1-v_2\rrvert \bigl\langle1+\llvert v\rrvert
^p,\hat \Gamma _t^{j1}+\hat
\Gamma_t^{j2}\bigr\rangle.
\end{eqnarray*}
By Lemmas \ref{DCL} and \ref{DCL2},
\begin{eqnarray*}
&&\bigl(f_0^1-f_0^2\bigr)
(v_1)-\bigl(f_0^1-f_0^2
\bigr) (v_2)\\
&&\qquad=E_{(0,v_1,v_2)}\bigl\langle f,\tilde \Lambda
_t^{11}-\tilde\Lambda_t^{12}-\tilde
\Lambda_t^{21}+\tilde\Lambda _t^{22}
\bigr\rangle
\\
&&\qquad =E_{(0,v_1,v_2)}\bigl(\bigl\langle f,\Xi_t^{11,*}-
\Xi_t^{12,*}-\Xi _t^{21,*}+
\Xi_t^{22,*}\bigr\rangle\\
&&\hspace*{54pt}\qquad{}+\bigl\langle f\circ
p_1-f_2\circ p_2,\Xi _t^{10,*}-
\Xi _t^{20,*}\bigr\rangle\bigr)
\\
&&\qquad \le E_{(0,v_1,v_2)}\bigl(\bigl\langle1+\llvert v\rrvert ^p,
\Xi_t^{11}+\Xi_t^{12}+\Xi
_t^{21}+\Xi_t^{22}\bigr\rangle\\
&&\qquad\qquad\hspace*{32pt}{}+p
\llvert v_1-v_2\rrvert \bigl\langle1+\llvert v\rrvert
^p,\hat\Gamma _t^{11}+\hat\Gamma
_t^{12}+\hat\Gamma_t^{21}+\hat
\Gamma_t^{22}\bigr\rangle\bigr)
\\
&&\qquad \le C(p)\kappa\bigl(1+\llvert v_1\rrvert ^{p+1}+\llvert
v_2\rrvert ^{p+1}\bigr)\\
&&\qquad\quad{}\times\llvert v_1-v_2
\rrvert \exp \biggl\{ \int_0^tC(p)\bar
m_{p+2}(s)\,ds \biggr\}\int_0^td_{p+1}(s)\,ds.
\end{eqnarray*}
Here, there are no terms in $\Xi_t^{0k,*}$ for $k=0,1,2$ because these
are the empirical distributions of particles, or pairs of particles,
with unbroken environment coupling, which cancel completely in the
considered integral.
On combining this estimate with (\ref{DEEE}), we deduce that
\[
\bigl\llVert f_0^1-f_0^2\bigr
\rrVert _{(p+1)}\le3C(p)\kappa\exp \biggl\{\int_0^tC(p)
\bar m_{p+2}(s)\,ds \biggr\}\int_0^td_{p+1}(s)\,ds.
\]
\upqed\end{pf*}

\section{Maximal inequalities for stochastic convolutions}\label{SCB}
The key formula for our analysis is shown in Proposition~\ref{FMU}.
For all $t\ge0$ and all functions $f$ in our weighted Lipschitz class
$\mathcal{F}$, we have
\[
\bigl\langle f,\mu_t^N-\mu_t^{N'}
\bigr\rangle=\bigl\langle f_{0t},\mu_0^N-\mu
_0^{N'}\bigr\rangle+\int_0^t
\bigl\langle f_{st},dM_s^N\bigr\rangle-\int
_0^t\bigl\langle f_{st},dM^{N'}_s
\bigr\rangle,
\]
where
\[
f_{st}(v)=E_{st}f(v)=E_{(s,v)}\langle f,\tilde
\Lambda_t\rangle
\]
and where $(\Lambda_t^*)_{t\ge s}$ is the linearized Kac process in
environment $((\mu_t^N+\mu_t^{N'})/ 2)_{t\ge0}$ starting from $v$ at
time $s$.

The notion of stochastic convolution has been extensively studied in
connection with infinite-dimensional stochastic evolution equations;
see, for example, \cite{MR1154532,MR1864042}.
The operator
\[
f\mapsto\int_0^t\bigl\langle
E_{st}f,dM_s^N\bigr\rangle
\]
shares some features with stochastic convolutions, namely that
$E_{st}E_{tu}=E_{su}$ for $s\le t\le u$ and that good estimates rely on
exploiting martingale properties of the integrator.
In this section, we prove a maximal inequality for this operator in
Wasserstein norms, in the case where the environment $((\mu_t^N+\mu
_t^{N'})/2)_{t\ge0}$ is replaced by a nonrandom process $(\rho
_t)_{t\ge0}$.
The proof of Proposition~\ref{DAT} below uses some of the same ideas
in a simpler context.

We will use the following inequality for a function $f$ on $\mathbb{R}^d$
which is Lipschitz of constant $1$.
For $B=[0,2^{-k}]^d$, we have
\begin{equation}
\label{LIPA} \bigl\llvert f(v)-\langle f\rangle_B\bigr\rrvert
\le2^{-k}c_d,\qquad v\in B,
\end{equation}
where $\langle f\rangle_B$ is the average value of $f$ on $B$ and
where $c_d=\mathbb{E}
|X|$, with $X$ uniformly distributed on $[0,1]^d$.
To see this, set $Y=2^{-k}X$, and note that $|f(v)-f(Y)|\le|v-Y|\le
|Y|$ so $|f(v)-\langle f\rangle_B|=|\mathbb{E}(f(v)-f(Y))|\le\mathbb
{E}|Y|=2^{-k}\mathbb{E}|X|$.
By a similar calculation, we have also
\begin{equation}
\label{LIPB} \bigl\llvert \langle f\rangle_B-\langle f
\rangle_{2B}\bigr\rrvert \le2^{-k}c_d.
\end{equation}
It is the scaling properties of inequalities (\ref{LIPA}) and (\ref
{LIPB}) which will be critical for our argument, rather than the value
of the constant $c_d$.

Let $(\rho_t)_{t\ge0}$ be a nonrandom process\footnote{We will in
fact use only the case where $(\rho_t)_{t\ge0}$ is constant.}
satisfying (\ref{SKP}).
Write, as above, $m_p(t)=\langle1+|v|^p,\rho_t\rangle$, and set
\[
m^*(p)=\sup_{t\ge0}m_p(t).
\]
Then, by Proposition~\ref{JPE}, for all $s\ge0$ and all $v_0\in
\mathbb{R}^d$,
and for the linearized Kac process $(\Lambda_t^*)_{t\ge s}$ in environment
$(\rho_t)_{t\ge0}$ starting from $v_0$ at time $s$, we have
\[
E_{(s,v_0)}\bigl\langle1+\llvert v\rrvert ^p,
\Lambda_t\bigr\rangle\le \bigl(1+\llvert v_0\rrvert
^p\bigr)e^{c(p)m^*(p+1)(t-s)}.
\]
Thus, whenever $m^*(p+1)<\infty$, we can define, for $s\le t$ and
$f\in\mathcal{F}(p)$,
\[
f_{st}(v)=E_{st}f(v)=E_{(s,v)}\langle f,\tilde
\Lambda_t\rangle.
\]
\begin{proposition}\label{WMME}
For all $d\ge3$, $p\in[2,\infty)$ and all $\delta\in(0,1]$, there
is a
constant $C(d,\delta, p)<\infty$ such that, for all $T\in[0,\infty)$,
we have
\begin{eqnarray}
&& \biggl\llVert \sup_{t\le T}\sup
_{f\in\mathcal{F}(p)}\int_0^t\bigl\langle
f_{st},dM^N_s\bigr\rangle \biggr\rrVert
_2
\nonumber
\\[-8pt]
\label{FETA}
\\[-8pt]
\nonumber
&&\qquad\le C\kappa N^{-1/d}e^{Cm^*(p+3+\delta)T} \biggl(\mathbb{E}\int
_0^T\bigl\langle\llvert v\rrvert
^{2p+5+2\delta
},\mu_s^N\bigr\rangle \,ds
\biggr)^{1/2}.
\end{eqnarray}
The same inequality holds for $d=2$ if we replace $N^{-1/d}$ by
$N^{-1/2}\log N$.
\end{proposition}
Here we have written $\|\cdot\|_2$ for the norm in $L^2(\mathbb{P})$.
This estimate will be applied in the next section, using the moment
estimates derived in Section~\ref{MOM} to control the right-hand side.
We will use also the following comparison estimate for two nonrandom
processes $(\rho_t^1)_{t\ge0}$ and $(\rho_t^2)_{t\ge0}$ satisfying
(\ref{SKP}).
Fix $p\in[2,\infty)$.
Write
\[
\bar m^*(p)=\sup_{t\ge0}\bigl\langle1+\llvert v\rrvert
^p,\rho_t^1+\rho_t^2
\bigr\rangle.
\]
We assume that $\bar m^*(p+1)<\infty$.
For $j=1,2$ and $f\in\mathcal{F}(p)$, define for $s,t\ge0$ with
$s\le t$
\[
f_{st}^j(v)=E_{st}^jf(v)=E_{(s,v)}
\bigl\langle f,\tilde\Lambda^j_t\bigr\rangle,
\]
where $(\Lambda_t^{j,*})_{t\ge s}$ is a linearized Kac process in
environment $(\rho^j_t)_{t\ge0}$ starting from~$v$.
\begin{proposition}\label{WMMF}
For all $d\ge3$, $p\in[2,\infty)$ and all $\delta\in(0,1]$, there
is a
constant $C(d,\delta, p)<\infty$ such that, for all $T\in[0,\infty)$,
we have
\begin{eqnarray}
&&\biggl\llVert \sup_{t\le T}\sup_{f\in\mathcal{F}(p)}
\int_0^t\bigl\langle f_{st}^1-f_{st}^2,dM^N_s
\bigr\rangle\biggr\rrVert _2
\nonumber
\\
\label{FETB}
&&\qquad \le C\kappa N^{-1/d}\bar m^*(p+2+\delta)\\
\nonumber
&&\qquad\quad{}\times Te^{C\bar
m^*(p+3+\delta
)T} \biggl(
\mathbb{E}\int_0^T\bigl\langle\llvert v\rrvert
^{2p+5+2\delta},\mu _s^N\bigr\rangle \,ds
\biggr)^{1/2}.
\end{eqnarray}
The same inequality holds for $d=2$ if we replace $N^{-1/d}$ by
$N^{-1/2}\log N$.
\end{proposition}
A small variation of the following proofs would allow the insertion of
a factor of $d^*(p+3+\delta)$ on the right in (\ref{FETB}),
where $d^*(p)=\sup_{t\ge0}\langle1+|v|^p,|\rho_t^1-\rho
_t^2|\rangle$, at the
cost of replacing $p$ in all other terms on the right by $p+1$.
We omit details as this variation is not needed for our main result.

\begin{pf*}{Proof of Proposition~\ref{WMME}}
Assume for now that $d\ge3$.
It will suffice to consider the case where $N\ge2^{2d}$.
We first prove a simpler estimate, where the function $f_{st}$ is
replaced by $f_{sT}$ on the left-hand side.
Set $L=\lfloor\log_2N/d\rfloor$, and note that $L\ge2$.
For $k\in\mathbb{Z}$, set $B_k=(-2^k,2^k]^d$.
Set $A_0=B_0$, and for $k\ge1$, set $A_k=B_k\setminus B_{k-1}$.
For $k\ge1$ and any integer $\ell\ge2$, there is a unique way to
partition $A_k$ by a set $\mathcal{P}_{k,\ell}$ of $2^{\ell
d}-2^{(\ell
-1)d}$ translates of $B_{k-\ell}$.
Also, there is a unique way to partition $A_0$ by a set $\mathcal
{P}_{0,\ell
}$ of $2^{\ell d}$ translates of $B_{-\ell}$.
Fix $p\in[2,\infty)$ and $f\in\mathcal{F}(p)$.
Then, for all $v,v'\in\mathbb{R}^d$, we have
\[
\bigl\llvert \hat f^{(p)}(v)\bigr\rrvert =\bigl\llvert f(v)/
\bigl(1+|v|^p\bigr)\bigr\rrvert \le1,\qquad  \bigl\llvert \hat
f^{(p)}(v)-\hat f^{(p)}\bigl(v'\bigr)\bigr\rrvert
\le\bigl\llvert v-v'\bigr\rrvert.
\]
For $B\in\mathcal{P}_{k,2}$, set $a_B=\langle\hat f^{(p)}\rangle
_B$, and note that
$|a_B|\le1$.
For $\ell\ge3$ and $B\in\mathcal{P}_{k,\ell}$, set $a_B=\langle
\hat f^{(p)}\rangle
_B-\langle\hat f^{(p)}\rangle_{\pi(B)}$,
where $\pi(B)$ is the unique element of $\mathcal{P}_{k,\ell-1}$ containing
$B$, and note that $|a_B|\le2^{k-\ell+1}c_d$.
Set $c_d'=4\vee(2c_d)$, then $|a_B|\le2^{k-\ell}c_d'$ for all $B\in
\mathcal{P}_{k,\ell}$, for all $k\ge0$ and all $\ell\ge2$.
Fix $\delta\in(0,1]$, and for $B\in\mathcal{P}_{k,\ell}$, set
\[
h^B(v)=2^{(1+\delta)k}\bigl(1+\llvert v\rrvert ^p
\bigr)1_B(v).
\]
Define a function $g_k$, supported on $A_k$, by
\[
f1_{A_k}=\sum_{\ell=2}^L\sum
_{B\in\mathcal{P}_{k,\ell
}}a_B\bigl(1+\llvert v\rrvert
^p\bigr)1_B(v)+g_k.
\]
Fix $K\in\mathbb{N}$, and set $g=\sum_{k=0}^Kg_k+f1_{B_K^c}$.
Note that $\hat g_k^{(p)}=\hat f^{(p)}-\langle\hat f^{(p)}\rangle_B$
on $B$ for
all $B\in\mathcal{P}_{k,L}$.
For $v\in A_k$, we have $|v|\ge2^{k-1}$, so
\[
\bigl\llvert \hat g_k^{(p)}(v)\bigr\rrvert
\le2^{k+1-L}c_d\le2^{-L+2}c_d\bigl(1+
\llvert v\rrvert \bigr).
\]
For $v\in B_K^c$, we have $|v|\ge2^{K-1}$, so $|\hat f^{(p)}(v)|\le
2^{-K+1}|v|$.
Set $c_d''=(8c_d)\vee4$.
Then, for all $v\in\mathbb{R}^d$, we have
\[
\bigl\llvert g(v)\bigr\rrvert \le\bigl\{2^{-L+2}c_d\bigl(1+
\llvert v\rrvert \bigr)+2^{-K+1}\llvert v\rrvert \bigr\}\bigl(1+\llvert v
\rrvert ^p\bigr)\le \bigl(2^{-K}+2^{-L}
\bigr)c_d''\bigl(1+\llvert v\rrvert
^{p+1}\bigr).
\]
Now
\[
f=\sum_{\ell=2}^L\sum
_{k=0}^K\sum_{B\in\mathcal{P}_{k,\ell
}}2^{-(1+\delta)k}a_Bh^B+g
\]
so
\begin{eqnarray}
&& \int_0^t\bigl\langle
f_{sT},dM^N_s\bigr\rangle
\nonumber
\\[-8pt]
\label{FEE}
\\[-8pt]
\nonumber
&&\qquad=\sum
_{\ell=2}^L\sum_{k=0}^K
\sum_{B\in
\mathcal{P}_{k,\ell}}2^{-(1+\delta)k}a_B\int
_0^t\bigl\langle h^B_{sT},dM^N_s
\bigr\rangle+\int_0^t\bigl\langle
g_{sT},dM^N_s\bigr\rangle,
\end{eqnarray}
where $h_{sT}^B=E_{sT}h^B$ and $g_{sT}=E_{sT}g$.
It will be convenient to set
\[
E(p)=\exp\bigl\{c(p)m^*(p+1)T\bigr\}
\]
and
\[
c(d,\delta)=\bigl(1-2^{-2\delta}\bigr)^{-1/2}c_d',\qquad
A=\bigl(2^{-K}+2^{-L}\bigr)c_d''E(p+1).
\]
Note that $2^{-(1+\delta)k}|a_B|\le2^{-\ell-\delta k}c_d'$ for all
$B\in\mathcal{P}
_{k,\ell}$, and $\mathcal{P}_{k,\ell}$ has cardinality at most
$2^{d\ell}$, so
\[
\sum_{k=0}^K\sum
_{B\in\mathcal{P}_{k,\ell}}\bigl(2^{-(1+\delta
)k}a_B
\bigr)^2\le 2^{(d-2)\ell}c(d,\delta)^2.
\]
Also, by Proposition~\ref{JPE}, for all $s\in[0,T]$, we have
\[
\bigl\llvert g_{sT}(v)\bigr\rrvert \le A\bigl(1+\llvert v\rrvert
^{p+1}\bigr).
\]
We use Cauchy--Schwarz in (\ref{FEE}) to obtain
\begin{eqnarray}
&&\sup_{t\le T}\sup_{f\in\mathcal{F}(p)}\int
_0^t\bigl\langle f_{sT},dM^N_s
\bigr\rangle
\nonumber
\\
\label{PEE}
&&\qquad \le\sum_{\ell=2}^L2^{(d/2-1)\ell}c(d,
\delta) \Biggl(\sum_{k=0}^K\sum
_{B\in\mathcal{P}_{k,\ell}} \sup_{t\le T}\biggl\llvert \int
_0^t\bigl\langle h^B_{sT},dM^N_s
\bigr\rangle\biggr\rrvert ^2 \Biggr)^{1/2}\\
\nonumber
&&\qquad\quad{}+A\int
_0^T\bigl\langle1+\llvert v\rrvert
^{p+1},\bigl\llvert dM^N_s\bigr\rrvert \bigr
\rangle.
\end{eqnarray}
Set
\begin{eqnarray*}
h(v)&=&\bigl(1+\llvert v\rrvert ^p\bigr)\sum
_{k=0}^K2^{(1+\delta)k}1_{A_k}(v)\\
&=&\sum
_{k=0}^K\sum
_{B\in\mathcal{P}_{k,\ell}}h^B(v).
\end{eqnarray*}
Note that $1\vee|v|\ge2^{k-1}$ for all $v\in A_k$ and all $k$.
Set $q=p+1+\delta$ and $A'=8E(q)$.
Then
\[
h(v)\le2^{1+\delta}\bigl(1\vee\llvert v\rrvert \bigr)^{1+\delta}
\bigl(1+\llvert v\rrvert ^p\bigr)\le8\bigl(1+\llvert v\rrvert
^q\bigr),
\]
so by Proposition~\ref{JPE},
\[
E_{(s,v)}\langle h,\Lambda_T\rangle\le A'
\bigl(1+\llvert v\rrvert ^q\bigr).
\]
Note that $|h_{sT}^B(v)|\le E_{(s,v)}\langle h^B,\Lambda_T\rangle$, so
\[
\sum_{k=0}^K\sum
_{B\in\mathcal{P}_{k,\ell}}\bigl(h_{sT}^B(v)
\bigr)^2\le \bigl(E_{(s,v)}\langle h,\Lambda_T
\rangle\bigr)^2\le\bigl(A'\bigr)^2\bigl(1+
\llvert v\rrvert ^q\bigr)^2.
\]
Hence, for some constant $C(q)<\infty$, we have
\begin{eqnarray*}
&&\sum_{k=0}^K\sum
_{B\in\mathcal{P}_{k,\ell}}\bigl\{ h^B_{sT}
\bigl(v'\bigr)+h^B_{sT}\bigl(v'_*
\bigr)-h_{sT}^B(v)-h_{sT}^B(v_*)
\bigr\}^2
\\
&&\qquad \le4\sum_{k=0}^K\sum
_{B\in\mathcal{P}_{k,\ell}}\bigl\{ h^B_{sT}
\bigl(v'\bigr)^2+h^B_{sT}
\bigl(v'_*\bigr)^2+h_{sT}^B(v)^2+h_{sT}^B(v_*)^2
\bigr\}
\\
&&\qquad \le4\bigl(A'\bigr)^2E(q)^2\bigl\{
\bigl(1+\bigl\llvert v'\bigr\rrvert ^q
\bigr)^2+\bigl(1+\bigl\llvert v'_*\bigr\rrvert
^q\bigr)^2+\bigl(1+\llvert v\rrvert ^q
\bigr)^2+\bigl(1+\llvert v_*\rrvert ^q\bigr)^2
\bigr\}
\\
&&\qquad \le C(q) \bigl(A'\bigr)^2\bigl(1+\llvert v\rrvert
^{2q}+\llvert v_*\rrvert ^{2q}\bigr).
\end{eqnarray*}
Then, by Doob's $L^2$-inequality,
\begin{eqnarray*}
&&\sum_{k=0}^K\sum
_{B\in\mathcal{P}_{k,\ell}}\mathbb{E} \biggl(\sup_{t\le T}\biggl
\llvert \int_0^t\bigl\langle
h^B_{sT},dM^N_s\bigr\rangle
\biggr\rrvert ^2 \biggr)
\\
&&\qquad=\frac{4}N\sum_{k=0}^K\sum
_{B\in\mathcal{P}_{k,\ell}} \mathbb{E}\int_0^T\!\!
\int\bigl\{ h^B_{sT}\bigl(v'
\bigr)+h^B_{sT}\bigl(v'_*
\bigr)-h_{sT}^B(v)-h_{sT}^B(v_*)
\bigr\} ^2\\
&&\qquad\quad\hspace*{90pt}{}\times B(v-v_*,d\sigma )\mu_{s-}^N(dv)
\mu_{s-}^N(dv_*)\,ds
\\
&&\qquad\le\frac{C(q)(A')^2}N\mathbb{E}\int_0^T\!\!\int
\bigl\{ 1+\llvert v\rrvert ^{2q}+\llvert v_*\rrvert ^{2q}
\bigr\} \llvert v-v_*\rrvert \mu_{s-}^N(dv)
\mu_{s-}^N(dv_*)\,ds
\\
&&\qquad\le\frac{C(q)(A')^2}N\mathbb{E}\int_0^T\bigl
\langle\llvert v\rrvert ^{2q+1},\mu ^N_s\bigr
\rangle \,ds.
\end{eqnarray*}
On the other hand, we have
\begin{eqnarray*}
&& \int_0^T\bigl\langle1+\llvert v\rrvert
^{p+1},\bigl\llvert dM^N_s\bigr\rrvert \bigr
\rangle \\
&& \qquad\le\int_0^T\!\!\int\bigl\{4+\bigl\llvert
v'\bigr\rrvert ^{p+1}+\bigl\llvert v_*'\bigr
\rrvert ^{p+1}+\llvert v\rrvert ^{p+1}+\llvert v_*\rrvert
^{p+1}\bigr\} \\
&&\qquad\qquad\hspace*{3pt}\quad{}\times(m+\bar m) (dv,dv_*,d\sigma,ds),
\end{eqnarray*}
where the measures $m$ and $\bar m$ are as defined in Section~\ref{MARTK}.
We split the integral using $m+\bar m=(m-\bar m)+2\bar m$ and use the
$L^2$-isometry for integrals with respect to the compensated measure
$m-\bar m$ to obtain
\begin{eqnarray}
\nonumber
&&\mathbb{E} \biggl(\biggl\llvert \int_0^T
\bigl\langle1+\llvert v\rrvert ^{p+1},\bigl\llvert dM^N_s
\bigr\rrvert \bigr\rangle \biggr\rrvert ^2 \biggr)
\\
&&\qquad \le C(p)\mathbb{E}\int_0^T\!\!\int
\bigl\{ 1+\llvert v\rrvert ^{p+1}+\llvert v_*\rrvert ^{p+1}\bigr
\} ^2\,d\bar m\nonumber\\
\label{GEE}
&&\qquad\quad {}+C(p)\mathbb{E} \biggl(\biggl\llvert \int
_0^T\!\!\int\bigl\{1+\llvert v\rrvert
^{p+1}+\llvert v_*\rrvert ^{p+1}\bigr\} \,d\bar m\biggr\rrvert
^2 \biggr)
\\
\nonumber
&&\qquad \le C(p)\mathbb{E}\int_0^T\bigl
\langle\llvert v\rrvert ^{2p+3},\mu ^N_s\bigr
\rangle \,ds+C(p)\mathbb{E} \biggl\llvert \int_0^T
\bigl\langle\llvert v\rrvert ^{p+2},\mu^N_s
\bigr\rangle \,ds\biggr\rrvert ^2
\\
\nonumber
 &&\qquad \le C(p) (1+T)\mathbb{E}\int_0^T
\bigl\langle\llvert v\rrvert ^{2p+3},\mu ^N_s
\bigr\rangle \,ds,
\end{eqnarray}
where the constant $C(p)<\infty$ varies from line to line.
In the final inequality, we dealt with the second term on the right by
writing $|v|^{p+2}=|v||v|^{p+1}$,
applying Cauchy--Schwarz and then using the fact that $\langle1,\mu
_t^N\rangle=\langle
|v|^2,\mu_t^N\rangle=1$.
We take $L^2$-norms in (\ref{PEE}) to obtain
\begin{eqnarray*}
&&\biggl\llVert \sup_{t\le T}\sup_{f\in\mathcal{F}(p)}\int
_0^t\bigl\langle f_{sT},dM^N_s
\bigr\rangle \biggr\rrVert _2
\\
&&\qquad \le \Biggl(\sum_{\ell=2}^L2^{(d/2-1)\ell}c(d,
\delta )A' \biggl(\frac{C(q)}N \biggr)^{1/2}+A
\bigl(C(p) (1+T)\bigr)^{1/2} \Biggr) \\
&&\qquad\quad{}\times\biggl(\mathbb{E}\int
_0^T\bigl\langle\llvert v\rrvert
^{2q+1},\mu_s^N\bigr\rangle \,ds
\biggr)^{1/2}.
\end{eqnarray*}
Recall that $A=(2^{-K}+2^{-L})c_d''E(p+1)$, $A'=8E(q)$, $L=\lfloor\log
_2N/d\rfloor$ and $q=p+1+\delta$.
Note that $2^{(d/2-1)L}N^{-1/2}\le N^{-1/d}$ and $2^{-L}\le2N^{-1/d}$.
Hence, on letting $K\to\infty$, we deduce that, for some constant
$C(d,\delta,p)<\infty$,
\begin{eqnarray}
&& \biggl\llVert \sup_{t\le T}\sup
_{f\in\mathcal{F}(p)}\int_0^t\bigl\langle
f_{sT},dM^N_s\bigr\rangle \biggr\rrVert
_2
\nonumber
\\[-8pt]
\label{FET}
\\[-8pt]
\nonumber
&&\qquad\le CN^{-1/d}e^{Cm^*(p+2+\delta)T} \biggl(\mathbb{E}\int
_0^T\bigl\langle \llvert v\rrvert
^{2p+3+2\delta},\mu _s^N\bigr\rangle \,ds
\biggr)^{1/2}.
\end{eqnarray}
This is not the inequality (\ref{FETA}) we seek because $f_{sT}$
rather than $f_{st}$ appears on the right-hand side.
However, it will prove to be a useful first step.

We now turn to the proof of (\ref{FETA}).
It will suffice to deal with the case where $T=2^{-J_0}$ for some
$J_0\in\mathbb{Z}$.
Set $\tau_j(t)=2^{-j}\lceil2^jt\rceil$.
Then, for all $t\in(0,T]$, we have $\tau_{J_0}(t)=T$ so, for $J\ge J_0$
and $s\in[0,t]$,
\[
f_{st}=f_{sT}+\sum_{j=J_0+1}^J(f_{s\tau_j(t)}-f_{s\tau
_{j-1}(t)})+(f_{st}-f_{s\tau_J(t)})
\]
and hence
\begin{eqnarray}
\nonumber
&&\biggl\llVert \sup_{t\le T}\sup_{f\in\mathcal{F}(p)}
\int_0^t\bigl\langle f_{st},dM^N_s
\bigr\rangle \biggr\rrVert _2\\
&&\qquad\le\biggl\llVert \sup
_{t\le T}\sup_{f\in\mathcal{F}(p)}\int_0^t
\bigl\langle f_{sT},dM^N_s\bigr\rangle\biggr
\rrVert _2
\nonumber
\\[-8pt]
\label{HEE}
\\[-8pt]
\nonumber
&&\qquad\quad{}+\sum_{j=J_0+1}^J\biggl\llVert \sup
_{t\le T}\sup_{f\in\mathcal
{F}(p)}\int_0^t
\bigl\langle f_{s\tau_j(t)}-f_{s\tau_{j-1}(t)},dM^N_s
\bigr\rangle\biggr\rrVert _2\\
\nonumber
&&\qquad\quad{}+\biggl\llVert \sup_{t\le
T}
\sup_{f\in\mathcal{F}(p)}\int_0^t\bigl
\langle f_{st}-f_{s\tau
_J(t)},dM^N_s\bigr
\rangle\biggr\rrVert _2.
\end{eqnarray}
Fix $j\ge J_0+1$ and, for $i=0,1,\dots,2^jT$, set $t_i=i2^{-j}$.
Note that, for $t\in(t_{i-1},t_i]$, we have
\[
f_{s\tau_j(t)}-f_{s\tau_{j-1}(t)}=
\cases{f_{st_i}-f_{st_{i+1}},&
\quad$\mbox{if $i$ is odd}$,\vspace*{3pt}
\cr
0,& \quad$\mbox{if $i$ is even}$.}
\]
Set $f^{(i)}=f-f_{t_it_{i+1}}$.
We can take $s=t_i$, $t=t_{i+1}$ and $\rho_r^1=\rho_r$, $\rho^2_r=0$
for all $r$ in Proposition~\ref{FSG} to obtain
\[
\bigl\llVert f^{(i)}\bigr\rrVert _{(p+1)}\le
A''2^{-j},
\]
where
\[
A''=C(p)\kappa e^{C(p)m^*(p+2)T}m^*(p+1).
\]
For $s\in[0,t_i]$, set $f^{(i)}_s=E_{st_i}f^{(i)}=f_{st_i}-f_{st_{i+1}}$.
Write
\[
X_{ij}=\sup_{t\in(t_{i-1},t_i]}\sup_{f\in\mathcal{F}(p)}\int
_0^t\bigl\langle f^{(i)}_s,dM^N_s
\bigr\rangle
\]
and note that
\[
\sup_{t\le T}\sup_{f\in\mathcal{F}(p)}\int
_0^t\bigl\langle f_{s\tau
_j(t)}-f_{s\tau
_{j-1}(t)},dM^N_s
\bigr\rangle\le\sup_{i\le2^jT}X_{ij}.
\]
Set
\[
A'''=CA''N^{-1/d}e^{Cm^*(p+3+\delta)T}
\biggl(\mathbb{E}\int_0^T\bigl\langle\llvert v
\rrvert ^{2p+5+2\delta
},\mu_s^N\bigr\rangle \,ds
\biggr)^{1/2},
\]
where $C$ is the constant $C(d,\delta,p+1)$ from (\ref{FET}).
We replace $T$ by $t_i$, $p$ by $p+1$ and $f$ by $f^{(i)}/\|f^{(i)}\|
_{(p+1)}$ in (\ref{FET}) to see that
\[
\llVert X_{ij}\rrVert _2\le2^{-j}A'''.
\]
Then
\begin{equation}
\label{LTI} \sum_{j=J_0+1}^J\Bigl\llVert
\sup_{i\le2^jT}X_{ij}\Bigr\rrVert _2 \le\sum
_{j=J_0+1}^J\bigl(2^jT
\bigr)^{1/2}2^{-j}A'''
\le3TA''',
\end{equation}
so there is a constant $C(d,\delta,p)<\infty$ such that
\begin{eqnarray}
\nonumber
&&\qquad\sum_{j=J_0+1}^J\biggl\llVert
\sup_{t\le T}\sup_{f\in\mathcal
{F}(p)}\int
_0^t\bigl\langle f_{s\tau_j(t)}-f_{s\tau_{j-1}(t)},dM^N_s
\bigr\rangle\biggr\rrVert _2
\nonumber
\\[-8pt]
\label{DEE}
\\[-8pt]
\nonumber
&&\qquad\qquad \le Cm^*(p+1)T\kappa N^{-1/d}e^{Cm^*(p+3+\delta)T} \biggl(\mathbb{E}\int
_0^T\bigl\langle \llvert v\rrvert
^{2p+5+2\delta},\mu_s^N\bigr\rangle \,ds
\biggr)^{1/2}.
\end{eqnarray}
Finally, we can take $t=\tau_J(t)$, $\rho_r^1=\rho_r$ and $\rho
^2_r=\rho_r1_{\{r\le t\}}$ for all $r$ in Proposition~\ref{FSG} to obtain
\[
\llVert f_{st}-f_{s\tau_J(t)}\rrVert _{(p+1)}\le
A''2^{-J}
\]
for all $s\le t\le T$. Hence
\[
\int_0^t\bigl\langle f_{st}-f_{s\tau_J(t)},dM^N_s
\bigr\rangle\le 2^{-J}A''\mathbb{E}\int
_0^T\bigl\langle 1+\llvert v\rrvert
^{p+1},\bigl\llvert dM^N_s\bigr\rrvert \bigr
\rangle,
\]
so estimate (\ref{GEE}) shows that, as $J\to\infty$,
\[
\biggl\llVert \sup_{t\le T}\sup_{f\in\mathcal{F}(p)}\int
_0^t\bigl\langle f_{st}-f_{s\tau
_J(t)},dM^N_s
\bigr\rangle\biggr\rrVert _2\to0.
\]
Hence (\ref{FETA}) follows from (\ref{FET}), (\ref{HEE}) and (\ref{DEE}).
The proof is the same for $d=2$ except that we get $N^{-1/2}\log_2 N$
in place of $N^{-1/d}$ in (\ref{FET}) and (\ref{DEE}).
\end{pf*}

\begin{pf*}{Proof of Proposition \protect\ref{WMMF}}
Fix $p\in[2,\infty)$ and $f\in\mathcal{F}(p)$.
We follow the preceding proof to obtain, for $t\le T$,
\begin{eqnarray}
&&\int_0^t\bigl\langle\tilde
f_{sT},dM^N_s\bigr\rangle
\nonumber
\\[-8pt]
\\[-8pt]
\nonumber
&&\qquad=\sum
_{\ell=2}^L\sum_{k=0}^K
\sum_{B\in\mathcal{P}_{k,\ell}}2^{-(1+\delta)k}a_B\int
_0^t\bigl\langle \tilde h^B_{sT},dM^N_s
\bigr\rangle+\int_0^t\bigl\langle\tilde
g_{sT},dM^N_s\bigr\rangle,
\end{eqnarray}
where $\tilde f_{sT}=(E^1_{sT}-E^2_{sT})f$, $\tilde
h_{sT}^B=(E^1_{sT}-E_{sT}^2)h^B$ and $\tilde g_{sT}=(E^1_{sT}-E_{sT}^2)g$.
By Proposition~\ref{FSG}, we have
\[
\bigl\llvert \tilde g_{sT}(v)\bigr\rrvert \le\tilde A\bigl(1+\llvert
v\rrvert ^{p+2}\bigr),
\]
where
\[
\tilde A=\bigl(2^{-K}+2^{-L}\bigr)c_d''CTd^*(p+2)e^{C\bar m^*(p+3)T}
\]
and $C=C(p+1)<\infty$. Note that
\[
\bigl\llvert \tilde h_{sT}^B(v)\bigr\rrvert =\bigl
\llvert E_{(s,v)}\bigl\langle h^B,\tilde\Lambda_T^1-
\tilde \Lambda_T^2\bigr\rangle \bigr\rrvert \le
E_{(s,v)}\bigl\langle h^B,\hat\Gamma_T^1+
\hat\Gamma_T^2\bigr\rangle,
\]
so by Lemma~\ref{DCL},
\[
\sum_{k=0}^K\sum
_{B\in\mathcal{P}_{k,\ell}}\bigl(\tilde h_{sT}^B(v_0)
\bigr)^2 \le64\bigl(E_{(s,v_0)}\bigl\langle1+\llvert v\rrvert
^q,\hat\Gamma_T^1+\hat\Gamma
_T^2\bigr\rangle\bigr)^2 \le\bigl(\tilde
A'\bigr)^2\bigl(1+\llvert v_0\rrvert
^{q+1}\bigr)^2,
\]
where
\[
\tilde A'=CTd^*(q+1)e^{C\bar m^*(q+2)T}
\]
and $C=C(q+1)<\infty$.
We continue to follow the steps of the preceding proof to arrive at
\begin{eqnarray*}
&&\biggl\llVert \sup_{t\le T}\sup_{f\in\mathcal{F}(p)}\int
_0^t\bigl\langle \tilde f_{sT},dM^N_s
\bigr\rangle\biggr\rrVert _2
\\
&&\qquad \le \Biggl(\sum_{\ell=2}^L2^{(d/2-1)\ell}c(d,
\delta )\tilde A' \biggl(\frac{C}N \biggr)^{1/2}+
\tilde A\bigl(C(1+T)\bigr)^{1/2} \Biggr)\\
&&\qquad\quad{}\times \biggl(\mathbb{E}\int
_0^T\bigl\langle\llvert v\rrvert
^{2q+3},\mu_s^N\bigr\rangle \,ds
\biggr)^{1/2}.
\end{eqnarray*}
Replace $\tilde A$, $\tilde A'$, $L$ and $q$ by their values and let
$K\to\infty$ to deduce that, for some constant $C(d,\delta,p)<\infty$,
\begin{eqnarray}
&& \biggl\llVert \sup_{t\le T}\sup
_{f\in\mathcal{F}(p)}\int_0^t\bigl\langle
\tilde f_{sT},dM^N_s\bigr\rangle\biggr\rrVert
_2
\nonumber
\\
\label{FETT}
&&\qquad \le CTd^*(p+2+\delta)N^{-1/d}e^{C\bar m^*(p+3+\delta)T}\\
\nonumber
&&\qquad\quad{}\times \biggl(
\mathbb {E}\int_0^T\bigl\langle \llvert v
\rrvert ^{2p+5+2\delta},\mu_s^N\bigr\rangle \,ds
\biggr)^{1/2}.
\end{eqnarray}
Now
\begin{eqnarray}
\nonumber
&&\biggl\llVert \sup_{t\le T}\sup_{f\in\mathcal{F}(p)}
\int_0^t\bigl\langle \tilde f_{st},dM^N_s
\bigr\rangle\biggr\rrVert _2\\
\nonumber
&&\qquad\le\biggl\llVert \sup_{t\le T}
\sup_{f\in
\mathcal{F}
(p)}\int_0^t\bigl
\langle\tilde f_{sT},dM^N_s\bigr\rangle\biggr
\rrVert _2
\\[-8pt]
\label{FETU}
\\[-8pt]
\nonumber
&&\qquad\quad{}+\sum_{j=J_0+1}^J\biggl\llVert \sup
_{t\le T}\sup_{f\in\mathcal
{F}(p)}\int_0^t
\bigl\langle \tilde f_{s\tau_j(t)}-\tilde f_{s\tau_{j-1}(t)},dM^N_s
\bigr\rangle\biggr\rrVert _2
\\
\nonumber
&&\quad\qquad {}+\biggl\llVert \sup_{t\le T}
\sup_{f\in\mathcal{F}(p)}\int_0^t\bigl
\langle\tilde f_{st}-\tilde f_{s\tau_J(t)},dM^N_s
\bigr\rangle\biggr\rrVert _2,
\end{eqnarray}
and the final term tends to $0$ as $J\to\infty$.
We consider the case where $\rho_t^1=\rho_t$ and $\rho^2_t=0$ for
all $t$, from which the general case follows by the triangle inequality.
Then $f^2_{st}=f$ for all $s$ and $t$, so $\tilde f_{s\tau
_j(t)}-\tilde
f_{s\tau_{j-1}(t)}=f_{s\tau_j(t)}-f_{s\tau_{j-1}(t)}$.
We then use (\ref{FETT}) for the first term on the right in (\ref
{FETU}), use (\ref{DEE}) for the sum over $j$ and let $J\to\infty$ to
obtain the claimed estimate.
\end{pf*}

\section{Proof of Theorem \texorpdfstring{\protect\ref{MR10}}{1.2}}\label{MRD}
We seek to show that, for $p>8$ and $\varepsilon>0$, for $N\le N'$ and
any two
Kac processes $(\mu_t^N)_{t\ge0}$ and $(\mu_t^{N'})_{t\ge0}$ with
collision kernel~$B$, which are adapted to a common filtration
$(\mathcal{F}
_t)_{t\ge0}$,
with probability exceeding $1-\varepsilon$, for all $t\in[0,T]$, we have
\[
W\bigl(\mu_t^N,\mu_t^{N'}\bigr)=
\sup_{f\in\mathcal{F}}\bigl\langle f,\mu_t^N-\mu
_t^{N'}\bigr\rangle\le C\bigl(W\bigl(\mu_0^N,
\mu_0^{N'}\bigr)+N^{-1/d}\bigr)
\]
for some constant $C<\infty$ depending only on $B,d,\varepsilon
,{\lambda},p$
and $T$, where ${\lambda}$ is an upper bound for $\langle|v|^p,\mu
_0^N\rangle$
and $\langle|v|^p,\mu_0^{N'}\rangle$.
Recall the representation formula of Proposition~\ref{FMU}.
For all $f\in\mathcal{F}$, we have
\[
\bigl\langle f,\mu_t^N-\mu_t^{N'}
\bigr\rangle=\bigl\langle f_{0t},\mu_0^N-\mu
_0^{N'}\bigr\rangle+\int_0^t
\bigl\langle f_{st},dM_s^N\bigr\rangle-\int
_0^t\bigl\langle f_{st},dM^{N'}_s
\bigr\rangle,
\]
where $M^N$ is given by (\ref{KLL}) and
$f_{st}(v)=E_{st}f(v)=E_{(s,v)}\langle f,\tilde\Lambda_t\rangle$,
with $(\Lambda^*_t)_{t\ge s}$ a linearized Kac process in environment
$\rho
_t=(\mu^N_t+\mu_t^{N'})/2$.
We showed a suitable bound for $\langle f_{0t},\mu_0^N-\mu
_0^{N'}\rangle$ in
Section~\ref{PMR}.
We now show that the stochastic convolution estimates just obtained
allow us to control $\int_0^t\langle f_{st},dM_s^N\rangle$ with rate
$N^{-1/d}$,
notwithstanding the fact that the functions $f_{st}$ depend on the
random environment $(\rho_r)_{r\in[s,t]}$ and therefore are anticipating.

It will suffice to consider the case where $p\in(8,9]$ and
$T=2^{-J_0}$ for some $J_0\in\mathbb{Z}$.
Set $\delta=(p-8)/6$.
Set $\sigma_j(t)=2^{-j}\lfloor2^jt\rfloor$, and note that $\sigma
_{J_0}(t)=0$ for all $t<T$.
Set $\rho_t^j=\rho_{\sigma_j(t)}$, and define
\[
E_{st}^jf(v)=E_{(s,v)}\bigl\langle f,\tilde
\Lambda_t^j\bigr\rangle,
\]
where $(\Lambda^{j,*}_t)_{t\ge s}$ is a linearized Kac process in
environment $(\rho_t^j)_{t\ge0}$ starting from~$v$.
Then, for $t\le T$ and $J\ge J_0$, we have
\begin{eqnarray}
\int_0^t\bigl\langle
E_{st}f,dM^N_s\bigr\rangle &=&\int
_0^t\bigl\langle E_{st}^{J_0}f,dM^N_s
\bigr\rangle +\sum_{j=J_0+1}^J\int
_0^t\bigl\langle \bigl(E_{st}^j-E_{st}^{j-1}
\bigr)f,dM^N_s\bigr\rangle
\nonumber
\\[-8pt]
\label{EBIT}
\\[-8pt]
\nonumber
&&{}+\int_0^t
\bigl\langle\bigl(E_{st}-E_{st}^J
\bigr)f,dM^N_s\bigr\rangle.
\end{eqnarray}
Note that $\rho_t^{J_0}=\rho_0$ for all $t<T$.
Take $p=2$ in Proposition~\ref{WMME} to see that, for some constant
$C(d,\delta)<\infty$, we have
\begin{eqnarray}
&& \biggl\llVert \sup_{t\le T}\sup
_{f\in\mathcal{F}}\int_0^t\bigl\langle
E_{st}^{J_0}f,dM^N_s\bigr\rangle
\biggr\rrVert _2
\nonumber
\\[-9pt]
\label{EA}\\[-9pt]
\nonumber
&&\qquad\le C\kappa N^{-1/d}e^{C\langle1+|v|^{5+\delta},\rho
_0\rangle T} \biggl(
\mathbb{E} \int_0^T\bigl\langle1+|v|^{9+2\delta},
\mu_s^N\bigr\rangle \,ds \biggr)^{1/2}.
\end{eqnarray}
Fix $j$, and set $t_i=i2^{-j}$.
Note that, for $t\in[t_i,t_{i+1})$,
\[
\rho_t^j-\rho_t^{j-1}=
\cases{0,& \quad$\mbox{if $i$ is even}$,\vspace*{3pt}
\cr
\rho_{t_i}-\rho_{t_{i-1}},&\quad $\mbox{if $i$ is odd}$.} %
\]
We have
\begin{eqnarray*}
&& \int_0^t\bigl\langle\bigl(E_{st}^j-E_{st}^{j-1}
\bigr)f,dM^N_s\bigr\rangle \\[-2pt]
&&\qquad=\sum
_{i=0}^{\lfloor2^jt\rfloor-1}\int_{t_i}^{t_{i+1}}
\bigl\langle \bigl(E_{st}^j-E_{st}^{j-1}
\bigr)f,dM^N_s\bigr\rangle+\int_{\sigma_j(t)}^t
\bigl\langle \bigl(E_{st}^j-E_{st}^{j-1}
\bigr)f,dM^N_s\bigr\rangle.
\end{eqnarray*}
For $s\le t_{i+1}\le t$, we have
$E_{st}^j=E_{st_{i+1}}^jE_{t_{i+1}t}^j$ so, for all $f\in\mathcal{F}$,
\begin{eqnarray*}
&&\int_0^t\bigl\langle\bigl(E_{st}^j-E_{st}^{j-1}
\bigr)f,dM^N_s\bigr\rangle\\
&&\qquad\le \sum
_{i=0}^{\lfloor2^jt\rfloor-1}\bigl\llVert E_{t_{i+1}t}^jf-E_{t_{i+1}t}^{j-1}f
\bigr\rrVert _{(3)}\sup_{f\in\mathcal
{F}(3)}\int
_{t_i}^{t_{i+1}}\bigl\langle E_{st_{i+1}}^jf,dM^N_s
\bigr\rangle
\\
&&\qquad\quad {}+\sum_{i=0}^{\lfloor2^jt\rfloor-1}\bigl\llVert
E_{t_{i+1}t}^{j-1}f\bigr\rrVert \sup_{f\in\mathcal{F}}\int
_{t_i}^{t_{i+1}}\bigl\langle \bigl(E_{st_{i+1}}^j-E_{st_{i+1}}^{j-1}
\bigr)f,dM^N_s\bigr\rangle\\
&&\qquad\quad{}+\int_{\sigma
_j(t)}^t
\bigl\langle \bigl(E_{st}^j-E_{st}^{j-1}
\bigr)f,dM^N_s\bigr\rangle.
\end{eqnarray*}
Fix $A\ge1$, and consider the event $\Omega_0=\Omega_1\cap\Omega
_2\cap\Omega_3$, where
\begin{eqnarray*}
\Omega_1&=& \Bigl\{\sup_{t\le T}\bigl
\langle1+|v|^{5+\delta},\rho _t\bigr\rangle\le A \Bigr\},
\\
\Omega_2&=& \biggl\{\int_0^T\bigl
\langle1+|v|^3,\bigl\llvert \rho_t^j-\rho
_t^{j-1}\bigr\rrvert \bigr\rangle \,dt\le A2^{-j(1-\delta)}
\mbox{ for all }j\ge J_0+1 \biggr\},
\\
\Omega_3&=& \biggl\{\int_0^T\bigl
\langle1+|v|^3,\bigl\llvert \rho_t^J-
\rho_t\bigr\rrvert \bigr\rangle \,dt\le A2^{-J} \biggr\}.
\end{eqnarray*}
By Proposition~\ref{ME}, there is an absolute constant $C<\infty$
such that
\[
\mathbb{E} \Bigl(\sup_{t\le T}\bigl\langle1+|v|^{5+\delta},
\rho _t\bigr\rangle \Bigr)\le C{\lambda }(1+T),\qquad \mathbb{E}\bigl(\bigl
\langle1+|v|^3,\llvert \rho_t-\rho_s\rrvert
\bigr\rangle\bigr)\le C{\lambda}|t-s|
\]
so
\begin{eqnarray*}
\mathbb{E}\int_0^T\bigl\langle1+|v|^3,
\bigl\llvert \rho_t^j-\rho_t^{j-1}
\bigr\rrvert \bigr\rangle \,dt &\le &  CT{\lambda }2^{-j}, \\
\mathbb{E}\int
_0^T\bigl\langle1+|v|^3,\bigl
\llvert \rho_t^J-\rho_t\bigr\rrvert \bigr
\rangle \,dt &\le &  CT{\lambda}2^{-J}.
\end{eqnarray*}
Hence
\[
\mathbb{P}(\Omega\setminus\Omega_1)\le C{\lambda}(1+T)A^{-1},\qquad
\mathbb{P}(\Omega\setminus\Omega_3)\le C{\lambda}TA^{-1}
\]
and
\[
\mathbb{P}(\Omega\setminus\Omega_2)\le\sum
_{j=J_0+1}^\infty CT{\lambda}A^{-1}2^{-j\delta
}=CT^{1-\delta}{
\lambda}A^{-1}\bigl(2^\delta-1\bigr)^{-1}.
\]
Hence we can choose $A(\varepsilon,{\lambda},p,T)<\infty$ such that
$\mathbb{P}(\Omega
_0)\ge1-\varepsilon/2$.
By Proposition~\ref{FSG}, there is an absolute constant $C<\infty$
such that, for $f\in\mathcal{F}$ and $i\le\sigma_j(t)-1$,
\begin{eqnarray*}
&& \bigl\llVert E_{t_{i+1}t}^jf-E_{t_{i+1}t}^{j-1}f
\bigr\rrVert _{(3)}\\
&&\qquad \le C\kappa\exp \biggl\{C\int_{t_{i+1}}^t
\bigl\langle1+|v|^4,\rho _s^j+\rho
_s^{j-1}\bigr\rangle \,ds \biggr\}\int_{t_{i+1}}^t
\bigl\langle1+|v|^3,\bigl\llvert \rho _s^j-
\rho _s^{j-1}\bigr\rrvert \bigr\rangle \,ds.
\end{eqnarray*}
Also, by Proposition~\ref{JPE},
\[
\bigl\llVert E_{t_{i+1}t}^{j-1}f\bigr\rrVert \le3\bigl(1+6
\kappa(t-t_{i+1})\bigr)\exp \biggl\{8\int_{t_{i+1}}^t
\bigl\langle1+|v|^3,\rho_s^{j-1}\bigr\rangle \,ds
\biggr\}.
\]
So, on $\Omega_0$, for some absolute constant $C<\infty$, we have
\begin{eqnarray*}
&&\sup_{t\le T}\sup_{f\in\mathcal{F}}\int
_0^t\bigl\langle \bigl(E_{st}^j-E_{st}^{j-1}
\bigr)f,dM^N_s\bigr\rangle
\\
&&\qquad\le CA\kappa e^{CAT}\sum_{i=0}^{2^jT-1}
\biggl(2^{-j(1-\delta)}\sup_{f\in
\mathcal{F}(3)}\int_{t_i}^{t_{i+1}}
\bigl\langle E_{st_{i+1}}^jf,dM^N_s
\bigr\rangle \\
&&\qquad\qquad\hspace*{68pt}{}+\sup_{t\in[t_i,t_{i+1}]}\sup_{f\in\mathcal{F}}\int
_{t_i}^t\bigl\langle \bigl(E_{st}^j-E_{st}^{j-1}
\bigr)f,dM^N_s\bigr\rangle \biggr).
\end{eqnarray*}
Set $\mathcal{F}_t=\sigma\{\mu_s^N,\mu_s^{N'} \dvtx  s\in[0,t]\}$.
We apply Propositions \ref{WMME} and \ref{WMMF} conditionally on
$\mathcal{F}
_{t_i}$ to obtain, for some constant $C(d,\delta)<\infty$,
\begin{eqnarray*}
&& \biggl\llVert \sup_{f\in\mathcal{F}(3)}\int_{t_i}^t
\bigl\langle E^j_{st}f,dM^N_s
\bigr\rangle1_{\{\langle
1+|v|^{5+\delta},\rho_{t_i}\rangle\le A\}}\biggr\rrVert _2 \\
&&\qquad\le C\kappa
e^{CA2^{-j}}N^{-1/d} \biggl(\mathbb{E}\int_{t_i}^{t_{i+1}}
\bigl\langle|v|^{9+2\delta
},\mu_s^N\bigr\rangle \,ds
\biggr)^{1/2}
\end{eqnarray*}
and
\begin{eqnarray*}
&&\biggl\llVert \sup_{t\in[t_i,t_{i+1}]}\sup_{f\in\mathcal{F}}\int
_{t_i}^t\bigl\langle \bigl(E_{st}^j-E_{st}^{j-1}
\bigr)f,dM^N_s\bigr\rangle1_{\{\langle1+|v|^{5+\delta
},\rho_{t_i}\rangle\le
A\}}\biggr\rrVert
_2
\\
&&\qquad \le CA2^{-j}\kappa e^{CA2^{-j}}N^{-1/d} \biggl(
\mathbb{E}\int_{t_i}^{t_{i+1}}\bigl\langle
|v|^{9+2\delta},\mu _s^N\bigr\rangle \,ds
\biggr)^{1/2}.
\end{eqnarray*}
By Proposition~\ref{ME}, there is a constant $C(B,p)<\infty$ such that
\[
\mathbb{E}\int_{t_i}^{t_{i+1}}\bigl
\langle|v|^{9+2\delta},\mu _s^N\bigr\rangle \,ds\le\int
_{t_i}^{t_{i+1}}C{\lambda}\bigl(1+t^{p-9-2\delta}
\bigr)\,dt=C{\lambda }\bigl(2^{-j}+2^{-4\delta j}/(4\delta)\bigr).
\]
Hence, for constants $C(B,d,p)<\infty$,
\begin{eqnarray}
\nonumber
&&\biggl\llVert \sup_{t\le T}\sup_{f\in\mathcal{F}}
\int_0^t\bigl\langle \bigl(E_{st}^j-E_{st}^{j-1}
\bigr)f,dM^N_s\bigr\rangle1_{\Omega_0}\biggr\rrVert
_2
\\
\label{EB}
&&\qquad\le CA^2{\lambda}^{1/2}\kappa
^2e^{CAT}2^jT\bigl(2^{-j(1-\delta)}+2^{-j}
\bigr) \bigl(2^{-j/2}+\delta ^{-1/2}2^{-2j\delta
}
\bigr)N^{-1/d}
\\
\nonumber
&&\qquad \le A{\lambda}^{1/2}\kappa ^2e^{CAT}
\bigl(2^{-j\delta
}+2^{-j/2}\bigr)N^{-1/d}.
\end{eqnarray}
Here, we absorbed $4\delta^{-1/2}CAT$ into $e^{CAT}$ in the second
inequality by changing the constant $C$.
By Proposition~\ref{FSG}, there is an absolute constant $C<\infty$
such that, for all $f\in\mathcal{F}$ and $s\le t\le T$,
\begin{eqnarray*}
&&\bigl\llvert E_{st}^Jf(v)-E_{st}f(v)\bigr
\rrvert\\
&&\qquad \le C\kappa\bigl(1+|v|^3\bigr)\exp \biggl\{C\int
_0^T\bigl\langle 1+|v|^4,
\rho_t^J+\rho_t\bigr\rangle \,dt \biggr\}\int
_s^t\bigl\langle 1+|v|^3,\bigl
\llvert \rho _r^J-\rho_r\bigr\rrvert
\bigr\rangle \,dr.
\end{eqnarray*}
So, on $\Omega_0$, we have
\[
\bigl\llvert E_{st}^Jf(v)-E_{st}f(v)\bigr
\rrvert \le CA\kappa e^{CAT}2^{-J}\bigl(1+|v|^3
\bigr),
\]
and so as $J\to\infty$,
\begin{eqnarray}
&&\biggl\llVert \sup_{t\le T}\sup
_{f\in\mathcal{F}}\int_0^t\bigl\langle
\bigl(E_{st}^J-E_{st}\bigr)f,dM^N_s
\bigr\rangle1_{\Omega_0}\biggr\rrVert _2
\nonumber
\\[-8pt]
\label{EC}
\\[-8pt]
\nonumber
&&\qquad\le CA\kappa
e^{CAT}2^{-J}\biggl\llVert \int_0^T
\bigl\langle1+|v|^3,\bigl\llvert dM_s^N\bigr
\rrvert \bigr\rangle\biggr\rrVert _2\to0.
\end{eqnarray}
Finally, we use estimates (\ref{EA}), (\ref{EB}) and (\ref{EC}) in
(\ref{EBIT}) and let $J\to\infty$ to obtain a constant
$C(B,d,\varepsilon
,{\lambda},p,T)<\infty$ such that
\[
\biggl\llVert \sup_{t\le T}\sup_{f\in\mathcal{F}}\int
_0^t\bigl\langle E_{st}f,dM^N_s
\bigr\rangle1_{\Omega
_0}\biggr\rrVert _2\le CN^{-1/d}.
\]
An analogous estimate holds for $N'$, and Theorem~\ref{MR10} then
follows by Chebyshev's inequality.

\section{Properties of the distance function}\label{W}
Recall that $W \dvtx  \mathcal{S}\times\mathcal{S}\to[0,4]$ is defined by
\[
W(\mu,\nu)=\sup_{f\in\mathcal{F}}\langle f,\mu-\nu\rangle,
\]
where $\mathcal{F}$ is the set of functions $f$ on $\mathbb{R}^d$
such that $|\hat
f(v)|\le1$ and $|\hat f(v)-\hat f(v')|\le|v-v'|$ for all $v,v'$,
where $\hat f(v)=f(v)/(1+|v|^2)$.
\begin{proposition}\label{CMS}
The metric space $(\mathcal{S},W)$ is complete and separable.
\end{proposition}
\begin{pf}
Write $\mathcal{P}$ for the set of Borel probability measures on
$\mathbb{R}^d$, and
define $\Phi \dvtx  \mathcal{S}\to\mathcal{P}$ by $\Phi(\mu)(dv)=\frac{1}2(1+|v|^2)\mu(dv)$.
Write $W_1^\rho$ for the Wasserstein-$1$ metric on $\mathcal{P}$ associated
with the bounded metric $\rho(v,v')=|v-v'|\wedge2$ on $\mathbb{R}^d$.
Then $W(\mu,\nu)=2W_1^\rho(\Phi(\mu),\Phi(\nu))$ for all $\mu
,\nu\in\mathcal{S}$.
Now $(\mathcal{P},W_1^\rho)$ is complete and separable, and $\Phi
(\mathcal{S})$ is
closed in $\mathcal{P}$ under $W_1^\rho$, so $(\mathcal{S},W)$ is
also complete and
separable.
\end{pf}

We will prove two approximation schemes for a measure $\mu$ in the
Boltzmann sphere $\mathcal{S}$, by empirical distributions of systems
of $N$ particles.
The first uses the empirical distribution $\mu^N=\frac{1}N\sum_{i=1}^N\delta_{V_i}$ of a sample of $N$ independent random variables
$V_1,\dots,V_N$ with distribution $\mu$.
The convergence of $\mu^N$ to $\mu$ has been extensively investigated
for standard Wasserstein distances; see \cite{MR3127919} or \cite{1312.2128}.
We modify some of the simpler ideas from \cite{1312.2128} to obtain
estimates for the weighted Wasserstein distance $W$ used in this paper.
The sample empirical distribution $\mu^N$ is not, however, a random
variable in the Boltzmann sphere.
Set
\[
\bar V_N=\frac{1}N\sum_{i=1}^NV_i,\qquad
S_N=\frac{1}N\sum_{i=1}^N
\llvert V_i-\bar V\rrvert ^2,\qquad \tilde
V_i=S^{-1/2}_N(V_i-\bar
V_N).
\]
On the event $\{S_N>0\}$, define the \textit{rescaled empirical
distribution} $\tilde\mu^N= \frac{1}N\sum_{i=1}^N\delta_{\tilde V_i}$.
On the event $\{S_N=0\}$, we take $\tilde\mu^N$ to be some arbitrary
element of $\mathcal{S}_N$.
Then $\tilde\mu^N$ is a random variable in $\mathcal{S}_N$.
We will quantify the convergence of $\tilde\mu^N$ to $\mu$ in
weighted Wasserstein distance, using the convergence of $\mu^N$ as an
intermediate step.

\begin{proposition}\label{DAT}
For all $d\ge3$ and all $\mu\in\mathcal{S}$, we have $\mathbb
{E}(W(\mu^N,\mu
))\to0$ as $N\to\infty$.
Set
\[
\beta=\beta(p)= %
\cases{ (p-2)/(p+d),&\quad$ \mbox{if $p\in\bigl(2,3d/(d-1)
\bigr)$,}$\vspace*{3pt}
\cr
1/d,& \quad$\mbox{if $p\in[3d/(d-1),\infty)$}$.}
\]
For all $p\in(2,\infty)\setminus\{3d/(d-1)\}$, there is a constant
$C(d,p)<\infty$ such that, for all $N\in\mathbb{N}$,
\begin{equation}
\label{F1} \mathbb{E}\bigl(W\bigl(\mu^N,\mu\bigr)\bigr)\le C\bigl
\langle|v|^p,\mu\bigr\rangle N^{-\beta}.
\end{equation}
For $d=2$ and $p\in(2,\infty)\setminus\{3d/(d-1)\}$, or for $d\ge3$ and
$p=3d/(d-1)$, the same estimate holds with an additional factor of
$\log(N+1)$ on the right-hand side.
In the case when both $d=2$ and $p=6$, the additional factor is squared.
\end{proposition}

\begin{proposition}\label{DAP}
The conclusions of Proposition~\ref{DAT} remain valid if $\mu^N$ is
replaced by $\tilde\mu^N$ and $\beta$ is replaced by $\tilde\beta
=\beta
\wedge((p-2)^2/(3p-4))$.
Moreover, the constant $C$ may be chosen so that, for all $N\in\mathbb{N}$,
there is an event $\Omega(1/4)$, of probability exceeding $1-C\langle
|v|^p,\mu
\rangle N^{-(p/4)\wedge(p/2-1)}$, such that
\[
\label{DA2} \mathbb{E}\bigl(\bigl\langle|v|^p,\tilde
\mu^N\bigr\rangle1_{\Omega(1/4)}\bigr)\le C\bigl
\langle|v|^p,\mu\bigr\rangle.
\]
\end{proposition}

It is simple to check that $\tilde\beta=\beta$ whenever $d\ge2$ and
$p\ge3$.
The following example shows that the exponent $\beta(p)$ cannot be
improved for $p\in(2,3d/(d-1))$ and hence that the moment threshold
$p=3d/(d-1)$ for convergence with optimal rate $N^{-1/d}$ also cannot
be improved.
Fix $p>2$ and $q>d+p$, and consider the measure $\mu(dv)=c1_{\{|v|>r\}
}|v|^{-q}\,dv$ where $c$ and $r$ are determined so that $\mu\in\mathcal{S}$.
Then $\langle|v|^p,\mu\rangle<\infty$.
Define
\[
f_N(v)=\bigl(\operatorname{dist}\bigl(v,\operatorname{supp}
\mu^N\bigr)\wedge1\bigr) \bigl(1+|v|^2\bigr).
\]
Then $f_N\in\mathcal{F}$, so
\[
\mathbb{E}\bigl(W\bigl(\mu,\mu^N\bigr)\bigr)\ge\mathbb{E}\bigl(\bigl
\langle f_N,\mu-\mu ^N\bigr\rangle\bigr)=\mathbb{E}\bigl(
\langle f_N,\mu\rangle\bigr).
\]
There are constants $a<\infty$ and $r_0\ge r$ such that $\mu(\{u\in
\mathbb{R}^d \dvtx  |u-v|\ge1\})\ge e^{-a|v|^{-q}}$ whenever $|v|\ge r_0$.
Then $\operatorname{dist}(v,\operatorname{supp}\mu^N)\ge1$ with
probability at least $e^{-Na|v|^{-q}}$.
Hence
\begin{eqnarray*}
\mathbb{E}\bigl(\langle f_N,\mu\rangle\bigr) &\ge &  c
\sigma_{d-1}\int_{r_0}^\infty
e^{-Nat^{-q}}t^{d+1-q}\,dt\\
&=& c\sigma_{d-1}N^{-1+(d+2)/q}\int
_{r_0N^{-1/q}}^\infty e^{-as^{-q}}s^{d+1-q}\,ds.
\end{eqnarray*}
Consider the limit $q\to p+d$ in the case $p<3d/(d-1)$. Then
$1-(d+2)/q\to(p-2)/(p+d)<1/d$, so we have justified the claims made above.

\begin{pf*}{Proof of Proposition~\ref{DAT}}
The following estimate is known for the $N$-sample empirical
distribution $\mu^N_0$ of a probability measure $\mu_0$ supported on
$B_0=(-1,1]^d$.
For all $d\ge3$, there is a constant $C(d)<\infty$ such that, for all
$N\in\mathbb{N}$, we have
\begin{equation}
\label{SWE} \mathbb{E}\bigl(W_1\bigl(\mu^N_0,
\mu_0\bigr)\bigr)\le C(d)N^{-1/d}.
\end{equation}
Here, $W_1$ denotes the Wasserstein-$1$ distance for the Euclidean
metric on $\mathbb{R}^d$.
For completeness, and since it may be read as a warm-up for the proof
of Proposition~\ref{WMME}, we give a proof.
Fix $L\in\mathbb{N}$.
For $\ell=0,1,\dots,L-1$, we can partition $B_0$ as a set $\mathcal
{P}_\ell
$ of $2^{\ell d}$ translates of $(-2^{-\ell},2^{-\ell}]^d$.
Fix a function $f$ on $\mathbb{R}^d$ with $f(0)=0$ and
$|f(v)-f(v')|\le
|v-v'|$ for all $v,v'\in B_0$.
Then we can write
\[
f=\sum_{\ell=0}^{L-1}\sum
_{B\in\mathcal{P}_\ell}a_B1_B+g,
\]
where $a_{B_0}=\langle f\rangle_{B_0}$ and $a_B=\langle f\rangle
_B-\langle f\rangle_{\pi(B)}$ for $B\in
\mathcal{P}_\ell$ and $\ell\ge1$.
Here we have written $\langle f\rangle_B$ for the average of $f$ over
$B$ and $\pi
(B)$ for the unique element of $\mathcal{P}_{\ell-1}$ containing $B$.
By\vspace*{1pt} (\ref{LIPA}) and (\ref{LIPB}), we have $|a_B|\le2^{-\ell+1}c_d$
for all $B\in\mathcal{P}_\ell$ and all $\ell$, and $|g(v)|\le2^{-L+2}c_d$.
So, by Cauchy--Schwarz,
\begin{eqnarray*}
\bigl\langle f,\mu_0^N-\mu_0\bigr\rangle &=&
\sum_{\ell=0}^{L-1}\sum
_{B\in\mathcal{P}_\ell}a_B\bigl(\mu _0^N(B)-
\mu _0(B)\bigr)+\bigl\langle g,\mu_0^N-
\mu_0\bigr\rangle
\\
&\le & 2c_d\sum_{\ell=0}^{L-1}2^{(d-2)\ell/2}
\biggl(\sum_{B\in
\mathcal{P}
_\ell}\bigl(\mu_0^N(B)-
\mu_0(B)\bigr)^2 \biggr)^{1/2}+8c_d2^{-L}.
\end{eqnarray*}
The right-hand side does not depend on $f$, so it is an upper bound for
$W_1(\mu_0^N,\mu_0)$, by duality.
Note that $\operatorname{var}(\mu_0^N(B))\le\mu_0(B)/N$.
Now take expectations and use Cauchy--Schwarz again to obtain
\[
\mathbb{E}\bigl(W_1\bigl(\mu_0^N,
\mu_0\bigr)\bigr)\le2c_d\sum_{\ell
=0}^{L-1}2^{(d-2)\ell
/2}N^{-1/2}+8c_d2^{-L}.
\]
We optimize at $L=\lceil\log_2(N+1)/d\rceil$ to obtain (\ref{SWE}).
The same argument produces $N^{-1/2}\log_2(N+1)$ on the right when $d=2$.

Set $B_k=2^kB_0$.
Fix $K\in\mathbb{N}\cup\{\infty\}$, and partition $\mathbb{R}^d$
as $\bigcup_{k=0}^KA_k$, where $A_0=B_0$, $A_k=B_k\setminus B_{k-1}$ for $1\le
k<K$, and
$A_K=\mathbb{R}^d\setminus(\bigcup_{k=0}^{K-1}A_k)$.
Set $p_k=\mu(A_k)$ and write $\mu_k$ for the conditional distribution
of $\mu$ on $A_k$.
Write $N_k$ for the number of elements of the sample falling in $A_k$
and write $\mu^{N_k}_k$ for the empirical distribution of this sub-sample.
Set $\hat p_k=N_k/N$.
Then
\[
\mu=\sum_{k=0}^K p_k
\mu_k, \qquad \mu^N=\sum_{k=0}^K
\hat p_k\mu^{N_k}_k.
\]
Fix a function $f$ on $\mathbb{R}^d$ such that $|\hat f(v)|\le1$ and
$|\hat
f(v)-\hat f(v')|\le|v-v'|$ for all $v,v'$, where $\hat f(v)=f(v)/(1+|v|^2)$.
Then, for all $k$ and all $v,v'\in B_k$, we have
\begin{eqnarray*}
\bigl\llvert f(v)\bigr\rrvert  &\le & 1+d2^{2k},\\
 \bigl\llvert f(v)-f
\bigl(v'\bigr)\bigr\rrvert  &\le & \bigl(2+|v|^2+\bigl\llvert
v'\bigr\rrvert ^2\bigr)\bigl\llvert v-v'
\bigr\rrvert \le2\bigl(1+d2^{2k}\bigr)\bigl\llvert v-v'
\bigr\rrvert.
\end{eqnarray*}
Hence
\begin{eqnarray*}
\bigl\langle f,\mu^N-\mu\bigr\rangle &=&\sum
_{k=0}^{K-1}\hat p_k\bigl\langle f,
\mu_k^{N_k}-\mu_k\bigr\rangle+(\hat
p_k-p_k)\langle f,\mu_k\rangle+\bigl
\langle f,\hat p_K\mu_K^{N_K}-p_K
\mu_K\bigr\rangle
\\
&\le&\sum_{k=0}^{K-1}\bigl(1+d2^{2k}
\bigr)\bigl\{2\hat p_kW_1\bigl(\mu_k^{N_k},
\mu _k\bigr)+|\hat p_k-p_k|\bigr\}\\
&&\quad\hspace*{7pt}{}+\bigl
\langle\bigl(1+|v|^2\bigr)1_{A_K},\hat p_K\mu
_K^{N_K}+p_K\mu_K\bigr\rangle.
\end{eqnarray*}
Note the inequalities
\[
\mathbb{E}\llvert \hat p_k-p_k\rrvert
\le(2p_k)\wedge(p_k/N)^{1/2}\le
2N^{-1/d}p_k^{1-1/d}, \mathbb{E}\bigl(\hat
p_k^{1-1/d}\bigr)\le p_k^{1-1/d}.
\]
Estimate (\ref{SWE}) scales from $B_0$ to $B_k$ to give, on the event
$\{N_k\ge1\}$,
\[
\mathbb{E}\bigl(W_1\bigl(\mu_k^{N_k},
\mu_k\bigr)|N_k\bigr)\le2^kC(d)N_k^{-1/d}.
\]
Hence, on taking the supremum over $f$ and then the expectation, we obtain
\[
\mathbb{E}\bigl(W\bigl(\mu^N,\mu\bigr)\bigr) \le\sum
_{k=0}^{K-1}2^{k+2}\bigl(1+d2^{2k}
\bigr)C(d)N^{-1/d}p_k^{1-1/d}+2\bigl\langle
\bigl(1+|v|^2\bigr)1_{A_K},\mu\bigr\rangle.
\]
Since $\mu\in\mathcal{S}$, the final term on the right is small for large
$K$, so $\mathbb{E}(W(\mu^N,\mu))\to0$ as $N\to\infty$.
If $\langle|v|^p,\mu\rangle<\infty$ for some $p>2$, we can control the
right-hand side using the bounds
\[
\sum_{k=1}^{K-1}2^{p(k-1)}p_k
\le\bigl\langle|v|^p,\mu\bigr\rangle, \bigl\langle
\bigl(1+|v|^2\bigr)1_{A_K},\mu\bigr\rangle
\le2^{-(K-1)(p-2)+1}\bigl\langle |v|^p,\mu\bigr\rangle.
\]
Finally, we optimize at $K=\lceil\log_2(N+1)/(d+p)\rceil$ when $p<3d/(d-1)$
and $K=\infty$ when $p>3d/(d-1)$ to obtain the claimed estimate.
\end{pf*}

\begin{pf*}{Proof of Proposition~\ref{DAP}}
Set $Q_N=N^{-1}\sum_{i=1}^N|V_i|^2$.
Fix $\delta\in(0,1/4]$, and consider the event
\[
\Omega(\delta)=\bigl\{\llvert Q_N-1\rrvert \le\delta\mbox{ and }
\llvert \bar V_N\rrvert \le\delta \bigr\}.
\]
Note that $Q_N=S_N+|\bar V_N|^2$.
On $\Omega(\delta)$, by some simple estimation, we have
$|S_N^{-1/2}-1|\le
4\delta$, so $|\tilde V_i-V_i|\le(4|V_i|+2)\delta$.
Hence, in particular, there is a constant $C(p)<\infty$ such that
\[
\label{F5} \mathbb{E}\bigl(\bigl\langle|v|^p,\tilde
\mu^N\bigr\rangle1_{\Omega(1/4)}\bigr)\le C(p)\mathbb{E}\bigl(\bigl
\langle|v|^p,\mu^N\bigr\rangle \bigr)=C(p)\bigl
\langle|v|^p,\mu\bigr\rangle.
\]
Now, for all $f\in\mathcal{F}$, we have
\[
f(\tilde V_i)-f(V_i)\le\bigl(\llvert \tilde
V_i-V_i\rrvert \wedge1\bigr) \bigl(2+\llvert \tilde
V_i\rrvert ^2+\llvert V_i\rrvert
^2\bigr)\le24\bigl(\bigl(\delta+\delta\llvert V_i
\rrvert \bigr)\wedge1\bigr) \bigl(1+\llvert V_i\rrvert
^2\bigr)
\]
and so
\[
\bigl\langle f,\tilde\mu^N-\mu^N\bigr\rangle =
\frac{1}N\sum_{i=1}^N\bigl(f(
\tilde V_i)-f(V_i)\bigr) \le\frac{24}N\sum
_{i=1}^N\bigl(\bigl(\delta+\delta\llvert
V_i\rrvert \bigr)\wedge1\bigr) \bigl(1+\llvert V_i
\rrvert ^2\bigr).
\]
Hence
\begin{equation}
\label{F2} \mathbb{E}\bigl(W\bigl(\tilde\mu^N,\mu^N
\bigr)1_{\Omega(\delta)}\bigr)\le24\bigl\langle\bigl\{ (\delta+\delta|v|)\wedge1
\bigr\} \bigl(1+\llvert v\rrvert ^2\bigr)\mu\bigr\rangle\to0
\end{equation}
as $\delta\to0$.

Since $\langle v,\mu\rangle=0$ and $\langle|v|^2,\mu\rangle=1$, we
have $\mathbb{P}(\Omega\setminus\Omega
(\delta))\to0$ as $N\to\infty$ for all $\delta>0$ by the weak law
of large numbers.
For $p\ge2$, there is a constant $C<\infty$, depending only on $d$
and $p$, such that $\mathbb{E}(|\bar V_N|^{p/2})^2\le\mathbb
{E}(|\bar V_N|^p)\le C\langle
|v|^p,\mu\rangle N^{-p/2}$.
Hence
\[
\mathbb{P}\bigl(\llvert \bar V_N\rrvert >\delta\bigr)\le C\bigl
\langle\llvert v\rrvert ^p,\mu\bigr\rangle\delta ^{-p/2}N^{-p/4}.
\]
For $p\ge4$, since $\langle|v|^2,\mu\rangle=1$, $C$ may be chosen
so that also
$\mathbb{E}(|Q_N-1|^{p/2})\le C\langle|v|^p,\mu\rangle N^{-p/4}$ and so
\[
\mathbb{P}\bigl(\llvert Q_N-1\rrvert >\delta\bigr)\le C\bigl\langle
\llvert v\rrvert ^p,\mu\bigr\rangle\delta ^{-p/2}N^{-p/4}.
\]
For $p\in(2,4]$ we use a different estimate.
Set $R=\sqrt{\delta N}$ and write $X_i=|V_i|^2\wedge R$ and $\bar
X=N^{-1}\sum_{i=1}^NX_i$ and $x=\mathbb{E}(X_1)$.
Then $\mathbb{E}(X_1^2)\le\langle|v|^p,\mu\rangle R^{4-p}$ and
\[
\llvert x-1\rrvert \le\mathbb{E}\llvert \bar X-Q_N\rrvert \le\bigl
\langle\llvert v\rrvert ^21_{\{\llvert v\rrvert \ge R\}},\mu \bigr\rangle\le\bigl
\langle \llvert v\rrvert ^p,\mu\bigr\rangle R^{2-p}
\]
so
\begin{eqnarray*}
\mathbb{P}\bigl(\llvert Q_N-1\rrvert >\delta\bigr) &\le & \mathbb{P}
\bigl(\llvert Q_N-\bar X\rrvert >\delta/3\bigr) +\mathbb{P}\bigl(
\llvert \bar X-x\rrvert >\delta/3\bigr) +\mathbb{P}\bigl(\llvert x-1\rrvert >
\delta/3\bigr)
\\
&\le & 12\bigl\langle\llvert v\rrvert ^p,\mu\bigr\rangle
\delta^{-p/2}N^{-(p/2-1)}.
\end{eqnarray*}
Here, for the second inequality, we estimated the first term using
Markov's inequality, the second using Chebyshev
and noted that the third term vanishes except in cases where the final
estimate exceeds $1$.
We combine these estimates to see that there is a constant $C<\infty$,
depending only on $d$ and $p$, such that, for all $\delta\in(0,1/4]$,
\begin{equation}
\label{F3} \mathbb{P}\bigl(\Omega\setminus\Omega(\delta)\bigr)\le C\delta
^{-p/2}\bigl\langle\llvert v\rrvert ^p,\mu\bigr\rangle
N^{-(p/4)\wedge(p/2-1)}.
\end{equation}

Now, from (\ref{F1}), (\ref{F2}) and (\ref{F3}), for all $p\in
(2,\infty)\setminus\{3d/(d-1)\}$, all $\delta\in(0,1/4]$ and all
$N\in\mathbb{N}$,
we have
\begin{eqnarray*}
\mathbb{E}\bigl(W\bigl(\tilde\mu^N,\mu\bigr)\bigr) &\le & \mathbb{E}
\bigl(W\bigl(\mu^N,\mu\bigr)\bigr)+\mathbb{E}\bigl(W\bigl(\tilde
\mu^N,\mu ^N\bigr)1_{\Omega(\delta)}\bigr)+4\mathbb{P}
\bigl(\Omega\setminus\Omega(\delta)\bigr)
\\
&\le&  C\bigl(N^{-\beta}+\delta^{(p-2)\wedge1}+\delta ^{-p/2}N^{-(p/4)\wedge
(p/2-1)}
\bigr)\bigl\langle\llvert v\rrvert ^p,\mu\bigr\rangle.
\end{eqnarray*}
Hence $\mathbb{E}(W(\tilde\mu^N,\mu))\to0$ as $N\to\infty$.
Moreover, on optimizing over $\delta$, the terms $\delta^{(p-2)\wedge
1}$ and
$\delta^{-p/2}N^{-(p/4)\wedge(p/2-1)}$ can be absorbed in the term
$N^{-\beta}$,
except possibly when $p\in(2,3)$,
and in that case we can take $\delta=(1/4)\times\break N^{-(p-2)/(3p-4)}$ for the
desired estimate.
\end{pf*}

\section{Spatially homogeneous Boltzmann equation}\label{SHBE}
Given an initial state $\mu_0$ in the Boltzmann sphere $\mathcal{S}$,
one can
ask\vspace*{1pt} whether there exists a process $(\mu_t)_{t\ge0}$ in $\mathcal
{S}$ such that,
for all bounded measurable functions $f$ of compact support in $\mathbb{R}^d$
and all $t\ge0$,
\begin{equation}
\label{WB} \langle f,\mu_t\rangle=\langle f,\mu_0
\rangle+\int_0^t\bigl\langle f,Q(
\mu_s,\mu_s)\bigr\rangle \,ds.
\end{equation}
Here $Q$ is the Boltzmann operator, defined in equation (\ref{TBO}).
Such a process would then be called
a measure solution of the spatially homogeneous Boltzmann equation.
While the existence and uniqueness (in law) of the Kac process is
elementary, the existence and uniqueness of measure solutions is
a hard question, but one which, extending a long line of prior works,
including \cite{MR1716814,MR1697562},
has been positively answered by Lu and Mouhot \cite{MR2871802}, Theorem~1.5.

After Kac \cite{MR0084985},
important contributions to understanding the behavior of versions of
the Kac process were made by McKean \cite{MR0214112} and Tanaka \cite
{MR512334,MR799949}.
Sznitman \cite{MR753814} gave the first proof for hard spheres that
the Kac process converges weakly to solutions of the Boltzmann equation.
Mischler and Mouhot \cite{MR3069113}, Theorem~6.2, proved a
quantitative refinement of Sznitman's result, using a Wasserstein
distance on the laws
of $k$-samples from the empirical distribution.
In recent work, Fournier and Mischler \cite{1302.5810} and Cortez and
Fontbona \cite{1406.2115} have proved Wasserstein estimates for some
other particle systems associated to the spatially homogeneous
Boltzmann equation.

Our consistency estimate allows a further strengthening of Sznitman's result.
In the convergence theorem below, we obtain a pathwise estimate,
expressed in terms of a Wasserstein distance on the empirical
distribution itself,
and we are able to show, under suitable moment conditions, that the
rate of convergence is the optimal one for discrete approximations in
Wasserstein distance.
The convergence results of both Sznitman and Mischler--Mouhot are
expressed in terms of propagation of chaos,
while our estimate may be applied to any initial $N$-particle system.
For $p\ge2$, define
\[
\mathcal{S}(p)=\bigl\{\mu\in\mathcal{S} \dvtx  \bigl\langle\llvert v\rrvert
^p,\mu\bigr\rangle <\infty\bigr\},
\]
and call a solution locally bounded in $\mathcal{S}(p)$ if $\langle
|v|^p,\mu_t\rangle$
is bounded on compact time intervals.
We know from \cite{MR2871802}, Theorem~1.5,  that, for all $\mu_0\in
\mathcal{S}$, there is a unique solution $(\mu_t)_{t\ge0}$ in
$\mathcal{S}$ to
(\ref{WB}).
Sznitman's theorem assumes $\mu_0\in\mathcal{S}(3)$.
The convergence result of Mischler and Mouhot, which has good long-time
properties, assumes compactly supported initial data or at least an
exponential moment.

\begin{theorem}\label{CBE}
Assume that the collision kernel $B$ satisfies conditions (\ref
{LBP}) and (\ref{LCP}).
Let $\mu_0\in\mathcal{S}(p)$ for some $p\in(2,\infty)$.
Then there exists a unique locally bounded solution $(\mu_t)_{t\ge0}$
to (\ref{WB}) in $\mathcal{S}(p)$ starting from $\mu_0$.
Let $\varepsilon\in(0,1]$, ${\lambda}\ge\langle|v|^p,\mu_0\rangle
$ and $T\in
[0,\infty)$.
Then there exists a constant $C(B,d,\varepsilon,{\lambda},p,T)<\infty
$ with
the following property.
For all $N\in\mathbb{N}$ and any Kac process $(\mu^N_t)_{t\ge0}$ in
$\mathcal{S}
_N$ with $\langle|v|^p,\mu_0^N\rangle\le{\lambda}$,
with probability exceeding $1-\varepsilon$, for all $t\in[0,T]$, we have
\[
W\bigl(\mu^N_t,\mu_t\bigr)\le C\bigl(W\bigl(
\mu_0^N,\mu_0\bigr)+N^{-\alpha}\bigr),
\]
where $\alpha(d,p)$ is given in Theorem~\ref{MR}.
For $p>8$ and $d\ge3$, we can take $\alpha=1/d$.
For $p>8$ and $d=2$ the estimate holds with $N^{-\alpha}$ replaced by
$N^{-1/2}\log N$.
\end{theorem}
\begin{pf}
We will prove the first assertion on existence and uniqueness for
completeness, while noting, as discussed above, that a stronger
statement is already known.
Let $(V_i \dvtx  i\in\mathbb{N})$ be a sequence of independent random
variables in
$\mathbb{R}^d$ of distribution $\mu_0$.
Write $\bar V_N$ for the sample mean and $S_N$ for the sample variance
of $V_1,\dots,V_N$.
For each $N\in\mathbb{N}$, set
\[
\nu_0^N=\frac{1}N\sum
_{i=1}^N\delta_{S_N^{-1/2}(V_i-\bar V_N)}
\]
on the event $\{S_N>0\}$, and set $\nu_0^N$ equal to some arbitrary
element of $\mathcal{S}_N$ otherwise.
Conditioning on $\nu^N_0$, let $(\nu^N_t)_{t\ge0}$ be a Kac process
in $\mathcal{S}_N$ starting from $\nu^N_0$.
Choose sequences $(\varepsilon_k \dvtx  k\in\mathbb{N})$ in $(0,1]$ and
$(T_k \dvtx  k\in\mathbb{N})$ in
$[0,\infty)$ such that $\sum_k\varepsilon_k<\infty$ and $T_k\to
\infty$.
By Proposition~\ref{DAP} and Theorem~\ref{MR}, there exists an
increasing sequence $(N_k \dvtx  k\in\mathbb{N})$ in $\mathbb{N}$ such
that, for all $k\in
\mathbb{N}$, with
probability exceeding $1-\varepsilon_k$,
\[
\bigl\langle\llvert v\rrvert ^p,\nu_0^{N_k}
\bigr\rangle\le C\bigl\langle\llvert v\rrvert ^p,\mu_0
\bigr\rangle,\qquad W\bigl(\nu_0^{N_k},\mu _0\bigr)
\le C\bigl\langle\llvert v\rrvert ^p,\mu_0\bigr\rangle
\varepsilon_k
\]
and then for all $t\le T_k$
\begin{equation}
\label{UCP} W\bigl(\nu_t^{N_k},\nu_t^{N_{k+1}}
\bigr)\le C\bigl(W\bigl(\nu_0^{N_k},\nu _0^{N_{k+1}}
\bigr)+\varepsilon_k\bigr).
\end{equation}
By Borel--Cantelli, almost surely, these inequalities hold for all
sufficiently large $k$,
so the sequence $((\nu_t^{N_k})_{t\ge0}\dvtx k\in\mathbb{N})$ is Cauchy
in the
Skorohod space $D([0,\infty),(\mathcal{S},W))$,
and hence converges, with limit $(\nu_t)_{t\ge0}$ say, since
$D([0,\infty),\break  (\mathcal{S},W))$ is complete.

By Fatou's lemma and the moment estimate (\ref{NME}),
\[
\mathbb{E} \Bigl(\sup_{s\le t}\bigl\langle\llvert v\rrvert
^p,\nu_s\bigr\rangle \Bigr) \le\liminf
_k\mathbb{E} \Bigl(\sup_{s\le t}\bigl\langle
\llvert v\rrvert ^p,\nu ^{N_k}_s\bigr
\rangle1_{\{\langle
\llvert v\rrvert ^p,\nu_0^{N_k}\rangle\le C\langle\llvert v\rrvert ^p,\mu_0\rangle\}} \Bigr)<\infty,
\]
so $(\nu_t)_{t\ge0}$ is locally bounded in $\mathcal{S}(p)$ almost surely.
Fix a function $f$ on $\mathbb{R}^d$ satisfying $|f(v)|\le1$ and
$|f(v)-f(v')|\le|v-v'|$ for all $v,v'\in\mathbb{R}^d$.
From (\ref{UCP}), since $\|f\|\le2$, we see that $\langle f,\nu
^{N_k}_t\rangle
\to\langle f,\nu_t\rangle0$ uniformly on compact time intervals
almost surely.
Consider the equation
\[
\bigl\langle f,\nu^N_t\bigr\rangle=\bigl\langle f,
\nu_0^N\bigr\rangle+M^{N,f}_t+\int
_0^t\bigl\langle f,Q\bigl(\nu_s^N,
\nu _s^N\bigr)\bigr\rangle \,ds
\]
with $N=N_k$ in the limit $k\to\infty$.
Estimate (\ref{DME}) implies that $M^{N_k,f}_t\to0$ uniformly on
compact time intervals in probability.
Moreover,
\[
\bigl\langle f,Q\bigl(\nu_t^{N_k},\nu_t^{N_k}
\bigr)\bigr\rangle-\bigl\langle f,Q(\nu_t,\nu _t)\bigr
\rangle=\bigl\langle g_t,\nu ^{N_k}_t-
\nu_t\bigr\rangle,
\]
where
\[
g_t(v)=\int_{\mathbb{R}^d\times S^{d-1}}\bigl\{f\bigl(v'
\bigr)+f\bigl(v'_*\bigr)-f(v)-f(v_*)\bigr\} B(v-v_*,d\sigma) \bigl(
\nu_t^{N_k}+\nu_t\bigr) (dv_*)
\]
and, by some straightforward estimation, $\|g_t\|\le\max\{
16,12+8\kappa\}
$ for all $t\ge0$.
Hence, we can pass to the limit uniformly on compact time intervals in
probability to obtain
\[
\langle f,\nu_t\rangle=\langle f,\mu_0\rangle+\int
_0^t\bigl\langle f,Q(\nu_s,
\nu_s)\bigr\rangle \,ds
\]
for all $t\ge0$, almost surely.
A separability argument shows that almost surely, this equation holds
for all such functions $f$ and all $t\ge0$.
So, almost surely, $(\nu_t)_{t\ge0}$ is a solution, and in
particular, a locally bounded solution in $\mathcal{S}(p)$ exists.

Now let $(\mu_t)_{t\ge0}$ be any locally bounded solution in
$\mathcal{S}
(p)$ starting from $\mu_0$, and let $(\mu^N_t)_{t\ge0}$ be any Kac
process in $\mathcal{S}_N$.
Then
\[
\mu_t^N-\mu_t=\bigl(\mu_0^N-
\mu_0\bigr)+M^N_t+\int_0^t2Q
\bigl(\rho_s,\mu _s^N-\mu_s
\bigr)\,ds,
\]
where now $\rho_t=(\mu_t+\mu_t^N)/2$.
The argument of Section~\ref{BPR} applies without essential change to
show that, for all $t\ge0$ and all functions $f_t$ on $\mathbb{R}^d$,
we have
\[
\bigl\langle f_t,\mu_t^N-\mu_t
\bigr\rangle=\bigl\langle f_0,\mu_0^N-
\mu_0\bigr\rangle +\int_0^t\bigl
\langle f_s,dM_s^N\bigr\rangle,
\]
where $f_s(v)=E_{(s,v)}\langle f_t,\tilde\Lambda_t\rangle)$ and
where $(\Lambda
^*_t)_{t\ge s}$ is a linearized Kac process in environment $(\rho
_t)_{t\ge0}$.
Next, the argument of Section~\ref{PMR} applies to show that, for all
$\varepsilon\in(0,1]$ and all $T\in[0,\infty)$, for all $N\in
\mathbb{N}$, with
probability exceeding $1-\varepsilon$, for all $t\le T$, we have
\begin{equation}
\label{WNM} W\bigl(\mu^N_t,\mu_t\bigr)\le
C\bigl(W\bigl(\mu_0^N,\mu_0
\bigr)+N^{-\alpha(d,p)}\bigr),
\end{equation}
where $C<\infty$ depends only on $B,d,\varepsilon,{\lambda},p$ and $T$,
where ${\lambda}$ is an upper bound for $\langle|v|^p,\mu_0\rangle$
and $\langle
|v|^p,\mu_0^N\rangle$.
Convergence at rate $N^{-1/d}$ could be proved for $p>8$ by checking
that the arguments leading to the estimate for $W(\mu_t^N,\mu_t^{N'})$
apply also when $(\mu_t^{N'})_{t\ge0}$ is replaced by $(\mu_t)_{t\ge0}$.
Alternatively, we can find $N'$ so that $(N')^{-\alpha(d,p)}\le N^{-1/d}$
and, by Proposition~\ref{DAP},
$W(\nu_0^{N'},\mu_0)\le CN^{-1/d}$ with probability exceeding
$1-\varepsilon$.
Then, by Theorem~\ref{MR} and (\ref{WNM}), with probability exceeding
$1-3\varepsilon$, for all $t\le T$, we have
\begin{eqnarray*}
W\bigl(\mu^N_t,\mu_t\bigr)&\le & W\bigl(
\mu^N_t,\nu_t^{N'}\bigr)+ W\bigl(
\nu^{N'}_t,\mu _t\bigr)\\
&\le & C\bigl(W\bigl(
\mu_0^N,\nu_0^{N'}\bigr)+W\bigl(
\nu_0^{N'},\mu _0\bigr)+N^{-1/d}+
\bigl(N'\bigr)^{-\alpha(d,p)}\bigr)
\\
&\le& C\bigl(W\bigl(\mu_0^N,\mu_0
\bigr)+4N^{-1/d}\bigr).
\end{eqnarray*}
Finally, we can take $\mu_t^{N_k}=\nu_t^{N_k}$ and let $k\to\infty$
to see that $\mu_t=\nu_t$ for all $t\ge0$, so $(\nu_t)_{t\ge0}$
is the only solution which is locally bounded in $\mathcal{S}(p)$.
\end{pf}

We can combine Theorem~\ref{CBE} with Proposition~\ref{DAP} to obtain
the following stochastic approximation for solutions to the spatially
homogeneous Boltzmann equation.
\begin{corollary}\label{COR}
Assume that the collision kernel $B$ satisfies conditions (\ref
{LBP}) and (\ref{LCP}).
Let $\mu_0\in\mathcal{S}(p)$ for some $p\in(2,\infty)$, and let
$(\mu
_t)_{t\ge0}$ be the unique locally bounded solution to (\ref
{WB}) in $\mathcal{S}(p)$ starting from $\mu_0$.
Write $\mu_0^N$ for the random variable in $\mathcal{S}_N$
constructed by
sampling from $\mu_0$ as in Proposition~\ref{DAP}, and conditioning
on $\mu_0^N$, let $(\mu_t^N)_{t\ge0}$ be a Kac process starting from
$\mu_0^N$.
Then, for all $\varepsilon\in(0,1]$, all ${\lambda}\ge\langle
|v|^p,\mu_0\rangle$ and
all $T\in[0,\infty)$, there are constants $\alpha(d,p)>0$ and
$C(B,d,\varepsilon
,{\lambda},p,T)<\infty$, such that
with probability exceeding $1-\varepsilon$, for all $t\le T$,
\[
W\bigl(\mu^N_t,\mu_t\bigr)\le
CN^{-\alpha}.
\]
For $p>8$, we can take $\alpha=1/d$ when $d\ge3$, and the estimate holds
with $N^{-1/2}\log N$ in place of $N^{-\alpha}$ when $d=2$.
\end{corollary}

On the other hand, if one views the spatially homogeneous Boltzmann
equation as a means to compute approximations to the Kac process,
the following estimate provides a measure of accuracy for this procedure.
\begin{corollary}\label{COR2}
Assume that the collision kernel $B$ satisfies conditions (\ref
{LBP}) and (\ref{LCP}).
Fix $d\ge3$, $\varepsilon\in(0,1]$ and $\tau,T\in(0,\infty)$ with
$\tau\le T$.
There is a constant $C<\infty$, depending only on $B$, $d$,
$\varepsilon$,
$\tau$ and $T$, with the following property.
Let $N\in\mathbb{N}$ and let $(\mu_t^N)_{t\ge0}$ be a Kac process
in $\mathcal{S}
_N$ with collision kernel $B$.
Denote by $(\mu_t)_{t\ge\tau}$ the solution to the spatially
homogeneous Boltzmann equation with collision kernel $B$ starting from
$\mu_\tau^N$ at time $\tau$.
Then, with probability exceeding $1-\varepsilon$, for all $t\in[\tau
,T]$, we
have $W(\mu_t^N,\mu_t)\le CN^{-1/d}$.
The same holds for $d=2$ if we replace $N^{-1/d}$ by $N^{-1/2}\log N$.
\end{corollary}
\begin{pf}
Use (\ref{NMG}) to find a constant ${\lambda}(B,\tau,\varepsilon
)<\infty$
such that $\langle|v|^9,\mu_\tau^N\rangle\le{\lambda}$ with probability
exceeding $1-\varepsilon/2$.
Then apply Theorem~\ref{CBE} with $\varepsilon/2$ in place of
$\varepsilon$ to find
the desired constant $C$.
\end{pf}
\begin{appendix}
\section*{Appendix}\label{APP}
We state and prove a basic lemma on the time-evolution of signed measures,
which allows us to control the evolution of the total variation when
the signed measures are given by an integral over time.
Let $(E,\mathcal{E})$ be a measurable space.
Write $\mathcal{M}^+$ (resp., $\mathcal{M}$) for the set of finite
measures (resp.,
signed measures of finite total variation) on $(E,\mathcal{E})$.
For $\mu\in\mathcal{M}$, write $|\mu|$ for the associated total variation
measure and $\|\mu\|$ for the total variation.

\setcounter{lemma}{0}
\begin{lemma}\label{BML}
Assume that $(E,\mathcal{E})$ is separable.
Let $T\in(0,\infty)$.
Let $\mu_0\in\mathcal{M}$ and ${\lambda}_0\in\mathcal{M}^+$ be
given, along with
a measurable map $t\mapsto\nu_t\dvtx [0,T]\to\mathcal{M}$
such that $\nu_t$ is absolutely continuous with respect to ${\lambda
}_0$ for all $t\in[0,T]$ and $ \int_0^T\|\nu_t\|\,dt<\infty$.
Set
\[
\mu_t=\mu_0+\int_0^t
\nu_s\,ds.
\]
Then there exists a measurable map $\sigma\dvtx [0,T]\times E\to\{-1,0,1\}$
such that, for all $t\in[0,T]$, we have $\mu_t=\sigma_t|\mu_t|$ and
\[
\llvert \mu_t\rrvert =\llvert \mu_0\rrvert +\int
_0^t\sigma_s\nu_s\,ds.
\]
\end{lemma}
A version of the lemma, without the hypothesis of separability and for
the case where $t\mapsto\nu_t\dvtx [0,T]\to\mathcal{M}$ is continuous in total
variation, was stated by Kolokoltsov in \cite{MR2223422}, Lemma~A.1.
The proof given in \cite{MR2223422} contains a gap, which we have not
been able to fill.
The case where $(E,\mathcal{E})$ is $\mathbb{R}^d$ with its Borel
$\sigma$-algebra and
where $t\mapsto\nu_t\dvtx [0,T]\to\mathcal{M}$ is continuous in total variation,
has been proved by Lu and Mouhot \cite{MR2871802}, Lemma~5.1.
We will use a substantially different argument, which allows us to
replace this hypothesis of continuity with the existence of a reference
measure ${\lambda}_0$.
\begin{pf*}{Proof of Lemma~\ref{BML}}
There exists an increasing sequence $(\mathcal{E}_n\dvtx n\in\mathbb{N})$
of finite $\sigma$-algebras generating $\mathcal{E}$.
Write $\mathcal{A}_n$ for the partition of $E$ generating $\mathcal{E}_n$.
Consider the finite measure ${\lambda}={\lambda}_0+|\mu_0|+\int_0^T|\nu_t|\,dt$ on $(E,\mathcal{E})$.
By scaling we reduce to the case where ${\lambda}$ is a probability measure.
For each $t\in[0,T]$, define $\mathcal{E}_n$-measurable functions
$\alpha_t^n$
and $\beta_t^n$ by on $E$ by setting
\[
\alpha_t^n(x)=\mu_t(A)/{\lambda}(A),\qquad
\beta_t^n(x)=\nu _t(A)/{\lambda}(A)
\]
if $x\in A$ for some $A\in\mathcal{A}_n$ with ${\lambda}(A)>0$ and setting
$\alpha_t^n(x)=\beta_t^n(x)=0$ if there is no such $A$.
Then, for all $x\in E$, the map $t\mapsto\beta_t^n(x)$ is integrable on
$[0,T]$ and
\[
\alpha_t^n(x)=\alpha_0^n(x)+
\int_0^t\beta_s^n(x)\,ds.
\]
For each $t\in[0,T]$, we have $|\mu_t|\le{\lambda}$ so $|\alpha
^n_t|\le1$ and $\alpha^n_t{\lambda}=\mu_t$ on $\mathcal{E}_n$.
Moreover, the sequence $(\alpha_t^n\dvtx n\in\mathbb{N})$ is a ${\lambda
}$-martingale in the filtration $(\mathcal{E}_n\dvtx n\in\mathbb{N})$.
So, by the martingale convergence theorem, there exists $\tilde\alpha
_t\in L^1({\lambda})$
such that $\alpha_t^n\to\tilde\alpha_t$ as $n\to\infty$,
${\lambda
}$-almost everywhere and in $L^1({\lambda})$.
Then $\tilde\alpha_t{\lambda}=\mu_t$ on $\bigcup_n\mathcal{E}_n$ and
hence on
$\mathcal{E}$ by uniqueness of extension.

For $\tau=\{t_0,\dots,t_N\}\subseteq[0,T]$ with $t_0<\cdots<t_N$ and any
function $(\alpha_t(x)\dvtx t\in[0,T],x\in E)$,
define a function $|\alpha|_\tau$ on $E$ by
\[
\llvert \alpha\rrvert _\tau=\llvert \alpha_0\rrvert +\sum
_{k=0}^{N-1}\llvert \alpha_{t_{k+1}}-
\alpha_{t_k}\rrvert.
\]
Then, for all $A\in\mathcal{A}_n$, on $A$, we have
\[
{\lambda}(A)\bigl\llvert \alpha^n\bigr\rrvert _\tau=
\llvert \mu_0\rrvert (A)+\sum_{k=0}^{N-1}
\llvert \mu _{t_{k+1}}-\mu_{t_k}\rrvert (A)\le{\lambda}(A)
\]
so $|\alpha^n|_\tau\le1$ everywhere.

Fix $E_0\in\mathcal{E}$ with ${\lambda}(E_0)=1$ such that
$\alpha^n_t(x)\to\tilde\alpha_t(x)$ as $n\to\infty$ for all $t\in
[0,T]\cap(T\mathbb{Q})$ and all $x\in E_0$.
Write $\mathcal{T}$ for the set of finite subsets of $[0,T]\cap
(T\mathbb{Q})$.
Then, for all $x\in E_0$ and $\tau\in\mathcal{T}$, we have $|\tilde
\alpha|_\tau
(x)\le1$,
so the map $t\mapsto\tilde\alpha_t(x)\dvtx [0,T]\cap(T\mathbb{Q})\to
[-1,1]$ has
total variation bounded by $1$.
Hence, for $x\in E_0$, we can define a c\`adl\`ag map $t\mapsto\alpha
_t(x)\dvtx [0,T]\to[-1,1]$ by
\[
\alpha_t(x)= %
\cases{\displaystyle\lim_{s\to t,s\in(t,T)\cap(T\mathbb{Q})}\tilde
\alpha _s(x),& \quad$t\in [0,T)$,\vspace*{3pt}
\cr
\tilde
\alpha_T(x),& \quad$t=T$.} %
\]
For $x\in E\setminus E_0$, set $\alpha_t(x)=0$ for all $t\in[0,T]$.
We have $\alpha_T{\lambda}=\tilde\alpha_T{\lambda}=\mu_T$ as we
showed above.
For $t\in[0,T)$ and $s\in(t,T)\cap(T\mathbb{Q})$, we have in the limit
$s\to t$
\[
\llVert \alpha_t{\lambda}-\mu_t\rrVert \le\llVert
\alpha_t{\lambda}-\tilde\alpha _s{\lambda}\rrVert +
\llVert \mu_s-\mu_t\rrVert \le\bigl\langle\llvert
\alpha_t-\tilde\alpha_s\rrvert ,{\lambda}\bigr\rangle+
\int_t^s\llVert \nu_r\rrVert \,dr
\to0
\]
so $\alpha_t{\lambda}=\mu_t$. Define $\sigma\dvtx [0,T]\times E\to\{
-1,0,1\}$
by $\sigma_t(x)=\operatorname{sgn}(\alpha_t(x))$. Then $\sigma$ is
measurable and
$\mu_t=\sigma_t\llvert \mu_t\rrvert $ for all $t\in[0,T]$.

For any function $\psi$ on $[-1,1]$ with continuous bounded
derivative, we have
\[
\psi\bigl(\alpha_t^n(x)\bigr)=\psi\bigl(
\alpha_0^n(x)\bigr)+\int_0^t
\psi'\bigl(\alpha _s^n(x)\bigr)
\beta^n_s(x)\,ds
\]
for all $t\in[0,T]$ and all $x\in E$.
Since $\nu_t$ is absolutely continuous with respect to ${\lambda}$
for all $t\in[0,T]$, we have on $\mathcal{E}_n$
\[
\psi\bigl(\alpha_t^n\bigr){\lambda}=\psi\bigl(
\alpha_0^n\bigr){\lambda}+\int_0^t
\psi '\bigl(\alpha _s^n\bigr)
\nu_s\,ds
\]
for all $t\in[0,T]$.
Since $\nu_s(dx)\,ds$ is absolutely continuous with respect to ${\lambda
}(dx)\,ds$, we have $\alpha_s^n(x)\to\alpha_s(x)$ as $n\to\infty$ almost
everywhere for $\nu_s(dx)\,ds$.
Hence, on letting $n\to\infty$, we obtain on $\bigcup_n\mathcal{E}_n$
\[
\psi(\alpha_t){\lambda}=\psi(\alpha_0){\lambda}+\int
_0^t\psi '(\alpha_s)
\nu_s\,ds
\]
for all $t\in[0,T]$. The identity then holds on $\mathcal{E}$ by uniqueness
of extension.
Set $\psi_k(x)=\sqrt{x^2+1/k}$.
Then\vspace*{1pt} $\psi_k(x)\to|x|$ and $\psi_k'(x)\to\operatorname{sgn}(x)$ as
$k\to\infty
$ for all $x\in[-1,1]$.
By dominated convergence, for all $A\in\mathcal{E}$ and all $t\in
[0,T]$, we have
\[
\bigl\langle\psi_k(\alpha_t)1_A,{\lambda}
\bigr\rangle\to\bigl\langle\llvert \alpha _t\rrvert
1_A,{\lambda}\bigr\rangle=\llvert \mu_t\rrvert (A)
\]
and
\[
\int_0^t\bigl\langle\psi_k'(
\alpha_s)1_A,\nu_s\bigr\rangle \,ds\to\int
_0^t\bigl\langle\operatorname{sgn}(\alpha
_s)1_A,\nu_s\bigr\rangle \,ds=\int
_0^t\langle\sigma_s1_A,
\nu_s\rangle \,ds.
\]
Hence, on taking $\psi=\psi_k$ above and letting $k\to\infty$, we
obtain the desired identity.
\end{pf*}
\end{appendix}

\section*{Acknowledgments}

I am grateful to Cl\'{e}ment Mouhot and Richard Nickl for several
discussions in the course of this work, and to a helpful referee whose
comments led to improvements in the paper.


%



\printaddresses
\end{document}